\documentclass{amsart}
\def\Defined{}
\ifx\UseRussian\Defined
	\usepackage[T2A,T2B]{fontenc}
	\usepackage[cp1251]{inputenc}
	\usepackage[english,russian]{babel}
\fi
\usepackage{footmisc}
\usepackage[all]{xy}
\usepackage{color}
\definecolor{UrlColor}{rgb}{.9,0,.3}
\definecolor{SymbColor}{rgb}{.4,0,.9}
\definecolor{IndexColor}{rgb}{1,.3,.6}
\definecolor{eml1}{rgb}{.8,.1,.1}
\definecolor{eml2}{rgb}{.1,.6,.6}

\usepackage{xr-hyper}
\usepackage[unicode]{hyperref}
\hypersetup{pdfdisplaydoctitle=true}
\hypersetup{colorlinks}
\hypersetup{citecolor=UrlColor}
\hypersetup{urlcolor=UrlColor}
\hypersetup{pdffitwindow=true}
\hypersetup{pdfnewwindow=true}
\hypersetup{pdfstartview={FitH}}

\def\hyph{\penalty0\hskip0pt\relax-\penalty0\hskip0pt\relax}
\def\Hyph{-\penalty0\hskip0pt\relax}

\newcommand{\Basis}[1]{\overline{\overline{#1}}{}}
\newcommand{\Vector}[1]{\overline{#1}{}}
\newcommand{\gi}[1]{\boldsymbol{\textcolor{IndexColor}{#1}}}
\makeatletter
\newcommand{\NameDef}[1]{%
	\expandafter\gdef\csname #1\endcsname%
}%
\newcommand{\ShowSymbol}[1]{%
	\@nameuse{ViewSymbol#1}%
}%
\newcommand{\symb}[3]{%
	\@ifundefined{ViewSymbol#3}{%
		\NameDef{ViewSymbol#3}{\textcolor{SymbColor}{#1}}%
		\NameDef{RefSymbol#3}{\pageref{symbol: #3}}%
		\@namedef{LabelSymbol#3}{\label{symbol: #3}}%
	}{%
		\NameDef{RefSymbol#3}{}%
		\@namedef{LabelSymbol#3}{}%
	}%
	\ifcase#2%0
	\or%1
		$\@nameuse{ViewSymbol#3}$%
	\or%2
		\[\@nameuse{ViewSymbol#3}\]%
	\else%
	\fi%
	\@nameuse{LabelSymbol#3}%
}%
\makeatother

\newcommand{\subs}{${}_*$\Hyph}
\newcommand{\sups}{${}^*$\Hyph}

\newcommand{\CRstar}{{}_*{}^*}
\newcommand{\RCstar}{{}^*{}_*}

\newcommand{\RC}{$\RCstar$\Hyph}
\newcommand{\CR}{$\CRstar$\Hyph}
\newcommand{\drc}{$D\RCstar$\Hyph}
\newcommand{\Drc}{$\mathcal D\RCstar$\Hyph}
\newcommand{\dcr}{$D\CRstar$\hyph}
\newcommand{\rcd}{$\RCstar D$\Hyph}
\newcommand{\crd}{$\CRstar D$\Hyph}

\newcommand\sT{$\star T$\Hyph}%
\newcommand\Ts{$T\star$\Hyph}%

\renewcommand{\uppercasenonmath}[1]{}

\makeatletter
\newcommand\@dotsep{4.5}
\def\@tocline#1#2#3#4#5#6#7
{\relax
		\par \addpenalty\@secpenalty\addvspace{#2}%
		\begingroup \hyphenpenalty\@M
		\@ifempty{#4}{%
			\@tempdima\csname r@tocindent\number#1\endcsname\relax
		}{%
			\@tempdima#4\relax
		}%
		\parindent\z@ \leftskip#3\relax \advance\leftskip\@tempdima\relax
		\rightskip\@pnumwidth plus1em \parfillskip-\@pnumwidth
		#5\leavevmode\hskip-\@tempdima #6\relax
		\leaders\hbox{$\m@th
		\mkern \@dotsep mu\hbox{.}\mkern \@dotsep mu$}\hfill
		\hbox to\@pnumwidth{\@tocpagenum{#7}}\par
		\nobreak
		\endgroup
}
\makeatother 

\ifx\PrintBook\undefined
	\usepackage{fancyhdr}
	\makeatletter
	\renewcommand{\@indextitlestyle}{%
		\twocolumn[\section{\indexname}]%
		\def\IndexSpace{off}%
	}
	\makeatother 
	\thanks{\href{mailto:Aleks\_Kleyn@MailAPS.org}{Aleks\_Kleyn@MailAPS.org}}
\else
	\makeatletter
	\renewcommand{\@indextitlestyle}{%
		\twocolumn[\chapter{\indexname}]%
		\def\IndexSpace{off}%
		\let\@secnumber\@empty
		\chaptermark{\indexname}%
	}
	\makeatother 
	\email{\href{mailto:Aleks\_Kleyn@MailAPS.org}{Aleks\_Kleyn@MailAPS.org}}
\fi

\ifx\SelectlEnglish\undefined
	\ifx\UseRussian\undefined
		\def\SelectlEnglish{}
	\fi
\fi

\ifx\SelectlEnglish\undefined
	\newcommand\CurrentLanguage{Russian.}%
	\author{Александр Клейн}
	\newtheorem{theorem}{Теорема}[section]
	\newtheorem{corollary}[theorem]{Следствие}
	\theoremstyle{definition}
	\newtheorem{definition}[theorem]{Определение}
	\newtheorem{example}[theorem]{Пример}
	\newtheorem{xca}[theorem]{Exercise}
	\theoremstyle{remark}
	\newtheorem{remark}[theorem]{Замечание}
	
	\newcommand\Gbasis{$G$\Hyph базис}
	\newcommand\Gcoords{$G$\Hyph координат}
	\newcommand\Gspace{$G$\Hyph пространств}
	\newcommand\xRefDef[2]
		{
		\externaldocument[#1-Russian-]{#1.Russian}[http://arxiv.org/PS_cache/#2.pdf]
		\NameDef{xRefDef#1}{}%
		}
	\makeatletter
	\newcommand\xRef[2]%
	{%
		\@ifundefined{xRefDef#1}{%
		\ref{#2}%
		}{% 
		\citeBib{#1}-\ref{#1-Russian-#2}%
		}%
	}%
	\newcommand\xEqRef[2]%
	{%
		\@ifundefined{xRefDef#1}{%
		\eqref{#2}%
		}{%
		\citeBib{#1}-\eqref{#1-Russian-#2}%
		}%
	}%
	\makeatother
	\ifx\PrintBook\undefined
		\newcommand{\BibTitle}{%
			\section{Список литературы}%
		}
	\else
		\newcommand{\BibTitle}{%
			\chapter{Список литературы}%
		}
	\fi
\else
	\newcommand\CurrentLanguage{English.}%
	\author{Aleks Kleyn}
	\newtheorem{theorem}{Theorem}[section]
	
	\theoremstyle{definition}
	\newtheorem{definition}[theorem]{Definition}
	\newtheorem{example}[theorem]{Example}
	
	\theoremstyle{remark}
	\newtheorem{remark}[theorem]{Remark}
	\newcommand\Gbasis{$G$\Hyph basis}
	\newcommand\Gcoords{$G$\Hyph coordinates}
	\newcommand\Gspace{$G$\Hyph space}
	\newcommand\xRefDef[2]
		{
		\externaldocument[#1-English-]{#1.English}[http://arxiv.org/PS_cache/#2.pdf]
		\NameDef{xRefDef#1}{}%
		}
	\makeatletter
	\newcommand\xRef[2]%
	{%
		\@ifundefined{xRefDef#1}{%
		\ref{#2}%
		}{% 
		\citeBib{#1}-\ref{#1-English-#2}%
		}%
	}%
	\newcommand\xEqRef[2]%
	{%
		\@ifundefined{xRefDef#1}{%
		\eqref{#2}%
		}{%
		\citeBib{#1}-\eqref{#1-English-#2}%
		}%
	}%
	\makeatother
	\ifx\PrintBook\undefined
		\newcommand{\BibTitle}{%
			\section{References}%
		}
	\else
		\newcommand{\BibTitle}{%
			\chapter{References}%
		}
	\fi
\fi

\ifx\PrintBook\undefined
	\numberwithin{Hfootnote}{section}
\else
	\numberwithin{section}{chapter}
	\numberwithin{footnote}{chapter}
	\numberwithin{Hfootnote}{chapter}
\fi

\numberwithin{equation}{section}
\numberwithin{figure}{section}
\numberwithin{table}{section}
\numberwithin{Item}{section}

\makeatletter
\newcommand\org@maketitle{}
\let\org@maketitle\maketitle
\def\maketitle{%
	\hypersetup{pdftitle={\@title}}%
	\hypersetup{pdfauthor={\authors}}%
	\hypersetup{pdfsubject=\@keywords}%
	\org@maketitle
}
\def\make@stripped@name#1{%
	\begingroup
		\escapechar\m@ne
		\global\let\newname\@empty
		\protected@edef\Hy@tempa{\CurrentLanguage #1}%
		\edef\@tempb{%
			\noexpand\@tfor\noexpand\Hy@tempa:=%
			\expandafter\strip@prefix\meaning\Hy@tempa
		}%
		\@tempb\do{%
			\if\Hy@tempa\else
				\if\Hy@tempa\else
					\xdef\newname{\newname\Hy@tempa}%
				\fi
			\fi
		}%
	\endgroup
}%
\newenvironment{enumBib}{%
	\BibTitle
	\advance\@enumdepth \@ne
	\edef\@enumctr{enum\romannumeral\the\@enumdepth}\list
	{\csname biblabel\@enumctr\endcsname}{\usecounter
	{\@enumctr}\def\makelabel##1{\hss\llap{\upshape##1}}}
}{%
	\endlist
}
\def\Chapters#1{\ChapterList#1,LastChapter,}%
\def\LastChapter{LastChapter}%
\def\ChapterList#1,{\def\temp{#1}%
	\ifx\temp\LastChapter
	\else
		\@ifundefined{#1}{%
		}{% 
			\def\Semafor{on}
		}
		\expandafter\ChapterList
	\fi
}%
\newcommand{\BiblioItem}[3]
{
	\def\Semafor{off}
	\Chapters{#1}
	\ifx\Semafor\ValueOn
		\ifx\IndexState\ValueOff
			\begin{enumBib}
			\def\IndexState{on}
		\fi
		\item \label{bibitem: #2}#3%
	\fi
}
\newcommand{\OpenBiblio}
{
	\def\IndexState{off}
}
\newcommand{\CloseBiblio}
{
	\ifx\IndexState\ValueOn
		\end{enumBib}
		\def\IndexState{off}
	\fi
}

\def\StartCite{[}%
\def\citeBib#1{\protect\showCiteBib#1,endCite,}%
\def\endCite{endCite}%
\def\showCiteBib#1,{\def\temp{#1}%
\ifx\temp\endCite
]%
\def\StartCite{[}%
\else
	\StartCite\ref{bibitem: #1}%
	\def\StartCite{, }%
\expandafter\showCiteBib%
\fi}%

\makeatother 
\newcommand{\arp}{\ar @{-->}}
\newcommand{\ars}{\ar @{.>}}
\newcommand{\bundle}[4]%
{%
	\def\tempa{}%
	\def\tempb{#3}%
	\def\tempc{#1}%
	\ifx\tempa\tempb%
		\ifx\tempa\tempc%
			#2%
		\else%
			\xymatrix{#2:#1\arp[r]&#4}%
		\fi%
	\else%
		\ifx\tempa\tempc%
			#2[#3]%
		\else%
			\xymatrix{#2[#3]:#1\arp[r]&#4}%
		\fi%
	\fi%
}%
\newcommand{\AddIndex}[2]%
{%
	{\bf #1}%
	\label{index: #2}%
}%
\newcommand{\Index}[3]%
{%
	\def\Semafor{off}%
	\Chapters{#1}%
	\ifx\Semafor\ValueOn%
		\def\tempa{}%
		\def\tempb{#3}%
		\ifx\IndexState\ValueOff%
			\begin{theindex}%
			\def\IndexState{on}%
		\fi%
		\ifx\IndexSpace\ValueOn%
			\indexspace%
			\def\IndexSpace{off}%
		\fi%
		\item #2%
		\ifx\tempa\tempb%
		\else%
			\ \pageref{index: #3}%
		\fi%
	\fi%
}%
\newcommand{\SubIndex}[3]
{
	\def\Semafor{off}
	\Chapters{#1}
	\ifx\Semafor\ValueOn
		\subitem #2 \pageref{index: #3}%\RefPage{#3}
	\fi
}%
\makeatletter
\newcommand{\Symb}[3]
{
	\def\Semafor{off}
	\Chapters{#1}
	\ifx\Semafor\ValueOn
		\ifx\IndexState\ValueOff
			\begin{theindex}
			\def\IndexState{on}
		\fi
		\ifx\IndexSpace\ValueOn
			\indexspace
			\def\IndexSpace{off}
		\fi
		\item $\@nameuse{ViewSymbol#3}$\ \ #2
		\@nameuse{RefSymbol#3}%
	\fi
}

\makeatother
\newcommand{\SetIndexSpace}%
{%
	\def\IndexSpace{on}%
}%
\def\ValueOff{off}
\def\ValueOn{on}

\newcommand{\OpenIndex}
{
	\def\IndexState{off}
}
\newcommand{\CloseIndex}
{
	\ifx\IndexState\ValueOn
		\end{theindex}
		\def\IndexState{off}
	\fi
}

\def\LastMemo{LastMemo}%
\def\MemoList#1//{\def\temp{#1}%
	\ifx\temp\LastMemo
	\else%
		\par\setlength{\parindent}{12pt}\textcolor{blue}{#1}%
		\expandafter\MemoList%
	\fi%
}%

\xRefDef{0701.238}{math/pdf/0701/0701238}
\xRefDef{0707.2246}{arxiv/pdf/0707/0707.2246v1}
\begin{document}
%auto-ignore
\def\texTstarMorphism{}
\ifx\PrintBook\undefined
\title{Morphism of \texorpdfstring{\Ts}{T*-}Representations}
\pdfbookmark[1]{Morphism of T*-Representations}{TitleEnglish}
\begin{abstract}
Importance of theorem dedicated to isomorphisms consist in statement
that they allow to identify different mathematical
objects which have something common from the point of view of certain model.
This paper cosiders morphisms of \Ts representation of $\mathfrak{F}$\Hyph algebra
and morphisms of \Ts representation of fibered $\mathfrak{F}$\Hyph algebra.
\end{abstract}

\maketitle

This paper appeared on intersection of two researches which I
make at the same time.
First half of paper is dedicated 
to morphisms of \Ts representations of $\mathfrak{F}$\Hyph algebra.
In second half I consider morphisms
of \Ts representations of fibered $\mathfrak{F}$\Hyph algebra.
Considered constructions appeared as result of study of
\drc linear mappings which are
morphisms of \Ts representations of skew field in Abelian group.
Therefore I use \drc linear mappings
for the purposes of illustration of stated theory.
\fi

			\section{Representation of \texorpdfstring{$\mathfrak{F}$\Hyph }{F-}Algebra}
			\label{section: Representation of F-Algebra}

	\begin{definition}
We call the map
\[
t:M\rightarrow M
\]
\AddIndex{transformation of set}{transformation of set} $M$.
	 \qed
	 \end{definition}

	 \begin{definition}
Transformations is
\AddIndex{left-side transformation}{left-side transformation} or
\AddIndex{\Ts transformation}{Tstar transformation}
if it acts from left
\[
u'=t u
\]
We denote
\symb{{}^\star M}1{set of Tstar transformations}
the set of \Ts transformations of set $M$.

Suppose we defined the structure of $\mathfrak{H}$\Hyph algebra on the set $M$
(\citeBib{Burris Sankappanavar}).
Then the set ${}^\star M$ consists from \Ts transformations which are homomorphisms
of $\mathfrak{H}$\Hyph algebra.
	 \qed
	 \end{definition}

	 \begin{definition}
Transformations is
\AddIndex{right-side transformations}{right-side transformation} or
\AddIndex{\sT transformation}{starT transformation}
if it acts from right
\[
u'= ut
\]
We denote
\symb{M^\star}1{set of starT transformations}
the set of nonsingular \sT transformations of set $M$.

Suppose we defined the structure of $\mathfrak{H}$\Hyph algebra on the set $M$
(\citeBib{Burris Sankappanavar}).
Then the set $M^\star$ consists from \sT transformations which are homomorphisms
of $\mathfrak{H}$\Hyph algebra.
	 \qed
	 \end{definition}

We denote
\symb{\delta}1{identical transformation}
identical transformation.

		\begin{definition}
		\label{definition: Tstar representation of algebra} 
Suppose we defined the structure of $\mathfrak{F}$\Hyph algebra on the set ${}^\star M$
(\citeBib{Burris Sankappanavar}).
Let $A$ be $\mathfrak{F}$\Hyph algebra.
We call homomorphism
	\begin{equation}
f:A\rightarrow {}^\star M
	\label{eq: Tstar representation of algebra}
	\end{equation}
\AddIndex{left-side}
{left-side representation of algebra} or
\AddIndex{\Ts representation of $\mathfrak{F}$\Hyph algebra $A$ in set $M$}
{Tstar representation of algebra}
		\qed
		\end{definition}

		\begin{definition}
		\label{definition: starT representation of algebra} 
Suppose we defined the structure of $\mathfrak{F}$\Hyph algebra on the set $M^\star$
(\citeBib{Burris Sankappanavar}).
Let $A$ be $\mathfrak{F}$\Hyph algebra.
We call homomorphism
	\[
f:A\rightarrow M^\star
	\]
\AddIndex{right-side}
{right-side representation of algebra} or
\AddIndex{\sT representation of $\mathfrak{F}$\Hyph algebra $A$ in set $M$}
{starT representation of algebra}
		\qed
		\end{definition}

We extend to representation theory convention
described in remark
\xRef{0701.238}{remark: left and right matrix notation}.
We can write duality principle in the following form

		\begin{theorem}[duality principle]
		\label{theorem: duality principle, algebra representation}
Any statement which holds for \Ts representation of $\mathfrak{F}$\Hyph algebra $A$
holds also for \sT representation of $\mathfrak{F}$\Hyph algebra $A$.
		\end{theorem}

Diagram
\[
\xymatrix{
M\ar[rr]^{f(a)}&&M\\
&A\ar@{=>}[u]^f&
}
\]
means that we consider the representation of $\mathfrak{F}$\Hyph algebra $A$.
The map $f(a)$ is image of $a\in A$.

	 \begin{definition}
	 \label{definition: effective representation of algebra}
Suppose map \eqref{eq: Tstar representation of algebra} is
an isomorphism of the $\mathfrak{F}$\Hyph algebra $A$ into ${}^\star M$.
Then the \Ts representation of the $\mathfrak{F}$\Hyph algebra $A$ is called
\AddIndex{effective}{effective representation of algebra}.
	 \qed
	 \end{definition}

		\begin{remark}
		\label{remark: notation for effective representation of algebra}
Suppose the \Ts representation of $\mathfrak{F}$\Hyph algebra is effective. Then we identify
an element of $\mathfrak{F}$\Hyph algebra and its image and write \Ts transormation
caused by element $a\in A$
as
\[v'=av\]
Suppose the \sT representation of $\mathfrak{F}$\Hyph algebra is effective. Then we identify
an element of $\mathfrak{F}$\Hyph algebra and its image and write \sT transormation
caused by element $a\in A$
as
\[v'=va\]
		\qed
		\end{remark}

	 \begin{definition}
	 \label{definition: transitive representation of algebra}
We call a \Ts representation of $\mathfrak{F}$\Hyph algebra
\AddIndex{transitive}{transitive representation of algebra}
if for any $a, b \in V$ exists such $g$ that
\[a=f(g)b\]
We call a \Ts representation of $\mathfrak{F}$\Hyph algebra
\AddIndex{single transitive}{single transitive representation of algebra}
if it is transitive and effective.
	 \qed
	 \end{definition}

		\begin{theorem}	%single transitive
\Ts representation is single transitive if and only if for any $a, b \in M$
exists one and only one $g\in A$ such that $a=f(g)b$
		\end{theorem}
		\begin{proof}
Corollary of definitions \ref{definition: effective representation of algebra}
and \ref{definition: transitive representation of algebra}.
		\end{proof}

Suppose we introduce additional structure on set $M$.
Then we create an additional requirement for the representation of $\mathfrak{F}$\Hyph algebra.

Since we defined the structure of algebra of type $\mathfrak{H}$ on the set $M$,
we suppose that \Ts transformation
	\[
u'=f(a)u
	\]
is automorphism of algebra of type $\mathfrak{H}$.
We also study \Ts representations which reflects symmetry of algebra of type $\mathfrak{H}$.

Since we introduce continuity on set $M$,
we suppose that \Ts transformation
	\[
u'=f(a)u
	\]
is continuous in $u$. 
Therefore, we get
\[\left|\frac { \partial u'} {\partial u}\right|\neq 0\]

\section{Morphism of \texorpdfstring{$T\star$}{T*}-Representations
of \texorpdfstring{$\mathfrak{F}$\Hyph }{F-}Algebra}

		\begin{theorem}
Let $A$ and $B$ be $\mathfrak{F}$\Hyph algebras.
\Ts representation of $\mathfrak{F}$\Hyph algebra $B$
	\[
f:B\rightarrow {}^\star M
	\]
and homomorphism of $\mathfrak{F}$\Hyph algebra
	\begin{equation}
h:A\rightarrow B
	\label{eq: homomorphism of algebra F}
	\end{equation}
define \Ts representation of $\mathfrak{F}$\Hyph algebra $A$
	\[
\xymatrix{
A\ar[dr]^h\ar[rr]^f&&{}^\star M\\
&B\ar[ur]^g&
}
	\]
		\end{theorem}
		\begin{proof}
Since mapping $f$ is homomorphism of $\mathfrak{F}$\Hyph algebra $B$ into
$\mathfrak{F}$\Hyph algebra ${}^\star M$,
the mapping $h$ is homomorphism of $\mathfrak{F}$\Hyph algebra $A$ into
$\mathfrak{F}$\Hyph algebra ${}^\star M$.
		\end{proof}

Considering representations of $\mathfrak{F}$\Hyph algebra in sets $M$ and $N$,
we are interested in a mapping that preserves the structure of representation.

		\begin{definition}
		\label{definition: morphism of representations of F algebra}
Let us consider \Ts representation
\[
f:A\rightarrow {}^\star M
\]
of $\mathfrak{F}$\Hyph algebra $A$ in $M$
and \Ts representation
\[
g:B\rightarrow {}^\star N
\]
of $\mathfrak{F}$\Hyph algebra $B$ in $N$.
Tuple of mapping $(r,R)$
	\begin{equation}
\xymatrix{
r:A\ar[r]&B&R:M\ar[r]&N
}
	\label{eq: morphism of representations of F algebra, definition, 1}
	\end{equation}
such, that $r$ is homomorphism of $\mathfrak{F}$\Hyph algebra and
	\begin{equation}
\textcolor{blue}{R(f(a)m)}=g(\textcolor{red}{r(a)})\textcolor{blue}{R(m)}
	\label{eq: morphism of representations of F algebra, definition, 2}
	\end{equation}
is called
\AddIndex{morphism of \Ts representations from $f$ into $g$}
{morphism of representations from f into g}.
We also say that
\AddIndex{morphism of \Ts representations of $\mathfrak{F}$\Hyph algebra}
{morphism of representations of F algebra} is defined.
		\qed
		\end{definition}

		\begin{remark}
Let us consider morphism of \Ts representations
\eqref{eq: morphism of representations of F algebra, definition, 1}.
We denote elements of the set $B$ by letter using pattern $b\in B$.
However if we want to show that $b$ is image of element $a\in A$,
we use notation $\textcolor{red}{r(a)}$.
Thus equation
\[
\textcolor{red}{r(a)}=r(a)
\]
means that $\textcolor{red}{f(a)}$ (in left part of equation)
is image $a\in A$ (in right part of equation).
Using such considerations, we denote
element of set $N$ as $\textcolor{blue}{R(m)}$.
We will follow this convention when we consider correspondences
between homomorphisms of $\mathfrak{F}$\Hyph algebra
and mappings between sets where we defined corresponding\Ts representations.

There are two ways to interpret
\eqref{eq: morphism of representations of F algebra, definition, 2}
\begin{itemize}
\item Let \Ts transformation $f(a)$ map $m\in M$ into $f(a)m$.
Then \Ts transformation $g(\textcolor{red}{r(a)})$ maps
$\textcolor{blue}{R(m)}\in N$ into $\textcolor{blue}{R(f(a)m)}$.
\item We represent morphism of representations from $f$ into $g$
using diagram
\[
\xymatrix{
&M\ar[dd]_(.3){f(a)}\ar[rr]^R&&N\ar[dd]^(.3){g(\textcolor{red}{r(a)})}\\
&\ar @{}[rr]|{(1)}&&\\
&M\ar[rr]^R&&N\\
A\ar[rr]^r\ar@{=>}[uur]^(.3)f&&B\ar@{=>}[uur]^(.3)g
}
\]
From \eqref{eq: morphism of representations of F algebra, definition, 2}
it follows that diagram $(1)$ is commutative.
\end{itemize}
		\qed
		\end{remark}

		\begin{theorem}
Let us consider \Ts representation
\[
f:A\rightarrow {}^\star M
\]
of $\mathfrak{F}$\Hyph algebra $A$ 
and \Ts representation
\[
g:B\rightarrow {}^\star N
\]
of $\mathfrak{F}$\Hyph algebra $B$.
Morphism
\[
\xymatrix{
h:A\ar[r]&B&H:M\ar[r]&N
}
\]
of \Ts representations from $f$ into $g$
satisfies equation
	\begin{equation}
H(\omega(f(a_1),...,f(a_n))m)=\omega(g(h(a_1)),...,g(h(a_n)))H(m)
	\label{eq: morphism of representations of F algebra, 1}
	\end{equation}
for any $n$-ary operation $\omega$ of $\mathfrak{F}$\Hyph algebra.
		\end{theorem}
		\begin{proof}
Since $f$ is homomorphism, we have
	\begin{equation}
H(\omega(f(a_1),...,f(a_n))m)=H(f(\omega(a_1,...,a_n))m)
	\label{eq: morphism of representations of F algebra, 2}
	\end{equation}
From \eqref{eq: morphism of representations of F algebra, definition, 2} and
\eqref{eq: morphism of representations of F algebra, 2} it follows that
	\begin{equation}
H(\omega(f(a_1),...,f(a_n))m)=g(h(\omega(a_1,...,a_n)))H(m)
	\label{eq: morphism of representations of F algebra, 3}
	\end{equation}
Since $h$ is homomorphism, from
\eqref{eq: morphism of representations of F algebra, 3} it follows that
	\begin{equation}
H(\omega(f(a_1),...,f(a_n))m)=g(\omega(h(a_1),...,h(a_n)))H(m)
	\label{eq: morphism of representations of F algebra, 4}
	\end{equation}
Since $g$ is homomorphism,
\eqref{eq: morphism of representations of F algebra, 1} follows from
\eqref{eq: morphism of representations of F algebra, 4}.
		\end{proof}

		\begin{theorem}
	\label{theorem: morphism of representations of F algebra}
Given single transitive \Ts representation
\[
f:A\rightarrow {}^\star M
\]
of $\mathfrak{F}$\Hyph algebra $A$ and single transitive \Ts representation
\[
g:B\rightarrow {}^\star N
\]
of $\mathfrak{F}$\Hyph algebra $B$,
there exists morphism
\[
\xymatrix{
p:A\ar[r]&B&P:M\ar[r]&N
}
\]
of \Ts representations from $f$ into $g$.
		\end{theorem}
		\begin{proof}
Let us choose homomorphism $h$.
Let us choose element $m\in M$
and element $n\in N$. 
To define mapping $H$, let us consider following diagram
\[
\xymatrix{
&M\ar[dd]^(.3)a\ar[rr]^H&&N\ar[dd]_(.3){p(a)}\\
&\ar @{}[rr]|{(1)}&&\\
&M\ar[rr]^H&&N\\
A\ar[rr]^p\ar@{=>}[uur]^(.3)f&&B\ar@{=>}[uur]^(.3)g
}
\]
From commutativity of diagram $(1)$, it follows that
\[
H(am)=p(a)H(m)
\]
For arbitrary $m'\in M$, we defined unambiguously $a\in A$ such
that $m'=am$. Therefore, we defined mapping $H$
which satisfies to equation
\eqref{eq: morphism of representations of F algebra, definition, 2}.
		\end{proof}

		\begin{theorem}
Given single transitive \Ts representation
\[
f:A\rightarrow {}^\star M
\]
of $\mathfrak{F}$\Hyph algebra $A$,
for any automorphism of $\mathfrak{F}$\Hyph algebra $A$ there exists morphism
\[
\xymatrix{
p:A\ar[r]&A&P:M\ar[r]&M
}
\]
of \Ts representations from $f$ into $f$.
		\end{theorem}
		\begin{proof}
Let us consider following diagram
\[
\xymatrix{
&M\ar[dd]^(.3)a\ar[rr]^H&&N\ar[dd]_(.3){p(a)}\\
&\ar @{}[rr]|{(1)}&&\\
&M\ar[rr]^H&&N\\
A\ar[rr]^p\ar@{=>}[uur]^(.3)f&&A\ar@{=>}[uur]^(.3)g
}
\]
Statement of theorem is corollary of theorem
\ref{theorem: morphism of representations of F algebra}.
		\end{proof}

		\begin{theorem}
		\label{theorem: product of morphisms of representations of F algebra}
Let
\[
f:A\rightarrow {}^\star M
\]
be \Ts representation of $\mathfrak{F}$\Hyph algebra $A$,
\[
g:B\rightarrow {}^\star N
\]
be \Ts representation of $\mathfrak{F}$\Hyph algebra $B$,
\[
h:C\rightarrow {}^\star L
\]
be \Ts representation of $\mathfrak{F}$\Hyph algebra $C$.
Given morphisms of \Ts representations of $\mathfrak{F}$\Hyph algebra
\[
\xymatrix{
p:A\ar[r]&B&P:M\ar[r]&N
}
\]
\[
\xymatrix{
q:B\ar[r]&C&Q:N\ar[r]&L
}
\]
There exists morphism of \Ts representations of $\mathfrak{F}$\Hyph algebra
\[
\xymatrix{
r:A\ar[r]&C&R:M\ar[r]&L
}
\]
where $r=qp$, $R=QP$.
We call morphism $(r,R)$ of \Ts representations from $f$ into $h$
\AddIndex{product of morphisms  $(p,P)$ and $(q,Q)$
of \Ts representations of $\mathfrak{F}$\Hyph algebra}
{product of morphisms of representations of F algebra}.
		\end{theorem}
		\begin{proof}
Map $r$ is homomorphism of $\mathfrak{F}$\Hyph algebra $A$ into
$\mathfrak{F}$\Hyph algebra $C$.
We need to show that tuple of maps $(r,R)$ satisfies to
\eqref{eq: morphism of representations of F algebra, definition, 2}:
\begin{align*}
\textcolor{blue}{R(f(a)m)}&=\textcolor{blue}{QP(f(a)m)}\\
&=\textcolor{blue}{Q(g(\textcolor{red}{p(a)})\textcolor{blue}{P(m)})}\\
&=h(\textcolor{red}{qp(a)})\textcolor{blue}{QP(m)})\\
&=h(\textcolor{red}{r(a)})\textcolor{blue}{R(m)}
\end{align*}
		\end{proof}

Representations and morphisms of representations of $\mathfrak{F}$\Hyph algebra
create category of representations of $\mathfrak{F}$\Hyph algebra.

\begin{definition}
Let us define equivalence $S$ on the set $M$.
\Ts transformation $f$ is called
\AddIndex{coordinated with equivalence}{transformation coordinated with equivalence} $S$,
when
$fm_1\equiv fm_2(\mathrm{mod} S)$ follows from condition $m_1\equiv m_2(\mathrm{mod} S)$.
\end{definition}

\begin{theorem}
\label{theorem: transformation correlated with equivalence}
Let us consider equivalence $S$ on set $M$.
Let us consider $\mathfrak{F}$\Hyph algebra on set ${}^\star M$.
Since \Ts transformations are coordinated with equivalence $S$,
we can define the structure of $\mathfrak{F}$\Hyph algebra
on the set ${}^\star(M/S)$.
\end{theorem}
\begin{proof}
Let $h=\mathrm{nat}\ S$. If $m_1\equiv m_2(\mathrm{mod} S)$,
then $h(m_1)=h(m_2)$. Since
$f\in{}^\star M$ is coordinated with equivalence $S$,
then $h(f(m_1))=h(f(m_2))$. This allows to define
\Ts transformation $F$ according to rule
\[
F([m])=h(f(m))
\]

Let $\omega$ be n\Hyph ary operation of $\mathfrak{F}$\Hyph algebra.
Suppose $f_1$, ..., $f_n\in{}^\star M$ and
\begin{align*}
F_1([m])&=h(f_1(m))&&...&F_n([m])&=h(f_n(m))
\end{align*}
We define operation on the set ${}^\star(M/S)$ according to rule
\[
\omega(F_1,...,F_n)[m]=h(\omega(f_1,...,f_n)m)
\]
This definition is proper because
$\omega(f_1,...,f_n)\in{}^\star M$ and is coordinated with equivalence $S$.
\end{proof}

		\begin{theorem}
		\label{theorem: decompositions of morphism of representations}
Let
\[
f:A\rightarrow {}^\star M
\]
be \Ts representation of $\mathfrak{F}$\Hyph algebra $A$,
\[
g:B\rightarrow {}^\star N
\]
be \Ts representation of $\mathfrak{F}$\Hyph algebra $B$.
Let
\[
\xymatrix{
r:A\ar[r]&B&R:M\ar[r]&N
}
\]
be morphism of representations from $f$ into $g$.
Suppose
\begin{align*}
s&=rr^{-1}&S&=RR^{-1}
\end{align*}
Then there exist decompositions of $r$ and $R$,
which we describe using diagram
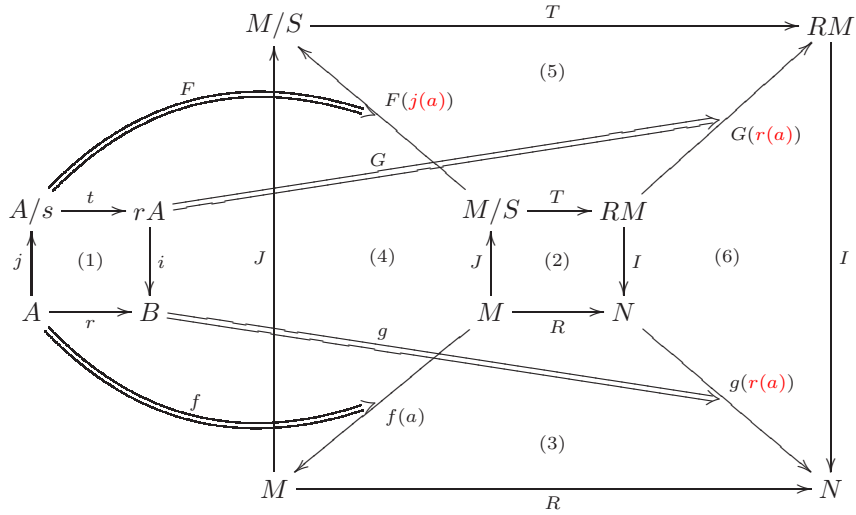
\[
\xymatrix{
&&M/S\ar[rrrrr]^T\ar@{}[drrrrr]|{(5)}&&\ar@{}[dddddll]|{(4)}&
\ar@{}[dddddrr]|{(6)}&&RM\ar[ddddd]^I\\
&&&&&&&\\
A/s\ar[r]^t\ar@/^2pc/@{=>}[urrr]^F&
rA\ar[d]^i\ar@{=>}[urrrrr]^(.4)G&&&
M/S\ar[r]^T\ar[lluu]_{F(\textcolor{red}{j(a)})}&
RM\ar[d]^I\ar[rruu]_{G(\textcolor{red}{r(a)})}\\
A\ar[r]_r\ar[u]^j\ar@{}[ur]|{(1)}\ar@/_2pc/@{=>}[drrr]^f&
B\ar@{=>}[drrrrr]^(.4)g&&&
M\ar[r]_R\ar[u]^J\ar@{}[ur]|{(2)}\ar[ddll]^{f(a)}&
N\ar[ddrr]^{g(\textcolor{red}{r(a)})}\\
&&&&&&&\\
&&M\ar[uuuuu]^J\ar[rrrrr]_R\ar@{}[urrrrr]|{(3)}&&&&&N
}
\]
\begin{itemize}
\item $s=\mathrm{ker}\ r$ is a congruence on $A$.
There exists decompositions of homomorphism $r$
	\begin{equation}
r=itj
	\label{eq: morphism of representations of algebra, homomorphism, 1}
	\end{equation}
$j=\mathrm{nat}\ s$ is the natural homomorphism
	\begin{equation}
\textcolor{red}{j(a)}=j(a)
	\label{eq: morphism of representations of algebra, homomorphism, 2}
	\end{equation}
$t$ is isomorphism
	\begin{equation}
\textcolor{red}{r(a)}=t(\textcolor{red}{j(a)})
	\label{eq: morphism of representations of algebra, homomorphism, 3}
	\end{equation}
$i $ is the inclusion mapping
	\begin{equation}
r(a)=i(\textcolor{red}{r(a)})
	\label{eq: morphism of representations of algebra, homomorphism, 4}
	\end{equation}
\item $S=\mathrm{ker}\ R$ is an equivalence on $M$.
There exists decompositions of homomorphism $R$
	\begin{equation}
R=ITJ
	\label{eq: morphism of representations of algebra, map, 1}
	\end{equation}
$J=\mathrm{nat}\ S$ is surjection
	\begin{equation}
\textcolor{blue}{J(m)}=J(m)
	\label{eq: morphism of representations of algebra, map, 2}
	\end{equation}
$T$ is bijection
	\begin{equation}
\textcolor{blue}{R(m)}=T(\textcolor{blue}{J(m)})
	\label{eq: morphism of representations of algebra, map, 3}
	\end{equation}
$I$ is the inclusion mapping
	\begin{equation}
R(m)=I(\textcolor{blue}{R(m)})
	\label{eq: morphism of representations of algebra, map, 4}
	\end{equation}
\item $F$ is \Ts representation of $\mathfrak{F}$\Hyph algebra $A/s$ in $M/S$
\item $G$ is \Ts representation of $\mathfrak{F}$\Hyph algebra $rA$ in $RM$
\item There exists decompositions of morphism of representations
\[
(r,R)=(i,I)(t,T)(j,J)
\]
\end{itemize}
		\end{theorem}
		\begin{proof}
Existence of diagram $(1)$ follows from theorem II.3.7
(\citeBib{Cohn: Universal Algebra}, p. 60).
Existence of diagram $(2)$ follows from theorem I.3.1
(\citeBib{Cohn: Universal Algebra}, p. 15).

We start from diagram $(4)$.

Let $m_1\equiv m_2(\mathrm{mod}\ S)$.
%Тогда $J(m_1)=J(m_2)$,
Then
\begin{equation}
\textcolor{blue}{R(m_1)}=\textcolor{blue}{R(m_2)}
\label{eq: morphism of representations of algebra, (4), 1}
\end{equation}
%$JT(m_1)=JT(m_2)$,
%так как $JT$ - отображение на $MR$.
Since $a_1\equiv a_2(\mathrm{mod}s)$, then
\begin{equation}
\textcolor{red}{r(a_1)}=\textcolor{red}{r(a_2)}
\label{eq: morphism of representations of algebra, (4), 2}
\end{equation}
Therefore,
$j(a_1)=j(a_2)$.
Since $(r,R)$ is morphism of representations, then
\begin{equation}
\textcolor{blue}{R(f(a_1)m_1)}=g(\textcolor{red}{r(a_1)})\textcolor{blue}{R(m_1)}
\label{eq: morphism of representations of algebra, (4), 3}
\end{equation}
\begin{equation}
\textcolor{blue}{R(f(a_2)m_2)}=g(\textcolor{red}{r(a_2)})\textcolor{blue}{R(m_2)}
\label{eq: morphism of representations of algebra, (4), 4}
\end{equation}
From \eqref{eq: morphism of representations of algebra, (4), 1},
\eqref{eq: morphism of representations of algebra, (4), 2},
\eqref{eq: morphism of representations of algebra, (4), 3},
\eqref{eq: morphism of representations of algebra, (4), 4},
it follows that
\begin{equation}
\textcolor{blue}{R(f(a_1)m_1)}=\textcolor{blue}{R(f(a_2)m_2)}
\label{eq: morphism of representations of algebra, (4), 5}
\end{equation}
From \eqref{eq: morphism of representations of algebra, (4), 5}
it follows
\begin{equation}
f(a_1)m_1\equiv f(a_2)m_2(\mathrm{mod} S)
\label{eq: morphism of representations of algebra, (4), 6}
\end{equation}
and, therefore,
\begin{equation}
\textcolor{blue}{J(f(a_1)m_1)}=\textcolor{blue}{J(f(a_2)m_2)}
\label{eq: morphism of representations of algebra, (4), 7}
\end{equation}
From \eqref{eq: morphism of representations of algebra, (4), 7}
it follows that we defined map
\begin{equation}
F(\textcolor{red}{j(a)})(\textcolor{blue}{J(m)})
=\textcolor{blue}{J(f(a)m))}
\label{eq: morphism of representations of algebra, (4), 8}
\end{equation}
reasonably and this map is \Ts transformation of set $M/S$.

From equation
\eqref{eq: morphism of representations of algebra, (4), 6}
(in case $a_1=a_2$) it follows that for any $a$ \Ts transformation
is coordinated with equivalence $S$.
From theorem
\ref{theorem: transformation correlated with equivalence} it follows that
we defined structure of $\mathfrak{F}$\Hyph algebra
on the set ${}^\star(M/S)$.
Let us consider $n$\Hyph ary operation $\omega$ and $n$
\Ts transformations
\begin{align*}
F(\textcolor{red}{j(a_i)})\textcolor{blue}{J(m)}&=
\textcolor{blue}{J(f(a_i)m))}
&i&=1,...,n
\end{align*}
of the set $M/S$. We assume
\[
\omega(F(\textcolor{red}{j(a_1)}),...,F(\textcolor{red}{j(a_n)}))
\textcolor{blue}{J(m)}=
J(\omega(f(a_1),...,f(a_n)))m)
\]
Therefore, map $F$ is representations of
$\mathfrak{F}$\Hyph algebra $A/s$.

From \eqref{eq: morphism of representations of algebra, (4), 8}
it follows that $(j,J)$ is morphism of representations $f$ and $F$.

Let us consider diagram $(5)$.

Since $T$ is bijection, then we identify elements of the set
$M/S$ and the set $MR$, and this identification has form
\begin{equation}
T(\textcolor{blue}{J(m)})=\textcolor{blue}{R(m)}
\label{eq: morphism of representations of algebra, (5), 1}
\end{equation}
We can write \Ts transformation $F(\textcolor{red}{j(a)})$
of the set $M/S$ as
\begin{equation}
F(\textcolor{red}{j(a)}):\textcolor{blue}{J(m)}\rightarrow
F(\textcolor{red}{j(a)})\textcolor{blue}{J(m)}
\label{eq: morphism of representations of algebra, (5), 2}
\end{equation}
%Из \eqref{eq: morphism of representations of algebra, (4), 8}
%следует, что мы можем записать
%\eqref{eq: morphism of representations of algebra, (5), 2} в виде
%\begin{equation}
%F(\textcolor{red}{j(a)}):\textcolor{blue}{J(m)}\rightarrow
%\textcolor{blue}{J(f(a)m))}
%\label{eq: morphism of representations of algebra, (5), 3}
%\end{equation}
Since $T$ is bijection, we define \Ts transformation
\begin{equation}
T(\textcolor{blue}{J(m)})\rightarrow
T(F(\textcolor{red}{j(a)})\textcolor{blue}{J(m)})
\label{eq: morphism of representations of algebra, (5), 3}
\end{equation}
of the set $RM$. \Ts transformation
\eqref{eq: morphism of representations of algebra, (5), 3}
depends on $\textcolor{red}{j(a)}\in A/s$.
Since $t$ is bijection, we identify elements of the set
$A/s$ and the set $rA$, and this identification has form
\[
t(\textcolor{red}{j(a)})=\textcolor{red}{r(a)}
\]
Therefore, we defined map
\[
G:rA\rightarrow{}^\star RM
\]
according to equation
\begin{equation}
G(t(\textcolor{red}{j(a)}))T(\textcolor{blue}{J(m)})=
T(F(\textcolor{red}{j(a)})\textcolor{blue}{J(m)})
\label{eq: morphism of representations of algebra, (5), 4}
\end{equation}

Let us consider $n$\Hyph ary operation $\omega$ and $n$
\Ts transformations
\begin{align*}
G(\textcolor{red}{r(a_i)})\textcolor{blue}{R(m)}&=
T(F(\textcolor{red}{j(a_i)})\textcolor{blue}{J(m)})
&i&=1,...,n
\end{align*}
of space $RM$. We assume
\[
\omega(G(\textcolor{red}{r(a_1)}),...,G(\textcolor{red}{r(a_n)}))
\textcolor{blue}{R(m)}=
T(\omega(F(\textcolor{red}{j(a_1)},...,F(\textcolor{red}{j(a_n)}))
\textcolor{blue}{J(m)})
\]
According to \eqref{eq: morphism of representations of algebra, (5), 4}
operation $\omega$ is defined reasonably on the set ${}^\star RM$.
Therefore, the map $G$ is representations of
$\mathfrak{F}$\Hyph algebra.

From \eqref{eq: morphism of representations of algebra, (5), 4}
it follows that $(t,T)$ is morphism of representations $F$ and $G$.

Diagram $(6)$ is the most simple case in our prove.
Since map $I$ is immersion and diagram $(2)$ is
commutative, we identify $n\in N$ and $\textcolor{blue}{R(m)}$
when $n\in\textrm{Im}R$. Similarly, we identify
corresponding \Ts transformations.
\begin{equation}
g'(i(\textcolor{red}{r(a)}))I(\textcolor{blue}{R(m)})=
I(G(\textcolor{red}{r(a)})\textcolor{blue}{R(m)})
\label{eq: morphism of representations of algebra, (6), 1}
\end{equation}
\[
\omega(g'(\textcolor{red}{r(a_1)}),...,g'(\textcolor{red}{r(a_n)}))
\textcolor{blue}{R(m)}=
I(\omega(G(\textcolor{red}{r(a_1)},...,G(\textcolor{red}{r(a_n)}))
\textcolor{blue}{R(m)})
\]
Therefore,
$(i,I)$ is morphism of representations $G$ and $g$.

To prove the theorem we need to show that defined
in the proof \Ts representation $g'$ is congruent
with representation $g$, and operations over transformations are congruent
with corresponding operations over ${}^\star N$.
\begin{align*}
g'(i(\textcolor{red}{r(a)}))I(\textcolor{blue}{R(m)})
&=I(G(\textcolor{red}{r(a)})\textcolor{blue}{R(m)})
&&\textrm{by}\ \eqref{eq: morphism of representations of algebra, (6), 1}\\
&=I(G(t(\textcolor{red}{j(a)}))T(\textcolor{blue}{J(m)}))
&&\textrm{by}\ \eqref{eq: morphism of representations of algebra, homomorphism, 3},
\eqref{eq: morphism of representations of algebra, map, 3}\\
&=IT(F(\textcolor{red}{j(a)})\textcolor{blue}{J(m)})
&&\textrm{by}\ \eqref{eq: morphism of representations of algebra, (5), 4}\\
&=ITJ(f(a)m)
&&\textrm{by}\ \eqref{eq: morphism of representations of algebra, (4), 8}\\
&=R(f(a)m)
&&\textrm{by}\ \eqref{eq: morphism of representations of algebra, map, 1}\\
&=g(r(a))R(m)
&&\textrm{by}\ \eqref{eq: morphism of representations of F algebra, definition, 2}
\end{align*}
\begin{align*}
\omega(G(\textcolor{red}{r(a_1)}),...,G(\textcolor{red}{r(a_n)}))
\textcolor{blue}{R(m)}
&=T(\omega(F(\textcolor{red}{j(a_1)},...,F(\textcolor{red}{j(a_n)}))
\textcolor{blue}{J(m)})\\
&=T(F(\omega(\textcolor{red}{j(a_1)},...,\textcolor{red}{j(a_n)}))
\textcolor{blue}{J(m)})\\
&=T(F(j(\omega(a_1,...,a_n)))
\textcolor{blue}{J(m)})\\
&=T(J(f(\omega(a_1,...,a_n))m))\\
\end{align*}
		\end{proof}

From theorem \ref{theorem: decompositions of morphism of representations}
it follows that we can reduce the problem of studying of morphism of \Ts representations
of $\mathfrak{F}$\Hyph algebra
to the case described by diagram
\begin{equation}
\xymatrix{
&M\ar[dd]_(.3){f(a)}\ar[rr]^J&&M/S\ar[dd]^(.3){F(\textcolor{red}{j(a)})}\\
&&&\\
&M\ar[rr]^J&&M/S\\
A\ar[rr]^j\ar@{=>}[uur]^(.3)f&&A/s\ar@{=>}[uur]^(.3)F
}
\label{eq: morphism of representations of algebra, reduce, 1}
\end{equation}

		\begin{theorem}
		\label{theorem: morphism of representations of algebra, reduce}
We can supplement diagram \eqref{eq: morphism of representations of algebra, reduce, 1}
with \Ts representation $F_1$ of
$\mathfrak{F}$\Hyph algebra $A$ into set $M/S$ such
that diagram
\begin{equation}
\xymatrix{
&M\ar[dd]_(.3){f(a)}\ar[rr]^J&&M/S\ar[dd]^(.3){F(\textcolor{red}{j(a)})}\\
&&&\\
&M\ar[rr]^(.65)J&&M/S\\
A\ar[rr]^j\ar@{=>}[uur]^(.3)f\ar@{=>}[uurrr]^(.7){F_1}&&A/s\ar@{=>}[uur]^(.3)F
}
\label{eq: morphism of representations of algebra, reduce, 2}
\end{equation}
is commutative. The set of \Ts transformations of \Ts representation $F$ and
the set of \Ts transformations of \Ts representation $F_1$ coincide.
		\end{theorem}
		\begin{proof}
To prove theorem it is enough to assume
\[F_1(a)=F(\textcolor{red}{j(a)})\]
Since map $j$ is surjection, then $\mathrm{Im}F_1=\mathrm{Im}F$.
Since $j$ and $F$ are homomorphisms of $\mathfrak{F}$\Hyph algebra, then $F_1$
is also homomorphism of $\mathfrak{F}$\Hyph algebra.
		\end{proof}

Theorem \ref{theorem: morphism of representations of algebra, reduce}
completes the series of theorems dedicated to the structure of morphism of \Ts representations
$\mathfrak{F}$\Hyph algebra.
From these theorems it follows that we can simplify task of studying of morphism of \Ts representations
$\mathfrak{F}$\Hyph algebra and not go beyond morphism of \Ts representations of form
\[
\xymatrix{
id:A\ar[r]&A&R:M\ar[r]&N
}
\]
In this cae we identify morphism of $(id,R)$ \Ts representations
of $\mathfrak{F}$\Hyph algebra and map $R$.

		\begin{definition}
Let
\[
f:A\rightarrow {}^\star M
\]
be \Ts representation of $\mathfrak{F}$\Hyph algebra $A$
in $\mathfrak{H}$\Hyph algebra $M$ and
\Ts representation 
\[
g:B\rightarrow {}^\star N
\]
be \Ts representation of $\mathfrak{F}$\Hyph algebra $B$
in $\mathfrak{H}$\Hyph algebra $N$.
Morphism $(h,H)$ of \Ts representations of algebra of type $\mathfrak{F}$
is called
\AddIndex{morphism of \Ts representations of $\mathfrak{F}$\Hyph algebra
in $\mathfrak{H}$\Hyph algebra}
{morphism of representations of F algebra in H algebra}.
		\qed
		\end{definition}

%auto-ignore
\def\texDrcMorphism{}

			\section{\texorpdfstring{$D\RCstar$}{Drc}\hyph Linear Map of Vector Spaces}
			\label{section: Linear Map of Vector Spaces}

		\begin{definition}
		\label{definition: src trc linear map of vector spaces}
Suppose $\mathcal{A}$ is $S\RCstar$\Hyph vector space.
Suppose $\mathcal{B}$ is $T\RCstar$\Hyph vector space.
Morphism
	\[
\xymatrix{
f:S\ar[r]&T&\Vector{A}:\mathcal{A}\ar[r]&\mathcal{B}
}
	\]
of \Ts representations of skew field 
in Abelian group is called
\AddIndex{$(S\RCstar,T\RCstar)$\Hyph linear map of vector spaces}
{src trc linear map of vector spaces}.
		\qed
		\end{definition}

By theorem
\ref{theorem: decompositions of morphism of representations}
studying $(S\RCstar,T\RCstar)$\Hyph linear map
we can consider case $S=T$.

		\begin{definition}
		\label{definition: drc linear map of vector spaces}
Suppose $\mathcal{A}$ and
$\mathcal{B}$ are \drc vector spaces.
We call map
	\begin{equation}
\Vector{A}:\mathcal{A}\rightarrow \mathcal{B}
	\label{eq: drc linear map of vector spaces}
	\end{equation}
\AddIndex{\drc linear map of vector spaces}{drc linear map of vector spaces}
if\footnote{Expression $a\RCstar \Vector{A}(\Vector{m})$ means expression
$a^{\gi a}\ \Vector{A}({}_{\gi a}\Vector{m})$}
	\begin{equation}
	\label{eq: drc linear map of vector spaces, product on scalar}
\Vector{A}(a\RCstar \Vector{m})=a\RCstar \Vector{A}(\Vector{m})
	\end{equation}
for any $a^{\gi a} \in D$, ${}_{\gi a}\Vector{m} \in \mathcal{A}$.
		\qed
		\end{definition}

		\begin{theorem}
		\label{theorem: drc linear map of drc vector space}
Let $\Vector{f}=({}_{\gi a}\Vector{f},\gi{a}\in \gi{I})$
be a \drc basis of vector space $\mathcal{A}$
and $\Vector{e}=({}_{\gi b}\Vector{e},\gi{b}\in \gi{J})$ 
be a \drc basis of vector space $\mathcal{B}$.
Then \drc linear map \eqref{eq: drc linear map of vector spaces}
of vector spaces has presentation
	\begin{equation}
b=a\RCstar A
	\label{eq: drc linear map of vector spaces, presentation}
	\end{equation}
relative to selected bases. Here
\begin{itemize}
\item $a$ is coordinate matrix of vector
$\Vector a$ relative the \drc basis $\Basis{f}$
\item $b$ is coordinate matrix of vector
$$\Vector b=\Vector{A}(\Vector a)$$ relative the \drc basis $\Basis{e}$
\item $A$ is coordinate matrix of set of vectors
$(\Vector A({}_{\gi a}\Vector{f}))$ in \drc basis $\Basis{e}$
called \AddIndex{matrix of \drc linear map}{matrix of drc linear map}
relative bases $\Basis{f}$ and $\Basis{e}$
\end{itemize}
		\end{theorem}
		\begin{proof}
Vector $\Vector{a}\in \mathcal{A}$ has expansion
	\[
\Vector{a}=a\RCstar \Vector{f}
	\]
relative to \drc basis $\Basis{f}$.
Vector $\Vector{b}=f(\Vector{a})\in \mathcal{B}$ has expansion
	\begin{equation}
\Vector{b}=b\RCstar \Vector{e}
	\label{eq: drc linear map of vector spaces, presentation, 1}
	\end{equation}
relative to \drc basis $\Basis{e}$.

Since $\Vector A$ is a \drc linear map, from
\eqref{eq: drc linear map of vector spaces, product on scalar} it follows
that
	\begin{equation}
\Vector{b}=\Vector A(\Vector{a})=\Vector A(a\RCstar \Vector{f})
=a\RCstar \Vector A(\Vector{f})
	\label{eq: drc linear map of vector spaces, presentation, 2}
	\end{equation}
$\Vector A({}_{\gi a}\Vector{f})$ is also a vector of $\mathcal{B}$ and has
expansion
	\begin{equation}
\Vector A({}_{\gi a}\Vector{f})
={}_{\gi a}A\RCstar \Vector{e}
={}_{\gi a}A^{\gi b}\ {}_{\gi b}\Vector{e}
	\label{eq: drc linear map of vector spaces, presentation, 3}
	\end{equation}
relative to basis $\Basis{e}$.
Combining \eqref{eq: drc linear map of vector spaces, presentation, 2}
and \eqref{eq: drc linear map of vector spaces, presentation, 3} we get
	\begin{equation}
\Vector{b}=a\RCstar A\RCstar\Vector{e}
	\label{eq: drc linear map of vector spaces, presentation, 4}
	\end{equation}
\eqref{eq: drc linear map of vector spaces, presentation}  follows
from comparison of \eqref{eq: drc linear map of vector spaces, presentation, 1}
and \eqref{eq: drc linear map of vector spaces, presentation, 4} and
theorem \xRef{0701.238}{theorem: expansion relative basis, vector space}.
		\end{proof}

On the basis of theorem \ref{theorem: drc linear map of drc vector space}
we identify the \drc linear map \eqref{eq: drc linear map of vector spaces} of vector spaces
and the matrix of its presentation \eqref{eq: drc linear map of vector spaces, presentation}.

		\begin{theorem}
		\label{theorem: product of drc linear maps, vector spaces}
Let $$\Basis{f}=({}_{\gi a}\Vector{f},\gi{a}\in \gi{I})$$ be a \drc basis of vector space $\mathcal{A}$,
$$\Basis{e}=({}_{\gi b}\Vector{e},\gi{b}\in \gi{J})$$ be a \drc basis of vector space $\mathcal{B}$,
and $$\Basis{g}=({}_{\gi c}\Vector{g},\gi{c}\in \gi{L})$$ be a \drc basis of vector space $\mathcal{C}$.
Suppose diagram of \drc linear maps
	\[
\xymatrix{
\mathcal{A}\ar[rr]^C\ar[dr]^A & & \mathcal{C}\\
&\mathcal{B}\ar[ur]^B &
}
	\]
is commutative diagram
where \drc linear map $A$  has presentation
	\begin{equation}
	\label{eq: product of drc linear maps, left vector spaces, A}
b=a\RCstar A
	\end{equation}
relative to selected bases
and \drc linear map $B$ has presentation
	\begin{equation}
	\label{eq: product of drc linear maps, left vector spaces, B}
c=b\RCstar B
	\end{equation}
relative to selected bases.
Then \drc linear map $C$ has presentation
	\begin{equation}
	\label{eq: product of drc linear maps, left vector spaces, C}
c=a\RCstar A\RCstar B
	\end{equation}
relative to selected bases.
		\end{theorem}
		\begin{proof}
Proof of the statement follows from substituting
\eqref{eq: product of drc linear maps, left vector spaces, A}
into \eqref{eq: product of drc linear maps, left vector spaces, B}.
		\end{proof}

Presenting \drc linear map as \RC product we can rewrite
\eqref{eq: drc linear map of vector spaces, product on scalar} as
	\begin{equation}
(ka)\RCstar A=k(a\RCstar A)
	\label{eq: drc linear map of vector spaces, product on scalar, image}
	\end{equation}
We can express the statement of the theorem
\ref{theorem: product of drc linear maps, vector spaces}
in the next form
	\begin{equation}
(a\RCstar A)\RCstar B=a\RCstar (A\RCstar B)
	\label{eq: product of drc linear maps, vector spaces, image}
	\end{equation}
Equations
\eqref{eq: drc linear map of vector spaces, product on scalar, image} and
\eqref{eq: product of drc linear maps, vector spaces, image}
represent the
\AddIndex{associative law for \drc linear maps of vector spaces}
{associative law for drc linear maps of vector spaces}.
This allows us writing of such expressions without using of brackets.

Equation \eqref{eq: drc linear map of vector spaces, presentation}
is coordinate notation for \drc linear map.
Based theorem \ref{theorem: drc linear map of drc vector space}
non coordinate notation also can be expressed using \RC product
	\begin{equation}
\Vector{b}=\Vector{a}\RCstar\Vector A=a\RCstar \Vector{f}\RCstar\Vector A
=a\RCstar A\RCstar \Vector{e}
	\label{eq: drc linear map of vector spaces, presentation, 5}
	\end{equation}
If we substitute equation \eqref{eq: drc linear map of vector spaces, presentation, 5}
into theorem
\ref{theorem: product of drc linear maps, vector spaces},
then we get chain of equations
	\begin{align*}
\Vector{c}&=\Vector{b}\RCstar\Vector B=b\RCstar \Vector{e}\RCstar\Vector B
=b\RCstar B\RCstar \Vector{g}\\
\Vector{c}&=\Vector{a}\RCstar\Vector A\RCstar\Vector B
=a\RCstar \Vector{f}\RCstar\Vector A\RCstar\Vector B
=a\RCstar A\RCstar B\RCstar \Vector{g}
	\end{align*}

\begin{remark}
One can easily see from the an example of \drc linear map
how theorem 
\ref{theorem: decompositions of morphism of representations}
makes our reasoning simpler in study of the morphism of \Ts representations
of $\mathfrak{F}$\Hyph algebra.
In the framework of this remark,
we agree to call the
theory of \drc linear mappings reduced theory,
and theory stated in this remark is called enhanced theory.

Suppose $\mathcal{A}$ is $S\RCstar$\Hyph vector space.
Suppose $\mathcal{B}$ is $T\RCstar$\Hyph vector space.
Suppose
	\[
\xymatrix{
f:S\ar[r]&T&\Vector{A}:\mathcal{A}\ar[r]&\mathcal{B}
}
	\]
is $(S\RCstar,T\RCstar)$\Hyph linear map of vector spaces.
Let $\Basis{f}=({}_{\gi a}\Vector{f},\gi{a}\in \gi{I})$
be a $S\RCstar$\Hyph basis of vector space $\mathcal{A}$
and $\Basis{e}=({}_{\gi b}\Vector{e},\gi{b}\in \gi{J})$
be a $T\RCstar$\Hyph basis of vector space $\mathcal{B}$.

From definitions
\ref{definition: src trc linear map of vector spaces} and
\ref{definition: morphism of representations of F algebra}
it follows
	\begin{equation}
\Vector{b}=\Vector A(\Vector{a})=\Vector A(a\RCstar \Vector{f})
=f(a)\RCstar \Vector A(\Vector{f})
	\label{eq: src trc linear map of vector spaces, presentation, 1}
	\end{equation}
$\Vector A({}_{\gi a}\Vector{f})$ is also a vector of $\mathcal{B}$ and has
expansion
	\begin{equation}
\Vector A({}_{\gi a}\Vector{f})
={}_{\gi a}A\RCstar \Vector{e}
={}_{\gi a}A^{\gi b}\ {}_{\gi b}\Vector{e}
	\label{eq: src trc linear map of vector spaces, presentation, 2}
	\end{equation}
relative to basis $\Basis{e}$.
Combining \eqref{eq: src trc linear map of vector spaces, presentation, 1}
and \eqref{eq: src trc linear map of vector spaces, presentation, 2}, we get
	\begin{equation}
\Vector{b}=f(a)\RCstar A\RCstar\Vector{e}
	\label{eq: src trc linear map of vector spaces, presentation, 3}
	\end{equation}

Suppose $\mathcal{C}$ is \drc vector space.
Suppose
	\[
\xymatrix{
g:T\ar[r]&D&\Vector{B}:\mathcal{B}\ar[r]&\mathcal{C}
}
	\]
is $(T\RCstar,D\RCstar)$\Hyph linear map of vector spaces.
Let $\Basis{h}=({}_{\gi a}\Vector{h},\gi{a}\in \gi{K})$
be \drc basis of vector space $\mathcal{C}$.
Then, according to
\eqref{eq: src trc linear map of vector spaces, presentation, 3},
the product of 
$(S\RCstar,T\RCstar)$\Hyph linear map $(f,\Vector{A})$ and
$(T\RCstar,D\RCstar)$\Hyph linear map $(g,\Vector{B})$
has form
	\begin{equation}
\Vector{c}=gf(a)\RCstar g(A)\RCstar B\RCstar\Vector{h}
	\label{eq: src trc linear map of vector spaces, presentation, 4}
	\end{equation}
Comparison of equations
\eqref{eq: product of drc linear maps, left vector spaces, C} and
\eqref{eq: src trc linear map of vector spaces, presentation, 4}
that extended theory of linear maps is more
complicated then reduced theory.

If we need we can use extended theory,
however we will not get new results comparing with
reduced theory. At the same time plenty of details
makes picture less clear and demands permanent attention.
		\qed
\end{remark}

%auto-ignore
\def\texBundle{}
\ifx\PrintBook\Defined
				\chapter{Cartesian Product of Bundles}
\fi

			\section{Bundle}

Let $M$ be a manifold and
\symb{\bundle{\mathcal{E}}pEM}0{bundle}
	\begin{equation}
	\label{eq: bundle, definition}
\ShowSymbol{bundle}
	\end{equation}
be a bundle over $M$ with fiber $E$.\footnote{Since I have
deal with different bundles I follow next agreement.
I use the same letter in different alphabets for notation of bundle and fiber.}
The symbol $\bundle{}pE{}$ means that
$E$ is a typical fiber of the bundle.
Set $\mathcal{E}$ is domain of
map $\bundle{}pE{}$.
Set $M$ is range of
map $\bundle{}pE{}$.
We identify the smooth map $\bundle{}pE{}$
and the bundle \eqref{eq: bundle, definition}.
Mapping $p$ is called
\AddIndex{projection of bundle $\mathcal{E}$ along fiber $E$}
{projection of bundle along fiber}.
Denote by
\symb{\Gamma(\bundle{}pE{})}1{set of sections of bundle}
the set of sections of bundle $\bundle{}pE{}$.

		\begin{definition}
		\label{definition: locally compact space}
A space is said to be \AddIndex{locally compact at point}{locally compact at point space} $p$
if there exists open set $U$, $p\in U$, whose closure $\overline U$ is compact.
A space is said to be \AddIndex{locally compact}{locally compact space}
if it is locally compact at each of its points.\footnote{I follow to
definition from \citeBib{Hocking Young Topology}, p. 71.}
		\qed
		\end{definition}

		\begin{definition}
		\label{definition: topology finer coarser}
Given topologies $\mathcal{T}_1$, $\mathcal{T}_2$ on the
set $X$, we say that $\mathcal{T}_1$ is finer than $\mathcal{T}_2$ and
that $\mathcal{T}_2$ is coarser than $\mathcal{T}_1$ if,
denoting by $X_i$ the set $X$ with the topology $\mathcal{T}_i$, $i=1,2$,
the identity mapping $X_1\rightarrow X_2$,
is continuous. If $\mathcal{T}_1\ne\mathcal{T}_2$,
we say that $\mathcal{T}_1$ is strictly finer than $\mathcal{T}_2$ and that
$\mathcal{T}_2$ is strictly coarser than $\mathcal{T}_1$.
		\qed
		\end{definition}

Let topology $\mathcal{T}_1$ in
the set $X$ be finer than topology $\mathcal{T}_2$. Let us consider
diagram
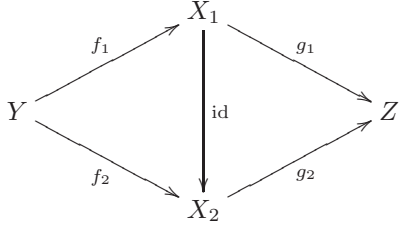
\[
\xymatrix{
&&X_1\ar[dd]^{\mathrm{id}}\ar[drr]^{g_1}\\
Y\ar[urr]^{f_1}\ar[drr]_{f_2}&&&&Z\\
&&X_2\ar[urr]_{g_2}
}
\]
According to definition \ref{definition: topology finer coarser},
if mapping $f_1$ is continues, than mapping $f_2$ is continues.
Similarly,
if mapping $g_2$ is continues, than mapping $g_1$ is continues.

Let $\bundle{\mathcal{E}}pEM$ be bundle.
Let us consider an open set $U\subset M$ such, that there exists
chart $\varphi$ over $U$
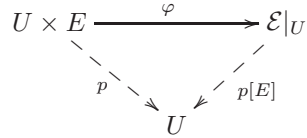
	\[
\xymatrix{
U\times E\ar[rr]^{\varphi}\arp[rd]_p
&&\mathcal{E}|_U\arp[ld]^{\bundle{}pE{}}\\
&U
}
	\]
Since $\varphi$ is homeomorphism, then topology in $U\times E$ and $\mathcal{E}|_U$
are comparable. Since $U\times E$ is Cartesian product of topological spaces
$U$ and $E$, then in $\mathcal{E}|_U$, we define coarsest topology
for which projection $\bundle{}pE{}$ is continuous
(\citeBib{Bourbaki: General Topology 1}, p. 31).

Let us consider relation $r$ in $\mathcal{E}$ such that $(p,q)\in r$ iff
$p$ and $q$ belong to the same fiber. Relation $r$ is equivalence.
$\bundle{}pE{}$ is natural mapping. Let us consider diagram
	
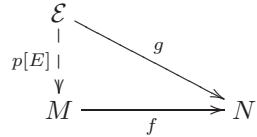
\[
\xymatrix{
\mathcal{E}\arp[d]_{\bundle{}pE{}}\ar[drr]^g\\
M\ar[rr]_f&&N
}
	\]
Continuity of mapping $g$ follows from continuity of mapping $f$.
Hence, we can define in $M$ quotient topology which is
the finest topology for which $\bundle{}pE{}$ is continuous
(\citeBib{Bourbaki: General Topology 1}, p. 39).

\AddIndex{Cartesian power $A$ of set $B$}{Cartesian power of set}
is the set \symb{B^A}1{Cartesian power of set}
of mappings $f:A\rightarrow B$
(\citeBib{Cohn: Universal Algebra}, page 5).
Let us consider subsets of $B^A$ of the form
\[
W_{K,U}=\{f:A\rightarrow B|f(K)\subset U\}
\]
where $K$ is compact subset of space $A$,
$U$ is open subset of space $B$.
Sets $W_{K,U}$ form base of topology
on space $B^A$. This topology is called 
\AddIndex{compact\hyph open topology}{compact open topology}.
Cartesian power $A$ of set $B$
equipped by compact\hyph open topology is called
\AddIndex{mapping space}{mapping space}
(\citeBib{Maunder: Algebraic Topology}, page 213).

According to \citeBib{Maunder: Algebraic Topology}, page 214,
given spaces $A$, $B$, $C$, $D$ and
mappings $f:A\rightarrow C$, $g:D\rightarrow B$
we define morphism of mapping spaces
\[
g^f:D^C\rightarrow B^A
\]
by law
\begin{align*}
g^f(h)&=fhg&h:&C\rightarrow D&g^f(h):&A\rightarrow B
\end{align*}
Thus, we can represent the morphism of mapping spaces
using diagram
	\[
\xymatrix{
A\ar[rr]^f\ar[d]_{g^f(h)}&&C\ar[d]^h\\
B&&D\ar[ll]^g
}
	\]

		\begin{remark}
		\label{remark: set of sections of bundle, topology}
Set $\Gamma(\mathcal{E})$ is subset
of set $\mathcal{E}^M$.
This is why for set of sections we can use definitions established
for mapping set.
We define sets $W_{K,U}$ by law
\[
W_{K,U}=\{f\in\Gamma(\mathcal{E})|f(K)\subset U\}
\]
where $K$ is compact subset of space $M$,
$U$ is open subset of space $\mathcal{E}$.
		\qed
		\end{remark}

	\begin{remark}
I use arrow $\xymatrix{\arp[r]&}$ to represent
projection of bundle on diagram.
	\qed
	\end{remark}

	\begin{remark}
I use arrow $\xymatrix{\ars[r]&}$ to represent
section of bundle on diagram.
	\qed
	\end{remark}

		\begin{definition}
		\label{definition: fibered morphism}
Let us consider bundles $$\bundle{\mathcal{A}}pAM$$
and $$\bundle{\mathcal{B}}qBN$$
Tuple of mapping
\symb{(\begin{array}{cc}
\mathcal{F}:\mathcal{A}\rightarrow\mathcal{B},&
f:M\rightarrow N
\end{array})}
0{fibered morphism from A into B}
	\begin{equation}
\ShowSymbol{fibered morphism from A into B}
	\label{eq: fibered morphism}
	\end{equation}
such, that diagram
	\[
\xymatrix{
\mathcal{A}\arp[d]^{\bundle{}pA{}}\ar[rr]^{\mathcal{F}}
&&\mathcal{B}\arp[d]^{\bundle{}qB{}}\\
M\ar[rr]^f&&N
}
	\]
is commutative, is called
\AddIndex{fibered morphism from bundle $\mathcal{A}$ into $\mathcal{B}$}
{fibered morphism from A into B}.
The map $f$ is the \AddIndex{base of map}{base of map}
$\mathcal{F}$.
The map $\mathcal{F}$ is the \AddIndex{lift of map}{lift of map} $f$.
		\qed
		\end{definition}

	\begin{theorem}
	\label{theorem: fibered morphism defines morphism of spaces of sections}
Suppose map $f$ is bijection.
Then the map $\mathcal{F}$ defines morphism $\mathcal{F}^{f^{-1}}$ of spaces of sections
$\Gamma(\bundle{}pA{}))$ to $\Gamma(\bundle{}qB{}))$
	\[
\xymatrix{
\mathcal{A}\ar[rr]^{\mathcal{F}}&&\mathcal{B}
&u'=\mathcal{F}^{f^{-1}}(u)=\mathcal{F}uf^{-1}\\
M\ars[u]^u="U"\ar[rr]^f&&N\ars[u]_{u'}="UP"
\ar @{=>}^{\mathcal{F}^{f^{-1}}} "U";"UP"
}
	\]%
		\end{theorem}
		\begin{proof}
It is enough to prove continuity of $f^{-1}$
to prove continuity of $u'$. However this is evident, because $f$
is continuous bijection.

We assume that we defined topology on the set
$\Gamma(\bundle{}pA{}))$ and $\Gamma(\bundle{}qB{}))$
according to remark
\ref{remark: set of sections of bundle, topology}.
Let us consider sections
$u$, $v\in\Gamma(\bundle{}pA{}))$, 
$u'=\mathcal{F}^{f^{-1}}(u)$, $v'=\mathcal{F}^{f^{-1}}(u')$
such that there exists
$W_{L,V}\subset\Gamma(\bundle{}qB{}))$ 
where $L$ is compact subset of space $N$,
$V$ is open subset of space $\mathcal{B}$,
$u'$, $v'\in W_{L,V}$.
Since $f$ is continuous bijection,
$K=f^{-1}(L)$ is compact subset of space $M$.
Since $\mathcal{F}$ is continuous,
$U=\mathcal{F}^{-1}(V)$ is open subset of space $\mathcal{A}$.

According to our design
	\begin{equation}
u'f=\mathcal{F}u
	\label{eq: fibered morphism defines morphism of spaces of sections, 1}
	\end{equation}
From \eqref{eq: fibered morphism defines morphism of spaces of sections, 1} it follows that
	\begin{equation}
\mathcal{F}u(K)=u'f(K)=u'(L)\subset V
	\label{eq: fibered morphism defines morphism of spaces of sections, 2}
	\end{equation}
From \eqref{eq: fibered morphism defines morphism of spaces of sections, 2} it follows that
	\begin{equation}
u(K)\subset \mathcal{F}^{-1}V=U
	\label{eq: fibered morphism defines morphism of spaces of sections, 3}
	\end{equation}
From \eqref{eq: fibered morphism defines morphism of spaces of sections, 3} it follows that
$u\in W_{K,U}$. Similarly we prove that $u'\in W_{K,U}$.
Therefore, for open set $W_{L,V}$ we found open set $W_{K,U}$
such that $\mathcal{F}^{f^{-1}}(W_{K,U})\subset W_{L,V}$.
This proves continuity of map $\mathcal{F}^{f^{-1}}$.
		\end{proof}

Since $f=\mathrm{id}$, then $\mathrm{id}^{-1}=\mathrm{id}$.
In this case
we use notation $\mathcal{F}^{\mathrm{id}}$ for morphism of spaces of sections.
It is evident, that
\[
\mathcal{F}^{\mathrm{id}}(u)=\mathcal{F}u
\]

		\begin{definition}
		\label{definition: fibered subset}
Let $\bundle{\mathcal{A}}aAN$ and $\bundle{\mathcal{B}}bBM$ be bundles.
Suppose fibered morphism $(\begin{array}{cc}
\mathcal{F}:\mathcal{A}\rightarrow\mathcal{B},&
f:M\rightarrow N
\end{array})$ is defined by diagram
	\[
\xymatrix{
\mathcal{A}\arp[d]^{\bundle{}aA{}}\ar[rr]^{\mathcal{F}}
&&\mathcal{B}\arp[d]^{\bundle{}{b}{B}{}}\\
N\ar[rr]^f&&M
}
	\]
where maps $\mathcal{F}$ and $f$ are injections. Then bundle $\bundle{}aA{}$ is called
\AddIndex{fibered subset}{fibered subset} or \AddIndex{subbundle}{subbundle} of $\bundle{}bB{}$.
We also use notation
\symb{\bundle{}{a}{A}{}\subseteq\bundle{}{b}{B}{}}1{fibered subset} or
\symb{\mathcal{A}\subseteq\mathcal{B}}1{subbundle}.

Without loss of generality we assume that $A\subseteq B$, $N\subseteq M$.
		\qed
		\end{definition}

Let us consider bundles $$\bundle{\mathcal{A}}pAM$$
and $$\bundle{\mathcal{B}}qBN$$
\AddIndex{Cartesian power $\mathcal{A}$ of bundle $\mathcal{B}$}{Cartesian power of bundle}
is the set \symb{\bundle{}pA{}^{\bundle{}qB{}}}1{Cartesian power of bundle}
of fibered morphisms $$(\begin{array}{cc}
\mathcal{F}:\mathcal{A}\rightarrow\mathcal{B},&
f:M\rightarrow N
\end{array})$$
At this time I do not see how we can define structure of bundle
in Cartesian power of bundle.
Although for given $m\in M$, $n\in N$
I can consider Cartesian power $B_n{}^{A_m}$.
Based on theorem
\ref{theorem: fibered morphism defines morphism of spaces of sections},
we can define map
\[
f:\bundle{}qB{}^{\bundle{}pA{}}\rightarrow\Gamma(\mathcal{B})^{\Gamma(\mathcal{A})}
\]
%Let us consider subbundles of $\bundle{}qB{}^{\bundle{}pA{}}$ of the form
%\[
%W_{K,U}=\{f:A\rightarrow B|f(K)\subset U\}
%\]
%where $K$ is compact subset of space $A$,
%$U$ is open subset of space $B$.
Let us consider subsets $W_{K,U}\subset\Gamma(\mathcal{B})^{\Gamma(\mathcal{A})}$
where $K$ is compact subset of sections of bundle $\mathcal{A}$,
$U$ is open subset of sections of bundle $\mathcal{B}$.
Sets
$W_{K,U}$
form base of topology
on space $\Gamma(\mathcal{B})^{\Gamma(\mathcal{A})}$.
We choose coarsest topology in $\bundle{}qB{}^{\bundle{}pA{}}$, 
for which mapping $f$ is continuous.

We considered the structure of open set of sections
of bundle $\mathcal{B}$
in remark \ref{remark: set of sections of bundle, topology}.
Since a set of sections of bundle $\mathcal{A}$
is set of mappings, we can
look for theorem similar to Arzel\`a's theorem in calculus
(\citeBib{Kolmogorov Fomin}, p. 54),
to answer the question when this set
is compact. At this time I leave this question open.

According to \citeBib{Maunder: Algebraic Topology}, p. 214,
given spaces $A$, $B$, $C$, $D$ and
mappings $f:A\rightarrow C$, $g:D\rightarrow B$
we define morphism of mapping spaces
\[
g^f:D^C\rightarrow B^A
\]
by law
\begin{align*}
g^f(h)&=fhg&h:&C\rightarrow D&g^f(h):&A\rightarrow B
\end{align*}
Thus, we can represent the morphism of mapping spaces
using diagram
	\[
\xymatrix{Yn j
A\ar[rr]^f\ar[d]_{g^f(h)}&&C\ar[d]^h\\
B&&D\ar[ll]^g
}
	\]

\ifx\texFuture\Defined
Прежде всего возникает вопрос. Что мы называем расслоением? Когда
мы говорим о расслоении, мы подразумеваем непрерывность базы и
непрерывность слоя. Если слой имеет дискретную топологию, то мы
называем соответствующее множество не расслоением, а слоением.
Расслоения и слоения очень сильно отличаются, и потому мы
используем различные методы для их изучения.
 
Тем не менее, если отвлечься от топологии слоя, эти объекты имееют
нечто общее. Кроме того, нужно иметь в виду следующее.
\begin{itemize}
\item Алгебра может быть определена не только на непрерывном множестве.
Конечные и дискретные группы, например, важны при изучении
многих геометрических объектов и физических явлений.
\item Фактор множество топологического пространства может оказаться
конечным множеством, тем самым разрушая топологию.
В то же время совсем не очевидно, почему мы не должны
рассматривать подобного рода отношение эквивалентности.
\end{itemize}

На данном этапе для меня топология слоя несущественна, но важна
непрерывность базы. Принимая во внимание высказанные выше соображения,
я не буду делать различия между расслоением и накрытием
и буду называть и то, и другое расслоением.
\fi

%auto-ignore
\def\texFiberedAlgebra{}
\ifx\PrintBook\Defined
				\chapter{Fibered \texorpdfstring{$\mathfrak{F}$}{F}\Hyph Algebra}
\fi

			\section{Fibered \texorpdfstring{$\mathfrak{F}$}{F}\Hyph Algebra}

	 \begin{definition}
An n-ary
\AddIndex{operation on bundle}{operation on bundle}
$\bundle{}pE{}$
is a fibered morphism
\[
\mathcal{F}:\mathcal{E}^n\rightarrow \mathcal{E}
\]
$n$ is \AddIndex{arity of operation}{arity of operation}.
$0$-arity operation is a section of $\mathcal{E}$.
	 \qed
	 \end{definition}

We can represent the operation using the diagram
	\[
\xymatrix{
\mathcal{E}^n\ar[rr]^\omega\arp@/^2pc/[dd]^p\arp@/_2pc/[dd]_p&&\mathcal{E}\arp[dd]^p\\
\bigodot...\bigodot\ar@{=>}@/^1pc/[rr]^\omega &&\\
M\ar[rr]^{id}&&M
}
	\]

		\begin{theorem}
		\label{theorem: continuity of operation}
Let $U$ be an open set of base $M$.
Suppose there exist trivialization of bundle $\bundle{}pE{}$ over $U$.
Let $x\in M$.
Let $\omega$ be $n$\Hyph ary operation on bundle $\bundle{}pE{}$ and
\[
\omega(p_1,...,p_n)=p
\]
in the fiber $E_x$.
Then there exist open sets $V\subseteq U$,
$W\subseteq E$, $W_1\subseteq E_1$, ..., $W_n\subseteq E_n$
such, that $x\in V$, $p\in W$, $p_1\in W_1$, ..., $p_n\in W_n$,
and for any $x'\in V$, $p'\in W\cap\omega V$ there exist
$p'_1\in W_1$, ..., $p'_n\in W_n$ such, that
\[
\omega(p'_1,...,p'_n)=p'
\]
in the fiber $E_{x'}$.
%Если множества $V$, $W$ связны,
%то также связны множества $W_1$, ..., $W_n$.
		\end{theorem}
		\begin{proof}
According to \citeBib{Bourbaki: General Topology 1}, page 44,
since $V$ belongs to the base
of topology of space $U$ and $W$ belongs to the base
of topology of space $E$, then
set $V\times W$ 
belongs the base of topology of space $\mathcal{E}$.
Similarly, since 
$V$ belongs to the base
of topology of space $U$ and
$W_1$, ..., $W_n$ belong to the base of topology of space $E$,
set $V\times W_1\times ...\times W_n$
belongs the base of topology of space $\mathcal{E}^n$.

Since mapping $\omega$ is continuous, then for an open set
$V\times W$ there exists an open set $S\subseteq \mathcal{E}^n$ such,
that $\omega S\subseteq V\times W$.
Suppose $x'\in V$.
Let
$(x',p')\in\omega S$ be an arbitrary point.
Then there exist such $p'_1\in E_{x'}$, ..., $p'_n\in E_{x'}$,
that
\[
\omega(p'_1,...,p'_n)=p'
\]
in fiber $E_{x'}$.
According to this there exist sets $R$, $R'$
from base of topology of space $U$,
and sets $T_1$, ..., $T_n$,
$T'_1$, ..., $T'_n$ from base of topology of space $E$,
such that $x\in R$, $x'\in R'$,
$p_1\in T_1$, $p'_1\in T'_1$, ..., $p_n\in T_n$, $p'_n\in T'_n$,
$R\times T_1\times ...\times T_n\subseteq S$,
$R'\times T'_1\times ...\times T'_n\subseteq S$.
We proved the theorem since
$W_1=T_1\cup T'_1$, ..., $W_n=T_n\cup T'_n$
are open sets.
%
%According to
%proposition 8
%(\citeBib{Bourbaki: General Topology 1}, page 110)
%if $V$ and $W$ are connected, then $V\times W$ is connected.
		\end{proof}

Theorem \ref{theorem: continuity of operation} tells
about continuity of operation $\omega$, however this
theorem tells nothing regarding sets
$W_1$, ..., $W_n$. In particular, it is possible that these sets
are not connected.

We suppose $W=\{p\}$, $W_1=\{p_1\}$, ..., $W_n=\{p_n\}$,
if topology on fiber $A$ is discrete.
This leads one to assume that  in the neighborhood $V$
the operation does not depend on a fiber.
We call the operation $\omega$ locally constant.
However, it is possible that a condition of constancy
is broken on bundle in general.
Thus  the covering space $R\rightarrow S^1$ of the circle $S^1$
defined by $p(t)=(\sin t,\cos t)$ for any $t\in R$ is
bundle over circle with fiber of group of integers.

Let us consider the alternative point of view on the
continuity of operation $\omega$
to get a better understanding of role of continuity
Let us consider the continuity of operation $\omega$
to better see what does it mean. We need to consider
sections, if we want to show that
infinitesimal change of operand when moving along base
causes infinitesimal change of operation.
This change is legal,
because we defined operation on bundle in fiber.

		\begin{theorem}
		\label{theorem: operation over sections}
An n-ary operation on bundle maps sections into section.
		\end{theorem}
		\begin{proof}
Suppose $f_1$, ..., $f_n$ are sections and we define map
	\begin{equation}
f=\omega^{\mathrm{id}}(f_1,...,f_n):M\rightarrow\mathcal{E}
	\label{eq: operation over sections}
	\end{equation}
as
	\begin{equation}
f(x)=\omega(f_1(x),...,f_n(x))
	\label{eq: operation over sections, x}
	\end{equation}
Let $x\in M$ and $u=f(x)$. Let $U$ be a neighborhood of the point $u$ in the range
of the map $f$.

Since $\omega$ is smooth map,
then according to \citeBib{Bourbaki: General Topology 1}, page 44,
for any $i$,
$1\le i\le n$ the set $U_i$ is defined in the range of section $f_i$
such, that $\prod_{i=1}^nU_i$
is open in the range of section $(f_1,...,f_n)$ of the bundle $\mathcal{E}^n$ and
\[
\omega(\prod_{i=1}^nU_i)\subseteq U
\]

Let $u'\in U$.
Since $f$ is a map, then there exist $x'\in M$ such that $f(x')=u'$.
From equation \eqref{eq: operation over sections, x} it follows
that there exist $u'_i\in U_i$, $p(u'_i)=x'$ such, that
$\omega(u'_1,...,u'_n)=u'$.
Since $f_i$ is a section, then there exist a set $V_i\subseteq M$ such, 
that $f_i(V_i)\subseteq U_i$
and $x\in V_i$, $x'\in V_i$. Therefore, the set
\[
V=\cap_{i=1}^nV_i
\]
is not empty, it is open in $M$ and $x\in V$, $x'\in V$.
Thus the map $f$ is smooth and $f$ is the section.
		\end{proof}

We can represent the operation using the diagram
	
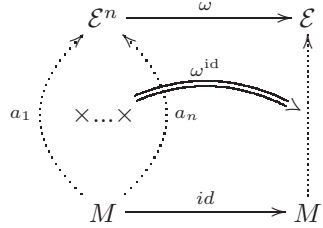
\[
\xymatrix{
\mathcal{E}^n\ar[rr]^\omega&&\mathcal{E}\\
\times...\times\ar@{=>}@/^1pc/[rr]^{\omega^{\mathrm{id}}} &&\\
M\ars@/^2pc/[uu]^{a_1}\ars@/_2pc/[uu]_{a_n}\ar[rr]^{id}&&M\ars[uu]
}
	\]

		\begin{theorem}
		\label{theorem: continuous operation over sections}
$\omega^{\mathrm{id}}$ is continuous on $\Gamma(\mathcal{E})$.
		\end{theorem}
		\begin{proof}
Let us consider a set $W_{K,U}$,
where $K$ is compact set of space $M$,
$U$ is open set of space $\mathcal{E}$.
We can represent set $U$ as $V\times E$,
where $V$ is open set of space $M$, $K\subset V$.
$\omega^{-1}(V\times E)=V\times E^n$ is open set.
Therefore,
\begin{equation}
(\omega^{\mathrm{id}})^{-1}W_{K,V\times E}=W_{K,V\times E^n}
\label{eq: continuous operation over sections}
\end{equation}
From \eqref{eq: continuous operation over sections}
continuity of $\omega^{\mathrm{id}}$ follows.
\end{proof}

	 \begin{definition}
Let $A$ be $\mathfrak{F}$\Hyph algebra
(\citeBib{Burris Sankappanavar}).
We can extend $\mathfrak{F}$\Hyph algebraic structure from fiber
$A$ to bundle
$\bundle{\mathcal{A}}pAM$.
If operation $\omega$ is defined on $\mathfrak{F}$\Hyph algebra $A$
\[
a=\omega(a_1,...,a_n)
\]
then operation $\omega$ is defined on bundle
\[
a(x)=\omega(a_1,...,a_n)(x)=\omega(a_1(x),...,a_n(x))
\]
We say that $\bundle{}pA{}$ is
a \AddIndex{fibered $\mathfrak{F}$\Hyph algebra}{fibered F-algebra}.
	 \qed
	 \end{definition}

Depending on the structure we talk for instance about
\AddIndex{fibered group}{fibered group},
\AddIndex{fibered ring}{fibered ring},
or \AddIndex{vector bundle}{vector bundle}.

Main properties of $\mathfrak{F}$\Hyph algebra
hold for fibered $\mathfrak{F}$\Hyph algebra as well. Proving
appropriate theorems we can refer on this statement.
However in certain cases the proof itself may be of deep interest,
allowing a better view
of the structure of the fibered $\mathfrak{F}$\Hyph algebra.
However properties of $\mathfrak{F}$\Hyph algebra on the set of sections
are different from properties of $\mathfrak{F}$\Hyph algebra in fiber.
For instance, if the product in fiber has inverse element,
it does not mean that the product of sections has inverse element.
Therefore, fibered continuous field generates
ring on the set of sections.
This is the advantage when we consider fibered algebra.
I want also to stress that the operation on bundle is not defined for
elements from different fibers.

Let transition functions $g_{\epsilon\delta}$ determine bundle $\mathcal{B}$
over base $N$.
Let us consider maps $V_\epsilon\in N$ and $V_\delta\in N$,
$V_\epsilon\cap V_\delta\ne\emptyset$.
Point $q\in\mathcal{B}$ has representation $(y,q_\epsilon)$
in map $V_\epsilon$ and representation
$(y,q_\delta)$ in map $V_\delta$.
Therefore,
\[
p_\alpha=f_{\alpha\beta}(p_\beta)
\]
\[
q_\epsilon=g_{\epsilon\delta}(q_\delta)
\]
When we move from map $U_\alpha$ to map $U_\beta$
and from map $V_\epsilon$ to map $V_\delta$,
representation of correspondence changes according to the law
\[
(x,y,p_\alpha,q_\epsilon)=(x,y,f_{\alpha\beta}(p_\beta),g_{\epsilon\delta}(q_\delta))
\]
This is consistent with the transformation
when we move from map $U_\alpha\times V_\epsilon$ to map $U_\beta\times V_\delta$
in the bundle $\mathcal{A}\times\mathcal{B}$.

		\begin{theorem}
		\label{theorem: operation in two neighborhood}
Let transition functions $f_{\alpha\beta}$ determine
fibered $\mathfrak{F}$\Hyph algebra $\bundle{\mathcal{A}}pAM$
over base $M$.
Then transition functions $f_{\alpha\beta}$ are
homomorphisms of $\mathfrak{F}$\Hyph algebra $A$.
		\end{theorem}
		\begin{proof}
Let $U_\alpha\in M$ and $U_\beta\in M$,
$U_\alpha\cap U_\beta\ne\emptyset$ be
neighborhoods where fibered $\mathfrak{F}$\Hyph algebra $\bundle{}pA{}$ is trivial.
Let
	\begin{equation}
a_\beta=f_{\beta\alpha}(a_\alpha)
	\label{eq: operation in two neighborhood}
	\end{equation}
be map from bundle $\bundle{}pA{}|_{U_\alpha}$ into
bundle $\bundle{}pA{}|_{U_\beta}$.
Let $\omega$ be $n$\hyph ary operation and points $e_1$, ..., $e_n$
belong to fiber $A_x$, $x\in U_1\cap U_2$. Suppose
	\begin{equation}
e=\omega(e_1,...,e_n)
	\label{eq: operation in two neighborhood, expression}
	\end{equation}

We represent point $e\in\bundle{}pA{}|_{U_\alpha}$ as
$(x,e_\alpha)$ and point $e_i\bundle{}pA{}|_{U_\alpha}$ as
$(x,e_{i\alpha})$.
We represent point $e\in\bundle{}pA{}|_{U_\beta}$ as
$(x,e_\beta)$ and point $e_i\in\bundle{}pA{}|_{U_\beta}$ as
$(x,e_{i\beta})$. According to \eqref{eq: operation in two neighborhood}
	\begin{equation}
e_\beta=f_{\beta\alpha}(e_\alpha)
	\label{eq: operation in two neighborhood, e}
	\end{equation}
	\begin{equation}
e_{i\beta}=f_{\beta\alpha}(e_{i\alpha})
	\label{eq: operation in two neighborhood, ei}
	\end{equation}
According to \eqref{eq: operation in two neighborhood, expression},
the operation in the bundle $A_x$ over neihgborhood $U_\beta$ is
	\begin{equation}
e_\beta=\omega(e_{1\beta},...,e_{n\beta})
	\label{eq: operation in two neighborhood, fiber}
	\end{equation}
Substituting \eqref{eq: operation in two neighborhood, e},
\eqref{eq: operation in two neighborhood, ei} into
\eqref{eq: operation in two neighborhood, fiber} we get
\[
f_{\beta\alpha}(e_\alpha)=\omega(f_{\beta\alpha}(e_{1\alpha}),...,f_{\beta\alpha}(e_{n\alpha}))
\]
This proves that $f_{\beta\alpha}$ is homomorphism of $\mathfrak{F}$\Hyph algebra.
		\end{proof}

	\begin{definition}
	\label{definition: homomorphism of fibered F-algebras} 
Let $\bundle{\mathcal{A}}{p}{A}M$
and $\bundle{\mathcal{A}'}{p'}{A'}{M'}$
be two fibered $\mathfrak{F}$\Hyph algebras. Fibered morphism
\symb{f}0{homomorphism of fibered F-algebras}
\[
\ShowSymbol{homomorphism of fibered F-algebras}:\mathcal{A}\rightarrow\mathcal{A}'
\]
is called
\AddIndex{homomorphism of fibered $\mathfrak{F}$\Hyph algebra}{homomorphism of fibered F-algebras}
if respective fiber map
\symb{f_x}0{fiber map}
\[
\ShowSymbol{fiber map}:A_x\rightarrow A_{x'}'
\]
is homomorphism of $\mathfrak{F}$\Hyph algebra $A$.
	 \qed
	 \end{definition}

	\begin{definition}
	\label{definition: isomorphism of fibered F-algebras} 
Let $\bundle{\mathcal{A}}{p}{A}M$
and $\bundle{\mathcal{A}'}{p'}{A'}{M'}$
be two fibered $\mathfrak{F}$\Hyph algebras. Homomorphism of fibered $\mathfrak{F}$\Hyph algebras
$f$
is called
\AddIndex{isomorphism of fibered $\mathfrak{F}$\Hyph algebras}{isomorphism of fibered F-algebras}
if respective fiber map
\symb{f_x}0{fiber map}
\[
\ShowSymbol{fiber map}:A_x\rightarrow A_{x'}'
\]
is isomorphism of $\mathfrak{F}$\Hyph algebra $A$.
	 \qed
	 \end{definition}

	\begin{definition}
		\label{definition: fibered subalgebra}
Let $\bundle{\mathcal{A}}pAM$ be
an $\mathfrak{F}$\Hyph fibered F-algebra and $A'$ be $\mathfrak{F}$\Hyph subalgebra
of the $\mathfrak{F}$\Hyph algebra $A$.
An fibered $\mathfrak{F}$\Hyph algebra
$\bundle{\mathcal{A}'}p{A'}M$
is a \AddIndex{fibered $\mathfrak{F}$\Hyph subalgebra}{fibered F-subalgebra}
of the fibered $\mathfrak{F}$\Hyph algebra
$\bundle{}pA{}$ if
homomorphism of fibered $\mathfrak{F}$\Hyph algebras
$\mathcal{A}'\rightarrow\mathcal{A}$
is fiber embedding.
	 \qed
	 \end{definition}

The homomorphism of fibered $\mathfrak{F}$\Hyph algebra is essential part of this definition.
We can break continuity,
if we just limit ourselves to the fact of the existence
of $\mathfrak{F}$\Hyph subalgebra in each fiber.

We defined an operation based reduced Cartesian product of bundles.
Suppose we defined an operation based Cartesian product of bundles.
Then the operation is defined for any elements of the bundle. However,
since $p(a_i)=p(b_i)$, $i=1$, ..., $n$, then $p(\omega(a_1,...,a_n))=p(\omega(b_1,...,b_n))$.
Therefore, the operation is defined between fibers. We can map this operation to base
using projection. This structure is not different from quotient $\mathfrak{F}$\Hyph algebra
and does not create new element in bundle theory. The same time 
mapping between different maps
of bundle and opportunity to define an operation over sections
create problems for this structure. 

			\section{Representation of Fibered \texorpdfstring{$\mathfrak{F}$}{F}\Hyph Algebra}
			\label{section: Representation of Fibered F-Algebra}

	\begin{definition}
	\label{definition: nonsingular transformation, bundle} 
We call the fibered map
\[
\mathcal{T}:\mathcal{E}\rightarrow\mathcal{E}
\]
\AddIndex{transformation of bundle}
{transformation of bundle},
if respective fiber map
\[
t_x:E_x\rightarrow E_x
\]
is transformation of a fiber.
	 \qed
	 \end{definition}

		\begin{theorem}
		\label{theorem: continuity of transformation}
Let $U$ be open set of base $M$ such
that there exists
a local chart of the bundle $\bundle{}pE{}$.
Let $t$ be transformation of bundle $\bundle{}pE{}$.
Let $x\in M$ and $p'=t_x(p)$ in fiber $E_x$.
Then for an open set $V\subseteq M$,
$x\in V$ and for an open set $W'\subseteq E$,
$p'\in W'$
there exists an open set $W\subseteq E$, $p\in W$
such that if $x_1\in V$, $p_1\in W$,
then $p'_1=t_{x_1}(p_1)\in W'$.
		\end{theorem}
		\begin{proof}
According to \citeBib{Bourbaki: General Topology 1}, page 44,
sets $V\times W$, where $V$ forms base
of topology of space $U$ and $W$ forms base
of topology of space $E$,
form base of topology of space $\mathcal{E}$. 

Since map $t$ is continuous, then for open set
$V\times W'$ there exists open set $V\times W$ such,
that $t( V\times W)\subseteq V\times W'$. This is the statement of theorem.
		\end{proof}

		\begin{theorem}
		%\label{theorem: continuity of operation}
Transformation of bundle $\bundle{}pE{}$
maps section into section.
		\end{theorem}
		\begin{proof}
We define the image of section $s$ over transformation $t$
using commutative diagram
	\[
\xymatrix{
\mathcal{E}\ar[rr]^t&&\mathcal{E}\\
&M\ars[ul]^s\ars[ur]_{s'}&
}
	\]
Continuity of map $s'$
follows from theorem
\ref{theorem: continuity of transformation}.
		\end{proof}

	 \begin{definition}
	\label{definition: Tstar transformation of bundle} 
Transformation of bundle is
\AddIndex{left-side transformation}
{left-side transformation of bundle} or
\AddIndex{\Ts transformation of bundle}{Tstar transformation of bundle}
if it acts from left
\[
u'=t u
\]
We denote
\symb{{}^\star \mathcal{E}}1
{set of Tstar nonsingular transformations of bundle} or
\symb{{}^\star \bundle{}pE{}}1
{set of Tstar nonsingular transformations of bundle, projection} or
the set of nonsingular \Ts transformations of bundle
$\bundle{}pE{}$.
	 \qed
	 \end{definition}

	 \begin{definition}
Transformations is
\AddIndex{right-side transformations}
{right-side transformation of bundle} or
\AddIndex{\sT transformation of bundle}{starT transformation of bundle}
if it acts from right
\[
u'= ut
\]
We denote
\symb{\mathcal{E}^\star}1{set of starT nonsingular transformations of bundle} or
\symb{\bundle{}pA{}^\star}1{set of starT nonsingular transformations of bundle, projection}
the set of nonsingular \sT transformations of bundle
$\bundle{}pE{}$.
	 \qed
	 \end{definition}

We denote
\symb{e}1{identical transformation of bundle}
identical transformation of bundle.

		\begin{remark}
		\label{remark: topology on set of transformations of bundle}
Since we define \Ts transformation of bundle by fiber, then set ${}^\star \bundle{}pE{}$
is bundle isomorphic to the bundle $\bundle{}p{{}^\star E}{}$.
Since $^\star E=E^E$,
we can define compact\hyph open topology in fiber.
This gives us an ability to answer on question: how close are transformations,
arising in neighboring fibers.
Let us assume that transformations $t(x)$, $t(x_1)$
are close if there exists open set $W_{K,U}\subset E^E$,
$t(x)\in W_{K,U}$, $t(x_1)\in W_{K,U}$.
%Согласно теореме \ref{theorem: continuity of transformation}
%как только точка 
%$x_1\in M$ оказывается в окрестности точки $x\in M$,
%а точка $p_1\in E_{x_1}$ оказывается в окрестности точки $p\in E_x$
%образ отображения $t_{x_1}(p_1)$ оказывается в окрестности образа отображения $t_x(p)$.
%Однако из теории метрических пространств известно,
%что выбор окрестности $W'$ существенно зависит от выбора точки $p$.
		\qed
		\end{remark}

		\begin{definition}
		\label{definition: Tstar representation of fibered F-algebra} 
Suppose we defined the structure of fibered $\mathfrak{F}$\Hyph algebra
on the set ${}^\star \bundle{}pA{}$
(\citeBib{Burris Sankappanavar}).
Let $\bundle{}pB{}$ be fibered $\mathfrak{F}$\Hyph algebra.
We call homomorphism of fibered $\mathfrak{F}$\Hyph algebras
	\begin{equation}
f:\bundle{}pB{}\rightarrow {}^\star \bundle{}pA{}
	\label{eq: Tstar representation of fibered F-algebra}
	\end{equation}
\AddIndex{left-side representation}
{left-side representation of fibered F-algebra} or
\AddIndex{\Ts representation of fibered $\mathfrak{F}$\Hyph algebra $\bundle{}pB{}$}
{Tstar representation of fibered F-algebra}.
		\qed
		\end{definition}

		\begin{definition}
		\label{definition: starT representation of fibered F-algebra} 
Suppose we defined the structure of fibered $\mathfrak{F}$\Hyph algebra
on the set $\bundle{}pA{}^\star$
(\citeBib{Burris Sankappanavar}).
Let $\bundle{}pB{}$ be fibered $\mathfrak{F}$\Hyph algebra.
We call homomorphism of fibered $\mathfrak{F}$\Hyph algebras
	\[
f:\bundle{}pB{}\rightarrow \bundle{}pA{}^\star
	\]
\AddIndex{right-side representation}
{right-side representation of fibered F-algebra} or
\AddIndex{\sT representation of fibered $\mathfrak{F}$\Hyph algebra $\bundle{}pB{}$}
{starT representation of fibered F-algebra}.
		\qed
		\end{definition}
We extend to bundle representation theory convention
described in remark
\xRef{0701.238}{remark: left and right matrix notation}.
We can write duality principle in the following form

		\begin{theorem}[duality principle]
		\label{theorem: duality principle, fibered F-algebra representation}
Any statement which holds for \Ts representation of fibered $\mathfrak{F}$\Hyph algebra $\bundle{}pA{}$
holds also for \sT representation of fibered $\mathfrak{F}$\Hyph algebra $\bundle{}pA{}$.
		\end{theorem}

There are two ways to define a \Ts representation of $\mathfrak{F}$\Hyph algebra $B$ in the bundle
$\bundle{}pA{}$. We can define or \Ts representation in the fiber,
either define \Ts representation in the set $\Gamma(\bundle{}pA{})$.
In the former case the representation defines the same transformation in
all fibers. In the later case the picture is less restrictive, however we do not have the whole
picture of the diversity of representations in the bundle. Studying the representation of
the fibered $\mathfrak{F}$\Hyph algebra, we point out that representations in different fibers are
independent. Demand of smooth dependence of transformation on fiber
put additional constrains for \Ts representation of fibered $\mathfrak{F}$\Hyph algebra.
The same time this constrain allows learn \Ts representation of the fibered $\mathfrak{F}$\Hyph algebra
when in the fiber there defined $\mathfrak{F}$\Hyph algebra with parameters
(for instance, the structure constants of Lie group) smooth dependent on fiber.

		\begin{remark}
		\label{remark: reduce details on the diagram, fibered F-algebra}
Using diagrams we can express definition
\ref{definition: Tstar representation of fibered F-algebra}
the following way.
	
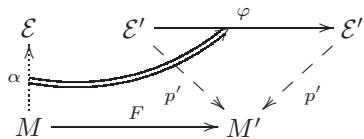
\[
\xymatrix{
\mathcal{E}&\mathcal{E}'\arp[dr]_{p'}\ar[rr]^\varphi="C2"
& & \mathcal{E}'\arp[dl]^{p'}\\
M\ars[u]^\alpha="C1"\ar[rr]^F& & M'\\
\ar @{=>}@/_1pc/ "C1";"C2"
}
	\]
Map $F$ is injection. Because we expect that representation of fibered $\mathfrak{F}$\Hyph algebra
acts in each fiber, then we see that map $F$ is bijection.
Without loss of generality, we assume that $M=M'$ and map $F$ is the identity
map. We tell that we define the representation  of the fibered $\mathfrak{F}$\Hyph algebra
$\bundle{}pB{}$ in the bundle $\bundle{}pA{}$ over the set $M$.
Since we know the base of the bundle, then to reduce details on the diagram we will
describe the representation using the following diagram
	
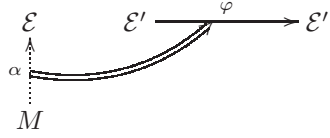
\[
\xymatrix{
\mathcal{E}&\mathcal{E}'\ar[rr]^\varphi="C2"
& & \mathcal{E}'\\
M\ars[u]^\alpha="C1"\\
\ar @{=>}@/_1pc/ "C1";"C2"
}
	\]
		\qed
		\end{remark}

	\begin{theorem}
	\label{theorem: continuous dependene of transformation, fiber}
Let $\mathcal F$ be \Ts representation of fibered $\mathfrak{F}$\Hyph algebra
$\mathcal A$ in bundle $\bundle{}qE{}$. Then for open
set $V\subset M$, $x\in V$ and for open set
$W_{K,U}\subset E^E$, $\mathcal F(x,p)\in W_{K,U}$
there exists open set $W\subset A$, $p\in W$ such
that $x_1\in V$, $p_1\in W$ as soon as
$\mathcal F(x_1,p_1)\in W_{K,U}$.
		\end{theorem}
		\begin{proof}
The statement of theorem is corollary of continuity of map $\mathcal F$
and definition of topology of bundle $\mathcal A$.
		\end{proof}

	\begin{theorem}
	\label{theorem: continuous dependene of transformation, section}
Let $\mathcal F$ be \Ts representation of fibered $\mathfrak{F}$\Hyph algebra
$\mathcal A$ in bundle $\bundle{}qE{}$. Let $a$ be section of bundle $\mathcal A$.
Then for open set
$V\subset M$, $x\in V$ and for open set
$W_{K,U}\subset E^E$, $\mathcal F^\mathrm{id}(a)(x)\in W_{K,U}$
there exists open set $W\subset A$, $a(x)\in W$ such
that $x_1\in V$, $a(x_1)\in W$ as soon as
$\mathcal F^\mathrm{id}(a)(x_1)\in W_{K,U}$.
		\end{theorem}
		\begin{proof}
This is corollary of theorem \ref{theorem: continuous dependene of transformation, fiber}.
		\end{proof}

	 \begin{definition}
	 \label{definition: effective representation of fibered F-algebra}
Suppose map \eqref{eq: Tstar representation of fibered F-algebra} is
an isomorphism of the fibered $\mathfrak{F}$\Hyph algebra $\bundle{}pB{}$ into ${}^\star \bundle{}pA{}$.
Then the \Ts representation of the fibered $\mathfrak{F}$\Hyph algebra $\bundle{}pB{}$ is called
\AddIndex{effective}{effective representation of fibered F-algebra}.
	 \qed
	 \end{definition}

		\begin{remark}
		\label{remark: notation for effective representation of fibered F-algebra}
Suppose the \Ts representation of fibered $\mathfrak{F}$\Hyph algebra
$\mathcal A$ is effective. Then we identify
a section of fibered $\mathfrak{F}$\Hyph algebra
$\mathcal A$ and its image and write \Ts transformation
caused by section $a\in\Gamma(\mathcal A)$
as
\[v'=av\]
Suppose the \sT representation of $\mathfrak{F}$\Hyph algebra
$\mathcal A$ is effective. Then we identify
a section of fibered $\mathfrak{F}$\Hyph algebra
$\mathcal A$ and its image and write \sT transformation
caused by section $a\in\Gamma(\mathcal A)$
as
\[v'=va\]
		\qed
		\end{remark}

	 \begin{definition}
	 \label{definition: transitive representation of fibered F-algebra}
We call a \Ts representation $\mathcal{F}$ of fibered $\mathfrak{F}$\Hyph algebra $\bundle{}pA{}$
\AddIndex{transitive}{transitive representation of fibered F-algebra}
if \Ts representation $F_x$ of $\mathfrak{F}$\Hyph algebra $A_x$
is transitive for any $x$.
We call a \Ts representation of fibered $\mathfrak{F}$\Hyph algebra
\AddIndex{single transitive}{single transitive representation of fibered F-algebra}
if it is transitive and effective.
	 \qed
	 \end{definition}

	 \begin{theorem}
	 \label{theorem: transitive representation of fibered F-algebra}
Let set $E$ be locally compact.
A \Ts representation $\mathcal{F}$ of fibered $\mathfrak{F}$\Hyph algebra
\[\bundle{\mathcal{A}}rAM\] in the bundle \[\bundle{\mathcal{E}}pEM\]
is transitive
if for any sections $a, b \in\Gamma(\mathcal{E})$ exists such section
$g\in\Gamma(\mathcal{A})$ that
\[b=\mathcal{F}^\mathrm{id}(g)a\]
	 \end{theorem}
		\begin{proof}
Let $a, b \in\Gamma(\mathcal{E})$ be sections. In fiber $E(x)$, these sections
define elements $a(x)$, $b(x) \in E(x)$. According to definition
\ref{definition: transitive representation of fibered F-algebra}
$g(x)\in A(x)$ is defined such that
\[b(x)=\mathcal{F}^\mathrm{id}(g(x))a(x)\]
%Наша задача - доказать, что отображения $g$ непрерывно как это
%описано в теореме \ref{theorem: continuity of transformation}.

Let $U_M$ be open set of base $M$ such
that there exists
a local chart of the bundle $\bundle{}pE{}$,
$x\in U_M$.
Let $W'\subseteq E$ is an open set,
$b(x)\in W'$.
Then according to theorem \ref{theorem: continuity of transformation}
there exists an open set $W\subseteq E$, $a(x)\in W$
such that if $x_1\in U_M$, $a(x_1)\in W$,
then $b(x_1)=\mathcal{F}^\mathrm{id}(g)(x_1)(a(x_1))\in W'$.

If closure $\overline W$ is compact, then let us assume $K=\overline W$.
Assume that closure $\overline W$ is not compact.
Then there exists open set $W(x)$, $a(x)\in W(x)$ such that
$\overline W(x)$ is compact.
Similarly there exists open set $W(x_1)$, $a(x_1)\in W(x_1)$ such that
$\overline W(x_1)$ is compact.
Set $V=W\cap(W(x)\cup W(x_1))$ is open, $a(x)\in V$, $a(x_1)\in V$,
closure $\overline V$ is compact. Assume that $K=\overline V$.

Let $U(x)$ be neighborhood of the set $g(x)K$.
Let $U(x_1)$ be neighborhood of the set $g(x_1)K$.
Assume that $U=W'\cup U(x)\cup U(x_1)$.
Then $W_{K,U}\subset E^E$ is open set,
$\mathcal{F}^\mathrm{id}(g)(x)\in W_{K,U}$,
$\mathcal{F}^\mathrm{id}(g)(x_1)\in W_{K,U}$.

According to theorem \ref{theorem: continuous dependene of transformation, section}
there exists open set $S\subset A$,
$g(x)\in S$, $g(x_1)\in S$.
		\end{proof}

		\begin{theorem}	%single transitive
Let set $E$ be locally compact.
A \Ts representation $\mathcal{F}$ of fibered $\mathfrak{F}$\Hyph algebra
\[\bundle{\mathcal{A}}rAM\] in the bundle \[\bundle{\mathcal{E}}pEM\]
is single transitive
if and only if for any sections $a, b \in\Gamma(\mathcal{E})$ exists one and only one section
$g\in\Gamma(\mathcal{A})$ such that
\[b=\mathcal{F}^\mathrm{id}(g)a\]
		\end{theorem}
		\begin{proof}
Corollary of definitions \ref{definition: effective representation of fibered F-algebra},
\ref{definition: transitive representation of fibered F-algebra}
and theorem \ref{theorem: transitive representation of fibered F-algebra}.
		\end{proof}

%auto-ignore
\def\texFiberedMorphism{}
\ifx\PrintBook\Defined
				\chapter{Fibered Morphism}
\fi

			\section{Fibered Morphism}

		\begin{theorem}
		\label{theorem: quotient bundle}
Let us consider fibered equivalence $\bundle{\mathcal{S}}sSM$
on the bundle $\bundle{\mathcal{E}}pEM$.
Then there exists bundle
\[\bundle{\mathcal{E}/\mathcal{S}}t{E/S}M\]
called \AddIndex{quotient bundle}{quotient bundle}
of bundle $\mathcal{E}$ by
the equivalence $\mathcal{S}$.
Fibered morphism
\[
\mathrm{nat}\mathcal{S}:\mathcal{E}\rightarrow
\mathcal{E}/\mathcal{S}
\]
is called \AddIndex{fibered natural morphism}{fibered natural morphism}
or \AddIndex{fibered identification morphism}{fibered identification morphism}.
		\end{theorem}
		\begin{proof}
Let us consider the commutative diagram
\begin{equation}
\label{eq: quotient bundle}
\xymatrix{
\mathcal{E}\ar[rr]^{\mathrm{nat}\mathcal{S}}\arp[rd]_{\bundle{}pE{}}
&&\mathcal{E}/\mathcal{S}\arp[dl]^{\bundle{}t{E/S}{}}\\
&M&
}
\end{equation}
We introduce in $\mathcal{E}/\mathcal{S}$ quotient topology
(\citeBib{Bourbaki: General Topology 1}, page 33),
demanding continuity of mapping $\mathrm{nat}\mathcal{S}$.
According to proposition \citeBib{Bourbaki: General Topology 1}-I.3.6
mapping $\bundle{}t{E/S}{}$ is continuous.

Because we defined equivalence $S$ only between points of
the same fiber $E$, equivalence classes belong to
the same fiber $E/S$ (compare with the remark to proposition
\citeBib{Bourbaki: General Topology 1}-I.3.6).
		\end{proof}

Let $f:\mathcal{A}\rightarrow\mathcal{B}$ be
fibered morphism, base of which is
identity mapping.
According to definition
\xRef{0707.2246}{definition: inverse reduced fibered correspondence}
there exists inverse reduced fibered correspondence $f^{-1}$.
According to theorems
\xRef{0707.2246}{theorem: composition of reduced fibered correspondences}
and \xRef{0707.2246}{theorem: 2 ary fibered relation}
$f^{-1}\circ f$ is $2$\Hyph ary fibered relation.

		\begin{theorem}
		\label{theorem: decomposition of fibered morphism}
Fibered relation $\mathcal{S}=f^{-1}\circ f$ is
fibered equivalence on the bundle $\mathcal{A}$.
There exists decomposition of fibered morphism $f$ into product of
fibered morphisms
	\[
f=itj
	\]
\[
\xymatrix{
\mathcal{A}/\mathcal{S}\ar[rr]^t&&f(\mathcal{A})\ar[d]^i\\
\mathcal{A}\ar[u]^j\ar[rr]^f&&\mathcal{B}
}
\]
$j=\mathrm{nat}\mathcal{S}$ is the natural homomorphism
	\[
\textcolor{red}{j(a)}=j(a)
	\]
$t$ is isomorphism
	\[
\textcolor{red}{r(a)}=t(\textcolor{red}{j(a)})
	\]
$i $ is the inclusion mapping
	\[
r(a)=i(\textcolor{red}{r(a)})
	\]
		\end{theorem}
		\begin{proof}
We verify the statement of theorem in fiber. We need also to check
that equivalence depends continuously on fiber.
		\end{proof}

\section{Free \texorpdfstring{$T\star$}{T*}-Representation of Fibered Group}

Mapping $\mathrm{nat}\mathcal{S}$ does not create bundle,
because different equivalence classes are not homeomorphic
in general.
However the proof of theorem \ref{theorem: quotient bundle}
suggests to the construction which reminds the construction
designed in \citeBib{Postnikov: Differential Geometry}, pages 16 - 17.

	 \begin{definition}
	 \label{definition: fibered little group}
Consider \Ts representation $f$ of fibered group $\bundle{}pG{}$
in bundle $\mathcal{M}$.
A \AddIndex{fibered little group}{fibered little group} or
\AddIndex{fibered stability group}{fibered stability group}
of $h\in \Gamma(\mathcal{M})$ is the set
\symb{\mathcal{G}_h}0{fibered little group}
\symb{\mathcal{G}_h}0{fibered stability group}
\[
\ShowSymbol{fibered little group}=\{g\in \Gamma(\mathcal{G}):f(g)h=h\}
\]
	 \qed
	 \end{definition}

	 \begin{definition}
	 \label{definition: free representation of fibered group}
\Ts representation $f$ of group $G$ is said to be
\AddIndex{free}{free representation of fibered group},
if for any $x\in M$ stability group $G_x=\{e\}$.
	 \qed
	 \end{definition}

		\begin{theorem}	%
		\label{theorem: free representation of group}
Given free \Ts representation $f$ of group $G$ in the set $A$,
there exist $1-1$ correspondence between orbits of representation,
as well between orbit of representation and group $G$.
		\end{theorem}
		\begin{proof}
		\end{proof}

Let us consider covariant free \Ts representation $f$
of fibered group $\bundle{}pG{}$ in fiber $\bundle{}pE{}$.
This \Ts representation determines fibered
equivalence $\mathcal{S}$ on $\bundle{}aE{}$,
$(p,q)\in S$ when $p$ and $q$ belong to
common orbit.
Since the representation in every fiber is free,
all equivalence classes are homeomorphic to group $G$.
Therefore, the mapping $\mathrm{nat}\mathcal{S}$
is projection of the bundle
$\bundle{\mathcal{E}}{\mathrm{nat}\mathcal{S}}G{\mathcal{E}/\mathcal{S}}$.
We also use notation $\mathcal{S}=\mathcal{G}\star$.
We may represent diagram \eqref{eq: quotient bundle}
in the following form
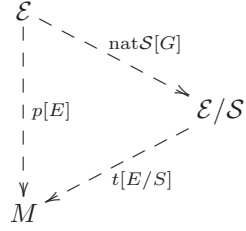
\[
\xymatrix{
\mathcal{E}\arp[drr]^{\bundle{}{\mathrm{nat}\mathcal{S}}G{}}
\arp[dd]^{\bundle{}pE{}}&&\\
&&\mathcal{E}/\mathcal{S}\arp[dll]^{\bundle{}t{E/S}{}}\\
M&&
}
\]
Bundle $\bundle{}{\mathrm{nat}\mathcal{S}}G{}$ is called
\AddIndex{bundle of level $2$}{bundle of level 2}.

\begin{example}
Let us consider the representation of rotation group $SO(2)$ in $R^2$.
All points except the point
$(0,0)$ have trivial little group.
Hence, we defined free representation of group $SO(2)$
in set $R^2\setminus \{(0,0)\}$.

We cannot use this idea in case of bundle
$\bundle{}p{R^2}{}$ and representation of fibered group $\bundle{}t{SO(2)}{}$.
Let $\mathcal{S}$ be relation of fibered equivalence.
The bundle $\bundle{}p{R^2\setminus \{(0,0)\}}{}/\bundle{}t{SO(2)}{}\star$
is not complete. As a consequence passage to the limit may bring into
non\Hyph existent fiber. Therefore we prefer to consider
bundle $\bundle{}p{R^2}{}/\bundle{}t{SO(2)}{}\star$,
keeping in mind, that fiber over point $(x,0,0)$ is degenerate.
	 \qed
\end{example}

We simplify the notation and represent this construction as
\symb{p[E_2,E_1]}0{bundle of level 2}
\[
\xymatrix{
\ShowSymbol{bundle of level 2}:\mathcal{E}_2\arp[r]&
\mathcal{E}_1\arp[r]&
M
}
\]
where we consider bundles
\begin{align*}
\bundle{\mathcal{E}_2}{p_2}{E_2}{\mathcal{E}_1}&&
\bundle{\mathcal{E}_1}{p_1}{E_1}M&
\end{align*}
Similarly we consider
\AddIndex{bundle of level $n$}{bundle of level n}
\symb{p[E_n,...,E_1]}0{bundle of level n}
\begin{equation}
\label{eq: bundle of level n}
\xymatrix{
\ShowSymbol{bundle of level n}:\mathcal{E}_n\arp[r]&\ar@{}[r]_{...}&\arp[r]&
\mathcal{E}_1\arp[r]&
M
}
\end{equation}
The sequence of bundles \eqref{eq: bundle of level n}
is called \AddIndex{tower of bundles}{tower of bundles}.
I made this definition by analogy with Postnikov tower
(\citeBib{Hatcher: Algebraic Topology}).
Postnikov tower is the tower of bundles. Fiber of bundle of level $n$ is
homotopy group of order $n$.
Such definitions are well known, however I gave definition of tower of
bundles, because it follows in a natural way from the text above.

One more example of tower of bundles
attracted my attention (\citeBib{geometry of differential equations},
\citeBib{cohomological analysis}, chapter 2).
We consider the set $J^0(n,m)$ of $0$\Hyph jets
of functions from $R^n$ to $R^m$ as base. We consider 
the set $J^p(n,m)$ of $p$\Hyph jets
of functions from $R^n$ to $R^m$ as bundle of level $p$.

\section{Morphism of \texorpdfstring{$T\star$}{T*}-Representations of Fibered
\texorpdfstring{$\mathfrak{F}$\Hyph }{F-}Algebra}

		\begin{definition}
		\label{definition: morphism of representations of fibered F algebra}
Let us consider \Ts representation
\[
\mathcal{F}:\mathcal{A}\rightarrow {}^\star\mathcal{P}
\]
of fibered $\mathfrak{F}$\Hyph algebra $\bundle{\mathcal{A}}aAM$ in bundle $\bundle{\mathcal{P}}pPM$
and \Ts representation
\[
\mathcal{G}:\bundle{}bB{}\rightarrow {}^\star \bundle{}qQ{}
\]
of fibered $\mathfrak{F}$\Hyph algebra $\bundle{}bB{}$ in bundle $\bundle{}qQ{}$.
Tuple of mappings
	\begin{equation}
(\begin{array}{cc}
\mathcal{C}:\mathcal{A}\rightarrow\mathcal{B},&
\mathcal{R}:\mathcal{P}\rightarrow\mathcal{Q}
\end{array})
	\label{eq: morphism of representations of fibered F algebra, definition, 1}
	\end{equation}
such, that $\mathcal{C}$ is homomorphism of fibered $\mathfrak{F}$\Hyph algebra and
	\begin{equation}
\textcolor{blue}{\mathcal{R}(\mathcal{F}(a)m)}=\mathcal{G}(\textcolor{red}{\mathcal{C}(a)})\textcolor{blue}{\mathcal{R}(m)}
	\label{eq: morphism of representations of fibered F algebra, definition, 2}
	\end{equation}
is called
\AddIndex{morphism of fibered \Ts representations from $\mathcal{F}$ into $\mathcal{G}$}
{morphism of fibered representations from f into g}.
We also say that
\AddIndex{morphism of \Ts representations of fibered $\mathfrak{F}$\Hyph algebra}
{morphism of representations of fibered F algebra} is defined.
		\qed
		\end{definition}

We represent
morphism of \Ts representations of fibered $\mathfrak{F}$\Hyph algebra
as diagram
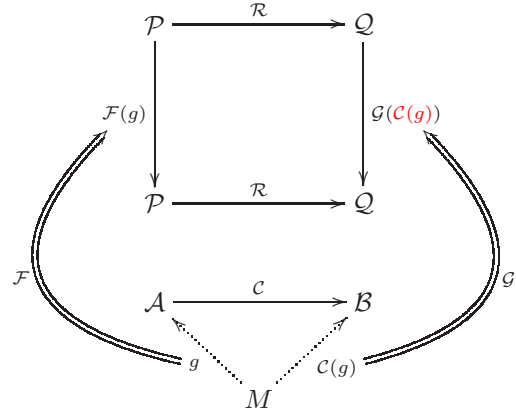
\[
\xymatrix{
\mathcal P\ar[dd]_{\mathcal{F}(g)}="D1"\ar[rr]^{\mathcal R}
&&\mathcal Q\ar[dd]^{\mathcal{G}(\textcolor{red}{\mathcal C(g)})}="D2"\\
\\
\mathcal P\ar[rr]^{\mathcal R}&&\mathcal Q\\
\mathcal{A}\ar[rr]^{\mathcal C}&&\mathcal B\\
&M\ars[ul]^g="C1"\ars[ur]_{\mathcal C(g)}="C2"\\
\ar @{=>}@/^4pc/^{\mathcal{F}} "C1";"D1"
\ar @{=>}@/_4pc/_{\mathcal{G}} "C2";"D2"
}
	\]
Therefore morphism of \Ts representations of fibered $\mathfrak{F}$\Hyph algebra
is morphism of \Ts representations of $\mathfrak{F}$\Hyph algebra in fiber.

		\begin{theorem}
	\label{theorem: morphism of representations of fibered F algebra}
Given single transitive \Ts representation
\[
\mathcal{F}:\mathcal{A}\rightarrow {}^\star\mathcal{P}
\]
of fibered $\mathfrak{F}$\Hyph algebra $\bundle{\mathcal{A}}aAM$ in bundle $\bundle{\mathcal{P}}pPM$
and single transitive \Ts representation
\[
\mathcal{G}:\bundle{}bB{}\rightarrow {}^\star \bundle{}qQ{}
\]
of fibered $\mathfrak{F}$\Hyph algebra $\bundle{}bB{}$ in bundle $\bundle{}qQ{}$,
there exists morphism
\[
(\begin{array}{cc}
\mathcal{C}:\mathcal{A}\rightarrow\mathcal{B},&
\mathcal{R}:\mathcal{P}\rightarrow\mathcal{Q}
\end{array})
\]
of fibered \Ts representations from $\mathcal{F}$ into $\mathcal{G}$.
		\end{theorem}
		\begin{proof}
Corollary of theorem \ref{theorem: morphism of representations of F algebra}
and definition \ref{definition: morphism of representations of fibered F algebra}.
		\end{proof}

		\begin{theorem}
		\label{theorem: product of morphisms of representations of fibered F algebra}
Let
\[
\mathcal F:\mathcal A\rightarrow {}^\star \mathcal M
\]
be \Ts representation of fibered $\mathfrak{F}$\Hyph algebra $\mathcal A$,
\[
\mathcal G:\mathcal B\rightarrow {}^\star \mathcal N
\]
be \Ts representation of fibered $\mathfrak{F}$\Hyph algebra $\mathcal B$,
\[
\mathcal H:\mathcal C\rightarrow {}^\star \mathcal L
\]
be \Ts representation of fibered $\mathfrak{F}$\Hyph algebra $\mathcal C$.
Given morphisms of \Ts representations of fibered $\mathfrak{F}$\Hyph algebra
\[
(\begin{array}{cc}
\mathcal U:\mathcal A\rightarrow\mathcal B,&
\mathcal P:\mathcal M\rightarrow\mathcal N
\end{array})
\]
\[
(\begin{array}{cc}
\mathcal V:\mathcal B\rightarrow\mathcal C,&
\mathcal Q:\mathcal N\rightarrow\mathcal L
\end{array})
\]
There exists morphism of \Ts representations of $\mathfrak{F}$\Hyph algebra
\[
(\begin{array}{cc}
\mathcal W:\mathcal A\rightarrow\mathcal C,&
\mathcal R:\mathcal M\rightarrow\mathcal L
\end{array})
\]
where $\mathcal W=\mathcal U\mathcal V$, $\mathcal R=\mathcal P\mathcal Q$.
We call morphism $(\mathcal W,\mathcal R)$ of
fibered \Ts representations from $\mathcal F$ into $\mathcal H$
\AddIndex{product of morphisms $(\mathcal U,\mathcal P)$ and $(\mathcal V,\mathcal Q)$
of \Ts representations of fibered $\mathfrak{F}$\Hyph algebra}
{product of morphisms of representations of fibered F algebra}.
		\end{theorem}
		\begin{proof}
The mapping $\mathcal W$ is homomorphism of fibered $\mathfrak{F}$\Hyph algebra $\mathcal A$ into
fibered $\mathfrak{F}$\Hyph algebra $\mathcal C$.
We need to show that tuple of mappings $(\mathcal W,\mathcal R)$ satisfies
\eqref{eq: morphism of representations of F algebra, definition, 2}:
\begin{align*}
\textcolor{blue}{\mathcal R(\mathcal F(a)m)}&=\textcolor{blue}{\mathcal Q\mathcal P(\mathcal F(a)m)}\\
&=\textcolor{blue}{\mathcal Q(\mathcal G(\textcolor{red}{\mathcal U(a)})\textcolor{blue}{\mathcal P(m)})}\\
&=\mathcal H(\textcolor{red}{\mathcal V\mathcal U(a)})\textcolor{blue}{\mathcal Q\mathcal P(m)})\\
&=\mathcal H(\textcolor{red}{\mathcal W(a)})\textcolor{blue}{R(m)}
\end{align*}
		\end{proof}

		\begin{theorem}
		\label{theorem: decompositions of morphism of fibered representations}
Let
\[
\mathcal F:\mathcal A\rightarrow {}^\star \mathcal D
\]
be \Ts representation of fibered $\mathfrak{F}$\Hyph algebra $\mathcal A$,
\[
\mathcal G:\mathcal B\rightarrow {}^\star \mathcal E
\]
be \Ts representation of fibered $\mathfrak{F}$\Hyph algebra $\mathcal B$.
Let
\[
(\begin{array}{cc}
\mathcal R_1:\mathcal A\rightarrow\mathcal B,&
\mathcal R_2:\mathcal D\rightarrow\mathcal E
\end{array})
\]
be morphism of fibered representations from $\mathcal F$ into $\mathcal G$.
Suppose
\begin{align*}
\mathcal S_1&=\mathcal R_1\mathcal R_1^{-1}&\mathcal S_2&=\mathcal R_2\mathcal R_2^{-1}
\end{align*}
Then there exist decompositions of $\mathcal R_1$ and $\mathcal R_2$,
which we describe using diagram
\[
\xymatrix{
&&&\mathcal D/\mathcal S_2\ar[rrr]^{\mathcal T_2}&&
&\mathcal R_2\mathcal D\ar[ddddd]^{\mathcal I_2}\\
\mathcal A/\mathcal S_1\ar[rr]^(.3){\mathcal T_1}&&\mathcal R_1\mathcal A\ar[d]^{\mathcal I_1}&&&&\\
\mathcal A\ar[rr]_(.6){\mathcal R_1}\ar[u]^{\mathcal J_1}\ar@{}[urr]_(.7){(1)}&&\mathcal B&&
\mathcal D/\mathcal S_2\ar[r]^{\mathcal T_2}\ar[luu]_{\mathcal{F}_1(\textcolor{red}{\mathcal J_1a})}="E2"&
\mathcal R_2\mathcal D\ar[d]^{\mathcal I_2}\ar[ruu]^{\mathcal{G}_1(\textcolor{red}{\mathcal R_1a})}="E3"\\
&
&&&
\mathcal D\ar[r]_{\mathcal R_2}\ar[u]^{\mathcal J_2}\ar@{}[ur]|{(2)}\ar[ddl]^{\mathcal{F}(a)}="E1"&
\mathcal E\ar[ddr]_{\mathcal{G}(\textcolor{red}{\mathcal R_1a})}="E4"\\
&&&&&&&\\
&M\ars@/^2pc/[uuul]^a="B1"\ars[uuuul]^{\mathcal J_1a}="B2"\ars[uuur]_{\mathcal R_1a}="B4"\ars[uuuur]^{\mathcal R_1a}="B3"
&&\mathcal D\ar[uuuuu]^{\mathcal J_2}\ar[rrr]_{\mathcal R_2}&&&\mathcal E\\
\ar @{=>}@/_1pc/^{\mathcal{F}} "B1";"E1"
\ar @{=>}@/^5pc/^{\mathcal{F}_1} "B2";"E2"
\ar @{=>}^{\mathcal{G}_1} "B3";"E3"
\ar @{=>}@/^1pc/^{\mathcal{G}} "B4";"E4"
}
\]
\begin{itemize}
\item $s=\mathrm{ker}\ \mathcal R_1$ is a congruence on $\mathcal A$.
There exists decompositions of homomorphism $\mathcal R_1$
	\begin{equation}
\mathcal R_1=\mathcal I_1\mathcal T_1\mathcal J_1
	\label{eq: morphism of representations of fibered algebra, homomorphism, 1}
	\end{equation}
$\mathcal J_1=\mathrm{nat}\ s$ is the natural homomorphism
	\begin{equation}
\textcolor{red}{\mathcal J_1(a)}=\mathcal J_1(a)
	\label{eq: morphism of representations of fibered algebra, homomorphism, 2}
	\end{equation}
$\mathcal T_1$ is isomorphism
	\begin{equation}
\textcolor{red}{\mathcal R_1(a)}=\mathcal T_1(\textcolor{red}{\mathcal J_1(a)})
	\label{eq: morphism of representations of fibered algebra, homomorphism, 3}
	\end{equation}
$\mathcal I_1 $ is the inclusion mapping
	\begin{equation}
\mathcal R_1(a)=\mathcal I_1(\textcolor{red}{\mathcal R_1(a)})
	\label{eq: morphism of representations of fibered algebra, homomorphism, 4}
	\end{equation}
\item $\mathcal S_2=\mathrm{ker}\ \mathcal R_2$ is an equivalence on $\mathcal D$.
There exists decompositions of homomorphism $\mathcal R_2$
	\begin{equation}
\mathcal R_2=\mathcal I_2\mathcal T_2\mathcal J_2
	\label{eq: morphism of representations of fibered algebra, map, 1}
	\end{equation}
$\mathcal J_2=\mathrm{nat}\ \mathcal S_2$ is surjection
	\begin{equation}
\textcolor{blue}{\mathcal J_2(m)}=\mathcal J_2(m)
	\label{eq: morphism of representations of fibered algebra, map, 2}
	\end{equation}
$\mathcal T_2$ is bijection
	\begin{equation}
\textcolor{blue}{\mathcal R_2(m)}=\mathcal T_2(\textcolor{blue}{\mathcal J_2(m)})
	\label{eq: morphism of representations of fibered algebra, map, 3}
	\end{equation}
$\mathcal I_2$ is the inclusion mapping
	\begin{equation}
\mathcal R_2(m)=\mathcal I_2(\textcolor{blue}{\mathcal R_2(m)})
	\label{eq: morphism of representations of fibered algebra, map, 4}
	\end{equation}
\item $\mathcal{F}_1$ is \Ts representation of $\mathfrak{F}$\Hyph algebra $\mathcal A/\mathcal S_1$ in $\mathcal D/\mathcal S_2$
\item $\mathcal G_1$ is \Ts representation of $\mathfrak{F}$\Hyph algebra $\mathcal R_1\mathcal A$ in $\mathcal R_2\mathcal D$
\item There exists decompositions of morphism of representations
\[
(\mathcal R_1,\mathcal R_2)=(\mathcal I_1,\mathcal I_2)(\mathcal T_1,\mathcal T_2)(\mathcal J_1,\mathcal J_2)
\]
\end{itemize}
		\end{theorem}
		\begin{proof}
From theorem
\ref{theorem: decomposition of fibered morphism}
it follows the validity of diagrams $(1)$, $(2)$.
We check the statement of theorem in fiber and it follows from theorem
\ref{theorem: decompositions of morphism of representations}.
		\end{proof}

%auto-ignore
\def\texVectorField{}
\ifx\PrintBook\Defined
				\chapter{Vector Bundle over Skew-Field}
				\label{chapter: Vector Bundle Over Skew-Field}
\fi

			\section{Vector Bundle over Skew-Field}
			\label{section: Vector Bundle Over Skew-Field}

To define \Ts representation
	\[
\mathcal F:\mathcal R\rightarrow {}^\star \mathcal{V}
	\]
of fibered ring $\mathcal R$ on the bundle $\mathcal{V}$
we need to define the structure of the fibered ring on the bundle
${}^\star \mathcal{V}$.

		\begin{theorem}
		\label{theorem: Ts representation of the ring}
\Ts representation $\mathcal F$ of the fibered ring
$\mathcal R$ on the bundle $\mathcal{V}$ is defined iff
\Ts representations of fibered multiplicative and additive groups
of the fibered ring $\mathcal R$ are defined
and these fibered \Ts representations hold relationship
\[
\mathcal F(a(b+c))=\mathcal F(a)\mathcal F(b)+\mathcal F(a)\mathcal F(c)
\]
		\end{theorem}
		\begin{proof}
Theorem follows from definition
\ref{definition: Tstar representation of fibered F-algebra}.
		\end{proof}

		\begin{theorem}
		\label{theorem: effective Ts representation of the fibered skew field}
\Ts representation of the fibered skew field $\mathcal D$ is
\AddIndex{effective}{effective representation of fibered skew field}
iff \Ts representation of its fibered multiplicative group is effective.
		\end{theorem}
		\begin{proof}
According to definitions \ref{definition: effective representation of fibered F-algebra}
and \ref{definition: isomorphism of fibered F-algebras} we check the statement of theorem
in fiber. The statement of theorem in fiber is corollary of theorem
\xRef{0701.238}{theorem: effective Ts representation of the skew field}.
		\end{proof}

According to the remark
\ref{remark: notation for effective representation of fibered F-algebra},
since the representation of the fibered skew field is effective, we identify
a section of the fibered skew field and
\Ts transformation corresponding to this section.

		\begin{definition}
		\label{definition: Tstar vector field}
$\mathcal{V}$ is a \AddIndex{$\mathcal D\star$\Hyph vector bundle}{Dstar vector bundle}
over a fibered skew field $\mathcal D$
if $\mathcal{V}$ is a fibered Abelian group and
there exists effective \Ts representation of fibered skew field $\mathcal D$.
Section of $\mathcal D\star$\Hyph vector bundle
is called a \AddIndex{$\mathcal D\star$\Hyph vector field}{Dstar vector field}.
		\qed
		\end{definition} 

		\begin{theorem}
Following conditions hold for $\mathcal D\star$\Hyph vector fields:
\begin{itemize}
\item 
\AddIndex{associative law}{associative law, Dstar vector fields}
	\begin{equation}
	\label{eq: associative law, Dstar vector fields}
(ab)\Vector{m}=a(b\Vector{m})
	\end{equation}
\item 
\AddIndex{distributive law}{distributive law, Dstar vector fields}
	\begin{align}
	\label{eq: distributive law, Dstar vector fields, 1}
a(\Vector{m}+\Vector{n})&=a\Vector{m}+a\Vector{n}\\
	\label{eq: distributive law, Dstar vector fields, 2}
(a+b)\Vector{m}&=a\Vector{m}+b\Vector{m}
	\end{align}
\item 
\AddIndex{unitarity law}{unitarity law, Dstar vector fields}
	\begin{equation}
	\label{eq: unitarity law, Dstar vector fields}
1\Vector{m}=\Vector{m}
	\end{equation}
\end{itemize}
for any $a, b \in \Gamma(\mathcal D)$, $\Vector{m}, \Vector{n} \in \Gamma(\mathcal{V})$.
We call \Ts representation
\AddIndex{$\mathcal D\star$\Hyph product of vector field over scalar}
{Dstar product of vector field over scalar, vector space}.
		\end{theorem}
		\begin{proof}
We check the statement of theorem in fiber and
this statement is corollary of theorem
\xRef{0701.238}{theorem: definition of vector space}.
		\end{proof}

		\begin{definition}
Let $\mathcal{V}$ be a $\mathcal D\star$\hyph vector bundle
over a fibered skew field $\mathcal D$.
Subbundle $\mathcal{N}\subset\mathcal{V}$ is a
\AddIndex{subbundle of $\mathcal D\star$\hyph vector space}
{subbundle of Dstar vector bundle} $\mathcal{V}$
if $\Vector{a}+\Vector{b}\in\Gamma(\mathcal{N})$
and $k\Vector{a}\in\Gamma(\mathcal{N})$ for any
$\Vector{a}$, $\Vector{b} \in\Gamma(\mathcal{N})$ and for any $k\in\Gamma(\mathcal D)$.
		\qed
		\end{definition}

		\begin{definition}
Suppose $\Vector u$, $\Vector v\in\Gamma(\mathcal{V})$ are $\mathcal D\star$\hyph vector fields. We call
vector field $\Vector w$ \AddIndex{$\mathcal D\star$\hyph linear composition of vector fields}
{linear composition of vector fields} $\Vector u$ and $\Vector v$
when we can write $\Vector w=a\Vector u+b\Vector v$ where $a$ and $b$ are scalars.
		\qed
		\end{definition}

		\begin{remark}
		\label{remark: left and right vector bundle notation}
We extend to vector bundle and its type convention
described in remark
\xRef{0701.238}{remark: left and right vector space notation}.
We assume that the fiber of \Drc vector bundle is
\drc vector space.
		\qed
		\end{remark}

			\section{\texorpdfstring{$\mathcal D\RCstar$}{Drc}-Basis of Vector Bundle}
			\label{section: Basis of Vector Bundle}

		\begin{definition}
Vector fields ${}_{\gi a}\Vector{a}$, $\gi{a}\in\gi{I}$ of \Drc vector bundle
$\mathcal{V}$ are \AddIndex{\Drc linearly independent}
{linearly independent vector fields}
if $c=0$ follows from the equation
\[c\RCstar\Vector{a}=0\]
Otherwise vector fields ${}_{\gi a}\Vector{a}$
are \AddIndex{\Drc linearly dependent}{linearly dependent vector fields}.
		\qed
		\end{definition}

		\begin{definition}
We call set of vector fields
\symb{\Basis{e}=({}_{\gi a}\Vector{e},\gi{a}\in \gi{I})}1
{drc basis, vector bundle}
a \AddIndex{\Drc basis for vector bundle}{drc basis, vector bundle}
if vectors ${}_{\gi a}\Vector{e}$ are
\Drc linearly independent and adding to this system any other vector
we get a new system which is \Drc linearly dependent.
		\qed
		\end{definition}

		\begin{theorem}
		\label{theorem: expansion relative basis, vector bundle}
If $\Basis{e}$ is a \Drc basis of vector space $\mathcal{V}$
then any vector field $\Vector{v}\in \mathcal{V}$
has one and only one expansion 
	\begin{equation}
\Vector{v}=v\RCstar\Vector{e}
	\label{eq: expansion relative basis, vector bundle}
	\end{equation}
relative to this \Drc basis.
		\end{theorem}
		\begin{proof}
We check the statement of theorem in fiber and this statement is corollary
of theorem \xRef{0701.238}{theorem: expansion relative basis, vector space}.
		\end{proof}

		\begin{definition}
We call the matrix
$v$ in expansion \eqref{eq: expansion relative basis, vector bundle}
\AddIndex{coordinate matrix of vector field $\Vector{v}$ in \Drc basis $\Basis{e}$}
{coordinate matrix of vector field in drc basis}
and we call its elements
\AddIndex{coordinates of vector field $\Vector{v}$ in \Drc basis $\Basis{e}$}
{coordinates of vector field in drc basis}.
		\qed
		\end{definition}

According to construction we execute all operation over vector field in
fiber. Therefore we can use methods of
theory of \drc vector space in most cases.
However there is difference.
Coordinate matrix as well coordinates of vector field
is function of coordinates of point on the base.\footnote{However
I keep open the question about continuity of this function.}
This leads to different circumstances. For instance,
\drc linear dependence of vector fields in fiber
is conditioned by the singularity of
coordinate matrix of corresponding vectors.
This is one of the reason why we have problem to define
dimention when we consider \Drc vector bundle over
ring of sections $\Gamma(\mathcal D)$.
Considering \Drc vector bundle as fibered representation,
we get more flexibility in definition of \Drc basis.

		\begin{theorem}
		\label{theorem: coordinate vector bundle}
Set of coordinates $a$ of vector field $\Vector{a}$ relative \Drc basis $\Basis{e}$
forms \Drc vector bundle \symb{\mathcal D^{\gi n}}1{coordinate drc vector bundle}
isomorphic \Drc vector bundle $\mathcal{V}$.
This \Drc vector bundle is called
\AddIndex{coordinate \Drc vector bundle}
{coordinate drc vector bundle}.
This isomorphism
is called \AddIndex{fibered coordinate \Drc isomorphism}{fibered coordinate drc isomorphism}.
		\end{theorem}
		\begin{proof}
Suppose vectors $\Vector{a}$ and $\Vector{b}\in\mathcal{V}$ have expansion
\[
\Vector{a}=a\RCstar \Vector{e}
\]
\[
\Vector{b}=b\RCstar \Vector{e}
\]
relative basis $\Basis{e}$. Then
\[
\Vector{a}+\Vector{b}=
a\RCstar \Vector{e}+b\RCstar \Vector{e}
=(a+b)\RCstar\Vector{e}
\]
\[
m\Vector{a}=m(a\RCstar \Vector{e})=(ma)\RCstar\Vector{e}
\]
for any $m\in D$.
Thus, operations in a vector bundle are defined by coordinates
\[
(a+b)^{\gi a}=a^{\gi a}+b^{\gi a}
\]
\[
(ma)^{\gi a}=ma^{\gi a}
\]
This completes the proof.
		\end{proof}

Since we defined linear composition in fiber
the duality principle stated in theorems
\xRef{0701.238}{theorem: duality principle, vector space},
\xRef{0701.238}{theorem: duality principle, vector space, quasideterminants}
holds for vector bundles.

			\section{\texorpdfstring{$\mathcal D\RCstar$}{Drc}\hyph Linear Map of Vector Bundles}
			\label{section: Linear Map of Vector Bundles}

		\begin{definition}
		\label{definition: src trc linear map of vector bundles}
Suppose $\mathcal{A}$ is $\mathcal S\RCstar$\Hyph vector bundle.
Suppose $\mathcal{B}$ is $\mathcal T\RCstar$\Hyph vector bundle.
Morphism
	\[
\xymatrix{
\mathcal F:\mathcal S\ar[r]&\mathcal T&\Vector{\mathcal H}:\mathcal{A}\ar[r]&\mathcal{B}
}
	\]
of \Ts representations of fibered skew field 
in fibered Abelian group is called
\AddIndex{$(\mathcal S\RCstar,\mathcal T\RCstar)$\Hyph linear map of vector bundles}
{src trc linear map of vector bundles}.
		\qed
		\end{definition}

By theorem
\ref{theorem: decompositions of morphism of fibered representations}
studying $(\mathcal S\RCstar,\mathcal T\RCstar)$\Hyph linear map
we can consider case $\mathcal S=\mathcal T$.

		\begin{definition}
		\label{definition: drc linear map of vector bundles}
Suppose $\mathcal{A}$ and
$\mathcal{B}$ are \Drc vector bundles.
We call map
	\begin{equation}
\Vector{\mathcal H}:\mathcal{A}\rightarrow \mathcal{B}
	\label{eq: drc linear map of vector bundles}
	\end{equation}
\AddIndex{\Drc linear map of vector bundles}{drc linear map of vector bundles}
if
	\begin{equation}
	\label{eq: drc linear map of vector bundles, product on scalar}
\Vector{\mathcal H}(a\RCstar \Vector{m})=a\RCstar \Vector{\mathcal H}(\Vector{m})
	\end{equation}
for any $a^{\gi a} \in \Gamma(D)$, ${}_{\gi a}\Vector{m} \in \Gamma(\mathcal{A})$.
		\qed
		\end{definition}

		\begin{theorem}
		\label{theorem: drc linear map of drc vector bundle}
Let $\Vector{f}=({}_{\gi a}\Vector{f},\gi{a}\in \gi{I})$
be a \Drc basis of vector bundle $\mathcal{A}$
and $\Vector{e}=({}_{\gi b}\Vector{e},\gi{b}\in \gi{J})$ 
be a \Drc basis of vector bundle $\mathcal{B}$.
Then \Drc linear map \eqref{eq: drc linear map of vector bundles}
of vector bundles has presentation
	\begin{equation}
b=a\RCstar \mathcal H
	\label{eq: drc linear map of vector bundles, presentation}
	\end{equation}
relative to selected bases. Here
\begin{itemize}
\item $a$ is coordinate matrix of vector field
$\Vector a$ relative the \Drc basis $\Basis{f}$
\item $b$ is coordinate matrix of vector field
$$\Vector b=\Vector{\mathcal H}(\Vector a)$$ relative the \Drc basis $\Basis{e}$
\item $\mathcal H$ is coordinate matrix of set of vectors fields
$(\Vector {\mathcal H}({}_{\gi a}\Vector{f}))$ in \Drc basis $\Basis{e}$
called \AddIndex{matrix of fibered \Drc linear map}{matrix of fibered drc linear map}
relative bases $\Basis{f}$ and $\Basis{e}$
\end{itemize}
		\end{theorem}
		\begin{proof}
Vector field $\Vector{a}\in \Gamma(\mathcal{A})$ has expansion
	\[
\Vector{a}=a\RCstar \Vector{f}
	\]
relative to \Drc basis $\Basis{f}$.
Vector field $\Vector{b}=f(\Vector{a})\in \Gamma(\mathcal{B})$ has expansion
	\begin{equation}
\Vector{b}=b\RCstar \Vector{e}
	\label{eq: drc linear map of vector bundles, presentation, 1}
	\end{equation}
relative to \Drc basis $\Basis{e}$.

Since $\Vector {\mathcal H}$ is a \Drc linear map, from
\eqref{eq: drc linear map of vector bundles, product on scalar} it follows
that
	\begin{equation}
\Vector{b}=\Vector {\mathcal H}(\Vector{a})=\Vector {\mathcal H}(a\RCstar \Vector{f})
=a\RCstar \Vector {\mathcal H}(\Vector{f})
	\label{eq: drc linear map of vector bundles, presentation, 2}
	\end{equation}
$\Vector {\mathcal H}({}_{\gi a}\Vector{f})$ is also a vector field of vector bundle $\mathcal{B}$ and has
expansion
	\begin{equation}
\Vector {\mathcal H}({}_{\gi a}\Vector{f})
={}_{\gi a}\mathcal H\RCstar \Vector{e}
={}_{\gi a}\mathcal H^{\gi b}\ {}_{\gi b}\Vector{e}
	\label{eq: drc linear map of vector bundles, presentation, 3}
	\end{equation}
relative to basis $\Basis{e}$.
Combining \eqref{eq: drc linear map of vector bundles, presentation, 2}
and \eqref{eq: drc linear map of vector bundles, presentation, 3} we get
	\begin{equation}
\Vector{b}=a\RCstar \mathcal H\RCstar\Vector{e}
	\label{eq: drc linear map of vector bundles, presentation, 4}
	\end{equation}
\eqref{eq: drc linear map of vector bundles, presentation}  follows
from comparison of \eqref{eq: drc linear map of vector bundles, presentation, 1}
and \eqref{eq: drc linear map of vector bundles, presentation, 4} and
theorem \ref{theorem: expansion relative basis, vector bundle}.
		\end{proof}

On the basis of theorem \ref{theorem: drc linear map of drc vector bundle}
we identify the fibered \Drc linear map \eqref{eq: drc linear map of vector bundles} of vector spaces
and the matrix of its presentation \eqref{eq: drc linear map of vector bundles, presentation}.

		\begin{theorem}
		\label{theorem: product of drc linear maps, vector bundles}
Let $$\Basis{f}=({}_{\gi a}\Vector{f},\gi{a}\in \gi{I})$$ be a \Drc basis of vector bundles $\mathcal{A}$,
$$\Basis{e}=({}_{\gi b}\Vector{e},\gi{b}\in \gi{J})$$ be a \Drc basis of vector bundles $\mathcal{B}$,
and $$\Basis{g}=({}_{\gi c}\Vector{g},\gi{c}\in \gi{L})$$ be a \Drc basis of vector bundles $\mathcal{C}$.
Suppose diagram of \Drc linear maps
	\[
\xymatrix{
\mathcal{A}\ar[rr]^{\mathcal H}\ar[dr]^{\mathcal F} & & \mathcal{C}\\
&\mathcal{B}\ar[ur]^{\mathcal G} &
}
	\]
is commutative diagram
where \drc linear map $\mathcal F$  has presentation
	\begin{equation}
	\label{eq: product of drc linear maps, tstar vector bundles, A}
b=a\RCstar \mathcal F
	\end{equation}
relative to selected bases
and \drc linear map $\mathcal G$ has presentation
	\begin{equation}
	\label{eq: product of drc linear maps, tstar vector bundles, B}
c=b\RCstar \mathcal G
	\end{equation}
relative to selected bases.
Then \drc linear map $\mathcal H$ has presentation
	\begin{equation}
	\label{eq: product of drc linear maps, tstar vector bundles, C}
c=a\RCstar \mathcal F\RCstar \mathcal G
	\end{equation}
relative to selected bases.
		\end{theorem}
		\begin{proof}
Proof of the statement follows from substituting
\eqref{eq: product of drc linear maps, tstar vector bundles, A}
into \eqref{eq: product of drc linear maps, tstar vector bundles, B}.
		\end{proof}

Presenting \drc linear map as \RC product we can rewrite
\eqref{eq: drc linear map of vector bundles, product on scalar} as
	\begin{equation}
(ka)\RCstar \mathcal F=k(a\RCstar \mathcal F)
	\label{eq: drc linear map of vector bundles, product on scalar, image}
	\end{equation}
We can expres the statement of the theorem
\ref{theorem: product of drc linear maps, vector bundles}
in the next form
	\begin{equation}
(a\RCstar \mathcal F)\RCstar \mathcal G=a\RCstar (\mathcal F\RCstar \mathcal G)
	\label{eq: product of drc linear maps, vector bundles, image}
	\end{equation}
Equations
\eqref{eq: drc linear map of vector bundles, product on scalar, image} and
\eqref{eq: product of drc linear maps, vector bundles, image}
represent the
\AddIndex{associative law for \Drc linear maps of vector bundles}
{associative law for drc linear maps of vector bundles}.
This allows us writing of such expressions without using of brackets.

Equation \eqref{eq: drc linear map of vector bundles, presentation}
is coordinate notation for fibered \Drc linear map.
Based theorem \ref{theorem: drc linear map of drc vector bundle}
non coordinate notation also can be expressed using \RC product
	\begin{equation}
\Vector{b}=\Vector{a}\RCstar\Vector {\mathcal F}=a\RCstar \Vector{f}\RCstar\Vector {\mathcal F}
=a\RCstar \mathcal F\RCstar \Vector{e}
	\label{eq: drc linear map of vector bundles, presentation, 5}
	\end{equation}
If we substitute equation \eqref{eq: drc linear map of vector bundles, presentation, 5}
into theorem
\ref{theorem: product of drc linear maps, vector spaces},
then we get chain of equations
	\begin{align*}
\Vector{c}&=\Vector{b}\RCstar\Vector {\mathcal G}=b\RCstar \Vector{e}\RCstar\Vector {\mathcal G}
=b\RCstar \mathcal G\RCstar \Vector{g}\\
\Vector{c}&=\Vector{a}\RCstar\Vector {\mathcal F}\RCstar\Vector {\mathcal G}
=a\RCstar \Vector{f}\RCstar\Vector {\mathcal F}\RCstar\Vector {\mathcal G}
=a\RCstar \mathcal F\RCstar \mathcal G\RCstar \Vector{g}
	\end{align*}

%auto-ignore
\OpenBiblio

%22. A.G. Kurosh, Lectures on General Algebra, Chelsea, New York, 1963

\BiblioItem{texSpaceTime}{Einstein: Geometry and Experience}
{
Albert Einstein, Geometry and Experience, (1921)\\
Albert Einstein, Sidelights on Relativity, 25 - 56,\\
Courier Dover Publications, 1983
}%

\BiblioItem{texSpaceTime}{Einstein: Main problems of general relativity}
{
Albert Einstein,
Grundgedanken und Probleme der Relativit\"atstheorie, (1923),\\
Nobelstiftelsen, Les Prix Nobel en 1921 - 1922,
Imprimerie Royale, Stockholm, 1923
}%

\BiblioItem{texSpaceTime}{Einstein: Isaak Newton}
{
Albert Einstein,
Isaak Newton,
Manchester Guardian, 16, 234 - 235 (1927)
}%

\BiblioItem{texSpaceTime}{Einstein: Foundations of general relativity}
{
Albert Einstein,
Die Grundlage der allgemeinen Relativit\"atstheorie,
Ann. Phys., 1916, {\bf 49}, 769 - 822,\\
Einstein's Annalen Papers: The Complete Collection 1901-1922,
edited by J\"urgen Renn, 517 - 571,\\
Wiley-VCH Verlag GmbH & Co. KGaA, 2005
}%

\BiblioItem{texGenRelativity}{Ghez}
{
Ghez et al.,
The First Measurement of Spectral Lines in a Short-Period Star Bound to the Galaxy's Central Black Hole: A Paradox of Youth,
\href{http://www.journals.uchicago.edu/ApJ/journal/issues/ApJL/v586n2/16990/brief/16990.abstract.html}{ApJL, 586, L127} (2003),
eprint \href{http://arxiv.org/abs/astro-ph/0302299}{arXiv:astro-ph/0302299} (2003)
}%

\BiblioItem{texGenRelativity}{Schodel}
{
R. Sch\"odel et al.,
A star in a 15.2-year orbit around the supermassive black hole at the centre of the Milky Way,
\href{http://www.nature.com/cgi-taf/DynaPage.taf?file=/nature/journal/v419/n6908/abs/nature01121_fs.html}{Nature 419, 694} (2002)
}%

\BiblioItem{texAffine,texGeomObject}{Mielke}
{
Eckehard W. Mielke, Affine generalization of the Komar complex of general relativity,
\href{http://prola.aps.org/searchabstract/PRD/v63/i4/e044018}{Phys. Rev. D 63, 044018} (2001)
}%

\BiblioItem{texAffine}{Obukhov}
{
Yu. N. Obukhov and J. G. Pereira, Metric\hyph affine approach to teleparallel gravity,
\href{http://scitation.aip.org/getabs/servlet/GetabsServlet?prog=normal&id=PRVDAQ000067000004044016000001&idtype=cvips&gifs=Yes}
{Phys. Rev. D 67, 044016} (2003),
eprint \href{http://arxiv.org/abs/gr-qc/0212080}{arXiv:gr-qc/0212080} (2002)
}%

\BiblioItem{texAffine}{Sardanashvily}
{
Giovanni Giachetta, Gennadi Sardanashvily, Dirac Equation in Gauge and Affine-Metric Gravitation Theories,
eprint \href{http://arxiv.org/abs/gr-qc/9511035}{arXiv:gr-qc/9511035} (1995)
}%

\BiblioItem{texAffine}{Gauge}
{
Frank Gronwald and Friedrich W. Hehl, On the Gauge Aspects of Gravity, eprint
\href{http://arxiv.org/abs/gr-qc/9602013}{arXiv:gr-qc/9602013} (1996)
}%

\BiblioItem{texAffine}{Neeman}
{
Yuval Neeman, Friedrich W. Hehl, Test Matter in a Spacetime with Nonmetricity, eprint
\href{http://arxiv.org/abs/gr-qc/9604047}{arXiv:gr-qc/9604047} (1996)
}%

\BiblioItem{texTidal,texAffine,texGeomObject}{torsion}
{
F. W. Hehl, P. von der Heyde, G. D. Kerlick, and J. M. Nester,
General relativity with spin and torsion: Foundations and prospects,\\
\href{http://prola.aps.org/abstract/RMP/v48/i3/p393_1}{Rev. Mod. Phys. 48, 393} (1976)
}%

\BiblioItem{texTidal,texNewton}{Megged}
{
O. Megged, Post-Riemannian Merger of Yang-Mills Interactions with Gravity,
eprint \href{http://arxiv.org/abs/hep-th/0008135}{arXiv:hep-th/0008135} (2001)
}%

%\BiblioItem{texNewton}{Hehl}
%{
%Friedrich W. Hehl, Uwe Muench,
%eprint \href{http://arxiv.org/abs/gr-qc/9708007}{arXiv:gr-qc/9708007} (1997)
%}%

\BiblioItem{texNewton}{gr-qc-9604027}
{
Yu.N. Obukhov, E.J. Vlachynsky, W. Esser, R. Tresguerres and F.W. Hehl,
An exact solution of the metric\hyph affine gauge theory with dilation, shear, and spin charges,
eprint \href{http://arxiv.org/abs/gr-qc/9604027}{arXiv:gr-qc/9604027} (1996)
}%

\BiblioItem{texLagrange}{Weinberg}
{
Steven Weinberg. The Quantum Theory of Fields. Cambridge university press.
}%

\BiblioItem{texLagrange}{Reinhardt}
{
Greiner Reinhardt. Field Quantization. Springer.
}%

\BiblioItem{texLagrange}{Landau}
{
L. D. Landau, E. M. Lifshich, The classical theory of fields.
Oxford, New York, Pergamon Press
}%

\BiblioItem{texTidal}{Wheeler}
{
Ignazio Ciufolini, John Wheeler. Gravitation and Inertia.
Princeton university press.
}%

\BiblioItem{texPrefaceTidal,texPrefaceRefernceFrame}{Anderson02}
{
J. D. Anderson, P. A. Laing, E. L. Lau, A. S. Liu, M. M. Nieto, and S. G. Turyshev,
Study of the anomalous acceleration of Pioneer 10 and 11,
\href{http://prola.aps.org/searchabstract/PRD/v65/i8/e082004}{Phys. Rev. D 65, 082004, 50 pp.}, (2002),
eprint \href{http://arxiv.org/abs/gr-qc/0104064}{arXiv:gr-qc/0104064} (2001)
}%

\BiblioItem{texTidal}{Anderson98}
{
J. D. Anderson, P. A. Laing, E. L. Lau, A. S. Liu, M. M. Nieto, and S. G. Turyshev,
Indication, from Pioneer 10/11, Galileo, and Ulysses Data, of an Apparent Anomalous, Weak, Long-Range Acceleration,
\href{http://prola.aps.org/abstract/PRL/v81/i14/p2858_1}{Phys. Rev. Lett. 81, 2858}, (1998),
eprint \href{http://arxiv.org/abs/gr-qc/9808081}{arXiv:gr-qc/9808081} (1998)
}%

%\BiblioItem{Havas} Peter Havas, The Classical Equations of Motion of Point Particles, I,
%{
%\href{http://prola.aps.org/abstract/PR/v87/i2/p309_1}{Phys. Rev. 87, 309} (1952)
%}%

\BiblioItem{texReferenceFrame,texFiberedAlgebra}{Serge Lang}
{
Serge Lang,
Algebra, Springer, 2002
}%

\BiblioItem{texFiberedAlgebra,texTstarMorphism}{Burris Sankappanavar}
{
S. Burris, H.P. Sankappanavar,
A Course in Universal Algebra, Springer-Verlag (March, 1982),
\\eprint
\href{http://www.math.uwaterloo.ca/~snburris/htdocs/ualg.html}
{http://www.math.uwaterloo.ca/~snburris/htdocs/ualg.html}
\\(The Millennium Edition)
}%

\BiblioItem{texGeomObject}{Shilov}
{
G. E. Shilov,
Calculus, Multivariable Functions,
Moscow, Nauka, 1972
}%

\BiblioItem{texBundle}{Kolmogorov Fomin}
{
A. N. Kolmogorov and S. V. Fomin,
Elements of the Theory of Functions and Functional Analysis,
Courier Dover Publication, 1999
}%

\BiblioItem{texAffine,texRepresentation,texBasis,texDrcBasis,texLinearMap,texVectorSpace,texPolyvector}{Rashevsky}
{
P. K. Rashevsky, Riemann Geometry and Tensor Calculus,\\
Moscow, Nauka, 1967
}%

\BiblioItem{texPolyvector}{Dubrovin Fomenko Novikov part 1}
{
B. A. Dubrovin, A. T. Fomenko, S. P. Novikov,
Modern Geometry - Methods and Applications,\\
Part 1, The Geometry of Surfaces, Transformation Groups, and Fields,\\
Translated by Robert G. Burns,\\
Springer - New York, 1992
}%

\BiblioItem{texDrcBasis,texBasis}{Korn}
{
Granino A. Korn, Theresa M. Korn,
Mathematical Handbook for Scientists and Engineer,
McGraw-Hill Book Company, New York, San Francisco,
Toronto, London, Sydney, 1968
}%

\BiblioItem{texBundle}{Hocking Young Topology}
{
John G. Hocking, Gail S. Young,
Topology,\\
Courier Dover Publications, 1988
}%

%\BiblioItem{Kleyn} http://www.geocities.com/aleks\_kleyn/Derivative/Derivative.htm

\BiblioItem{texPrefaceRefernceFrame}{Tartaglia}
{
Angelo Tartaglia and Matteo Luca Ruggiero,
Angular Momentum Effects in Michelson\Hyph Morley Type Experiments,
Gen.Rel.Grav. 34, 1371-1382 (2002),\\
eprint \href{http://arxiv.org/abs/gr-qc/0110015}{arXiv:gr-qc/0110015} (2001)
}%

\BiblioItem{texPrefaceRefernceFrame}{Tomozawa}
{
Yukio Tomozawa, Speed of Light in Gravitational Fields, eprint
\href{http://arxiv.org/abs/astro-ph/0303047}{arXiv:astro-ph/0303047} (2004)
}%

\BiblioItem{texPrefaceRefernceFrame}{Magueijo}
{
Joao Magueijo,
Covariant and locally Lorentz-invariant varying speed of light theories,
\href{http://prola.aps.org/abstract/PRD/v62/i10/e103521}{Phys. Rev. D 62, 103521} (2000),
eprint \href{http://arxiv.org/abs/gr-qc/0007036}{arXiv:gr-qc/0007036} (2000)
}%

\BiblioItem{texPrefaceRefernceFrame}{Bassett}
{
Bruce A. Bassett, Stefano Liberati, Carmen Molina-Paris, and Matt Visser,
Geometrodynamics of variable-speed-of-light cosmologies,
\href{http://prola.aps.org/abstract/PRD/v62/i10/e103518}{Phys. Rev. D 62}, 103518 (2000),
eprint \href{http://arxiv.org/abs/astro-ph/0001441}{arXiv:astro-ph/0001441} (2000)
}%

\BiblioItem{texPrefaceRefernceFrame}{Straumann}
{
Lochlainn O'Raifeartaigh and Norbert Straumann,
Gauge theory: Historical origins and some modern developments,
\href{http://prola.aps.org/abstract/RMP/v72/i1/p1_1}{Rev. Mod. Phys. 72, 1} (2000)
}%

\BiblioItem{texPrefaceRefernceFrame}{Lammerzahl}
{
Claus L\"ammerzahl, Mark P. Haugan,
On the interpretation of Michelson\Hyph Morley experiments,
%\href{http://www.sciencedirect.com/science?\_ob=ArticleURL&\_udi=B6TVM-42WP7CR-1&\_user=10&\_handle=W-WA-A-A-AZ-MsSAYZW-UUW-AUDDYZYZAU-WZCBYCEDW-AZ-U&\_fmt=summary&\_coverDate=04%2F23%2F2001&\_rdoc=1&\_orig=browse&\_srch=%23toc%235538%232001%23997179995%23246657!&\_cdi=5538&view=c&\_acct=C000050221&\_version=1&\_urlVersion=0&\_userid=10&md5=385478cda8c5568dea1aeaf0c43669da}
{Phys. Lett. A282 223-229} (2001),\\
eprint \href{http://arxiv.org/abs/gr-qc/0103052}{arXiv:gr-qc/0103052} (2001)
}%

\BiblioItem{texPrefaceRefernceFrame}{Muller}
{
Holger Muller et al.,
Modern Michelson-Morley Experiment using Cryogenic Optical Resonators,
\href{http://prola.aps.org/searchabstract/PRL/v91/i2/e020401}{Phys. Rev. Lett. 91, 020401} (2003),
eprint \href{http://arxiv.org/abs/physics/0305117}{arXiv:physics/0305117} (2000)
}%

\BiblioItem{texPrefaceRefernceFrame,texPrefaceTidal}{Ranada}
{
Antonio F. Ranada,
Pioneer acceleration and variation of light speed: experimental situation,
eprint \href{http://arxiv.org/abs/gr-qc/0402120}{arXiv:gr-qc/0402120} (2004)
}%

\BiblioItem{texBiring,texVectorSpace}{math.QA-0208146}
{
I. Gelfand, S. Gelfand, V. Retakh, R. Wilson,
Quasideterminants,\\
eprint \href{http://arxiv.org/abs/math.QA/0208146}{arXiv:math.QA/0208146} (2002)
}%

\BiblioItem{texBiring,texVectorSpace}{q-alg-9705026}
{
I.Gelfand, V.Retakh,
Quasideterminants, I,\\
eprint \href{http://arxiv.org/abs/q-alg/9705026}{arXiv:q-alg/9705026} (1997)
}%

\BiblioItem{texVectorSpace}{Gelfand Retakh 1991}
{
I. Gelfand and V. Retakh, Determinants of Matrices over Noncommutative Rings, Funct.
Anal. Appl. 25 (1991), no. 2, 91-102
}%

\BiblioItem{texVectorSpace}{Gelfand Retakh 1992}
{
I. Gelfand and V. Retakh, A Theory of Noncommutative Determinants and Characteristic
Functions of Graphs, Funct. Anal. Appl. 26 (1992), no. 4, 1-20
}%

\BiblioItem{texVectorSpace}{hep-th-9407124}
{
I. M. Gelfand, D. Krob, A. Lascoux, B. Leclerc, V.S. Retakh and J.-Y. Thibon,
Noncommutative symmetric functions,\\
eprint \href{http://arxiv.org/abs/hep-th/9407124}{arXiv:hep-th/9407124} (1994)
}%

\BiblioItem{texVectorSpace}{Carl Faith 1}
{
Carl Faith, Algebra: Rings, Modules and Categories I,
Springer - Verlag, Berlin - Heidelberg - New York, 1973
}%

%\BiblioItem{texVectorSpace}{Pareigis}
%{
%Bodo Pareigis, Categories and Functors,
%Academic Press - New York - London, 1970
%}%

%\BiblioItem{texVectorSpace}{Beachy}
%{
%John A. Beachy, Introductory Lectures on Rings i Modules,
%Cambridge University Press, 1999
%}%

\BiblioItem{texDrcReferenceFrame,texRefernceFrame,texLie,texLieRepresentation}{0412.391}
{
Aleks Kleyn,
Basis Manifold,
eprint \href{http://arxiv.org/abs/math.DG/0412391}{arXiv:math.DG/0412391} (2004)
}%

\BiblioItem{texAffine}{0405.027}
{
Aleks Kleyn,
Reference Frame in General Relativity,\\
eprint \href{http://arxiv.org/abs/gr-qc/0405027}{arXiv:gr-qc/0405027} (2004)
}%

\BiblioItem{texTidal}{0405.028}
{
Aleks Kleyn, Metric\hyph Affine Manifold,
eprint \href{http://arxiv.org/abs/gr-qc/0405028}{arXiv:gr-qc/0405028} (2004)
}%

\BiblioItem{texFiberedAlgebra,texBundleRelation,texTstarMorphism}{0701.238}
{
Aleks Kleyn,
Lectures on Linear Algebra over Skew Field,\\
eprint \href{http://arxiv.org/abs/math.GM/0701238}{arXiv:math.GM/0701238} (2007)
}%

\BiblioItem{texBundleRelation,texPrefaceRelation}{0702.561}
{
Aleks Kleyn,
Algebra Bundle,\\
eprint \href{http://arxiv.org/abs/math.DG/0702561}{arXiv:math.DG/0702561} (2007)
}%

%\BiblioItem{texBasis,texLinearMap}{math.RA-0501237}
%{
%Aleks Kleyn,
%Vector Space Over Skew-Field,\\
%eprint \href{http://arxiv.org/abs/math.RA/0412391}{arXiv:math.RA/0501237} (2005)
%}%

\BiblioItem{texPolymodule}{math.RA-0501237v1}
{
Aleks Kleyn,
Module Over Skew-Field, version 1,\\
eprint \href{http://arxiv.org/abs/math/0501237v1}{arXiv:math.RA/0501237v1} (2005)
}%

\ifx\texBiring\Defined
\else
\BiblioItem{texVectorSpace,texFiberedAlgebra}{0612.111}
{
Aleks Kleyn,
Biring of Matrices,\\
eprint \href{http://arxiv.org/abs/math.OA/0612111}{arXiv:math.OA/0612111} (2006)
}%
\fi

\ifx\texBundleRelation\Defined
\else
\BiblioItem{texFiberedMorphism}{0707.2246}
{
Aleks Kleyn,
Fibered Correspondence,\\
eprint \href{http://arxiv.org/abs/0707.2246}{arXiv:0707.2246} (2007)
}%
\fi

\ifx\PrintBook\Defined
\BiblioItem{texPrefaceRelation}{0707.2246}
{
Aleks Kleyn,
Fibered Correspondence,\\
eprint \href{http://arxiv.org/abs/0707.2246}{arXiv:0707.2246} (2007)
}%
\fi

\BiblioItem{texHomotopy}{q-alg-9705009}
{
John C. Baez,
An Introduction to n-Categories,\\
eprint \href{http://arxiv.org/abs/q-alg/9705009}{arXiv:q-alg/9705009} (1997)
}%

\BiblioItem{texPrefaceRelation}{Tolstoi about Anna Karenina}
{
Tolstoi about Anna Karenina,
in book A Karenina Companion, by C. J. G. Turner,
published by Wilfrid Laurier University Press (August 1993)
}%

\BiblioItem{texBundleRelation,texPrefaceRelation,texTstarMorphism,texBundle}
{Cohn: Universal Algebra}
{
Paul M. Cohn,
Universal Algebra,
Springer, 1981
}%

\BiblioItem{texBundle}
{Maunder: Algebraic Topology}
{
C. R. F. Maunder,
Algebraic Topology,
Dover Publications, Inc, Mineola, New York, 1996
}%

\BiblioItem{texFiberedAlgebra}{Pommaret: Partial Differential Equations}
{
J.-F. Pommaret,
Partial Differential Equations and Group Theory,
Springer, 1994
}%

\BiblioItem{texBundleRelation}{Bourbaki: Set Theory}
{
N. Bourbaki,
Theory of sets,
Springer, 2004
}%

\BiblioItem{texBundle,texCartesian,texFiberedAlgebra,texBundleRelation,texFiberedMorphism}
{Bourbaki: General Topology 1}
{
N. Bourbaki,
General Topology, Chapters 1 - 4,
Springer, 1989
}

\BiblioItem{texCalculus}{Bourbaki: General Topology: Chapter 5 - 10}
{
N. Bourbaki,
General Topology, Chapters 5 - 10,
Springer, 1989
}

\BiblioItem{texCalculus}{Bourbaki: Topological Vector Space}
{
N. Bourbaki,
Topological Vector Spaces, Chapters 1 - 5,
Springer, 1987
}

\BiblioItem{texCalculus}{Pontryagin: Topological Group}
{
L. S. Pontryagin,
Selected Works, Volume Two, Topological Groups,
Gordon and Breach Science Publishers, 1986
}

\BiblioItem{texFiberedMorphism}{Postnikov: Differential Geometry}
{
Postnikov M. M.,
Geometry IV: Differential geometry,
Moscow, Nauka, 1983
}

\BiblioItem{texFiberedAlgebra,texFiberedMorphism}{Hatcher: Algebraic Topology}
{
Allen Hatcher,
Algebraic Topology,
Cambridge University Press, 2002
}

\BiblioItem{texFiberedMorphism}{geometry of differential equations}
{
Vinogradov, A. M., Krasil'shchik, I. S., and Lychagin, V. V.,
Introduction to geometry of nonlinear differential equations,
Nauka, Moscow, 1986
}

\BiblioItem{texFiberedMorphism}{cohomological analysis}
{
A. M. Vinogradov,
Cohomological Analysis of Partial Differential Equations
and Secondary Calculus,
American Mathematical Society, 2001
}

\BiblioItem{texPolyvector}{0801.1734}
{
Brandon S. DiNunno, Richard A. Matzner,
The Volume Inside a Black Hole,\\
eprint \href{http://arxiv.org/abs/0801.1734v1}{arXiv:0801.1734v1} (2008)
}

\CloseBiblio

%auto-ignore
\OpenIndex
\SetIndexSpace%
\Index{texLinearMap}%1%1
   {$1$-\drc form}%
   {1-drc form, vector spaces}%
\SetIndexSpace%
\Index{texPolymodule}%2%66
   {$(2)$\hyph vector space}%
   {(2)-vector space}%
\Index{texBundleRelation}%2%122
   {$2$\Hyph ary fibered relation}%
   {2 ary fibered relation}%
\SetIndexSpace%
\Index{texCalculus}%A%67
   {$A$\Hyph valued function}%
   {A valued function}%
\Index{texBiring}%A%2
   {$(^{\gi a}_{\gi b})$\hyph \CR quasideterminant}%
   {a b cr-quasideterminant}%
\Index{texBiring}%A%34
   {$(^{\gi a}_{\gi b})$\hyph \RC quasideterminant}%
   {a b RC-quasideterminant}%
\Index{texCalculus}%A%100
   {absolute value on skew field}%
   {absolute value on skew field}%
\Index{texBasis}%A%262
   {active representation}%
   {active representation}%
\Index{texDrcBasis}%A%57
   {active \sT representation}%
   {active representation, vector space}%
\Index{texBasis}%A%263
   {active transformation on basis manifold}%
   {active transformation}%
\Index{texBasis}%A%275
   {affine basis}%
   {Affine Basis}%
\Index{texBasis}%A%274
   {affine transformation group}%
   {AffineTransformationGroup}%
\Index{texTypeBasis}%A%77
   {affine transformation group}%
   {AffineTransformationGroup}%
\Index{texBasis}%A%277
   {affine transformation on basis manifold}%
   {affine transformation}%
\Index{texBiring}%A%59
   {alternative representation of matrix}%
   {Alternative representation}%
\Index{texRefernceFrame}%A%125
   {anholonomic coordinate}%
   {anholonomic coordinate}%
\Index{texRefernceFrame}%A%110
   {anholonomic coordinates of connection}%
   {anholonomic coordinates of connection}%
\Index{texRefernceFrame}%A%126
   {anholonomic coordinates of vector}%
   {vector anholonomic coordinates}%
\Index{texRefernceFrame}%A%127
   {anholonomic coordinates on manifold}%
   {anholonomic coordinates on manifold}%
\Index{texRefernceFrame}%A%130
   {anholonomity object}%
   {anholonomity object}%
\Index{texFiberedGroup}%A%60
   {antihomomorphism of fibered groups}%
   {antihomomorphism of fibered groups}%
\Index{texBundleRelation}%A%156
   {antisymmetric $2$\Hyph ary fibered relation}%
   {antisymmetric 2 ary fibered relation}%
\Index{texFiberedAlgebra}%A%61
   {arity of operation}%
   {arity of operation}%
\Index{texTstarRepresentation}%A%80
   {associative law for covariant \Ts representation}%
   {associative law for Tstar covariant representation}%
\Index{texVectorField}%A%385
   {associative law for \Drc linear maps of vector bundles}%
   {associative law for drc linear maps of vector bundles}%
\Index{texDrcMorphism}%A%312
   {associative law for \drc linear maps of vector spaces}%
   {associative law for drc linear maps of vector spaces}%
\Index{texVectorField}%A%367
   {associative law for $\mathcal D\star$\Hyph vector fields}%
   {associative law, Dstar vector fields}%
\Index{texVectorSpace}%A%310
   {associative law for $D\star$\Hyph vector space}%
   {associative law, Dstar vector space}%
\Index{texLinearMap}%A%311
   {associative law for twin representations}%
   {associative law for twin representations}%
\Index{texBundleRelation}%A%155
   {associative law of composition of fibered correspondences}%
   {associative law, composition of fibered correspondences}%
\Index{texAffine}%A%56
   {auto parallel line}%
   {auto parallel line}%
\SetIndexSpace%
\Index{texBundleRelation}%B%233
   {base of fibered correspondence}%
   {base of fibered correspondence}%
\Index{texBundle}%B%172
   {base of map}%
   {base of map}%
\Index{texTypeBasis}%B%16
   {basis}%
   {}%
\SubIndex{texTypeBasis}%1
   {affine}%
   {Affine Basis}%
\SubIndex{texTypeBasis}%3
   {central affine}%
   {Central Affine Basis}%
\SubIndex{texTypeBasis}%2
   {orthonornal}%
   {Orthonornal Basis}%
\Index{texDrcBasis}%B%120
   {basis manifold}%
   {}%
\SubIndex{texTypeBasis}%25
   {of affine space}%
   {Basis Manifold, Affine Space}%
\SubIndex{texTypeBasis}%27
   {of central affine space}%
   {Basis Manifold, Central Affine Space}%
\SubIndex{texDrcBasis}%24
   {of \drc vector space}%
   {basis manifold of vector space}%
\SubIndex{texTypeBasis}%26
   {of Euclid space}%
   {Basis Manifold, Euclid Space}%
\Index{texBasis}%B%276
   {basis manifold of affine space}%
   {Basis Manifold, Affine Space}%
\Index{texBasis}%B%280
   {basis manifold of central affine space}%
   {Basis Manifold, Central Affine Space}%
\Index{texBasis}%B%282
   {basis manifold of Euclid space}%
   {Basis Manifold, Euclid Space}%
\Index{texBasis}%B%268
   {basis manifold of vector space}%
   {basis manifold of vector space}%
\Index{texBasis}%B%266
   {basis of vector space}%
   {Basis}%
\Index{texLieRepresentation}%B%62
   {basis vector}%
   {}%
\SubIndex{texLieRepresentation}%11
   {of \sT representation}%
   {basis vector of starT representation}%
\SubIndex{texLieRepresentation}%12
   {of \Ts representation}%
   {basis vector of Tstar representation}%
\Index{texBiring}%B%65
   {biring}%
   {biring}%
\Index{texFiberedMorphism}%B%349
   {bundle of level $2$}%
   {bundle of level 2}%
\Index{texFiberedMorphism}%B%350
   {bundle of level $n$}%
   {bundle of level n}%
\SetIndexSpace%
\Index{texVectorSpace}%C%14
   {\subs rows \dcr vector space}%
   {subs rows dcr vector space}%
\Index{texBiring}%C%3
   {\subs row of matrix}%
   {c row}%
\Index{texBiring}%C%10
   {$c$\hyph row of matrix}%
   {c-row}%
\Index{texAffine}%C%189
   {Cartan connection}%
   {Cartan connection}%
\Index{texAffine}%C%94
   {Cartan curvature}%
   {Cartan curvature}%
\Index{texAffine}%C%174
   {Cartan derivative}%
   {Cartan derivative}%
\Index{texAffine}%C%190
   {Cartan symbol}%
   {Cartan symbol}%
\Index{texAffine}%C%143
   {Cartan transport}%
   {Cartan transport}%
\Index{texBundle}%C%357
   {Cartesian power $\mathcal{A}$ of bundle $\mathcal{B}$}%
   {Cartesian power of bundle}%
\Index{texBundle}%C%331
   {Cartesian power $A$ of set $B$}%
   {Cartesian power of set}%
\Index{texCartesian}%C%176
   {Cartesian power of bundle}%
   {Cartesian power of bundle}%
\Index{texCartesian}%C%180
   {direct product of bundles}%
   {Cartesian product of bundles}%
\Index{texCartesian}%C%181
   {direct product of total spaces}%
   {Cartesian product of total spaces}%
\Index{texHomotopy}%C%4
   {category of \drc vector spaces}%
   {category of drc vector spaces}%
\Index{texBundleRelation}%C%112
   {category of fibered correspondences over diagonal}%
   {category of fibered correspondences over diagonal}%
\Index{texBundleRelation}%C%111
   {category of reduced fibered correspondences}%
   {category of reduced fibered correspondences}%
\Index{texBasis}%C%279
   {central affine basis}%
   {Central Affine Basis}%
\Index{texRepresentation}%C%261
   {column vector}%
   {column vector}%
\Index{texBundleRelation}%C%228
   {commutative diagram of correspondences}%
   {commutative diagram of correspondences}%
\Index{texBundle}%C%332
   {compact\hyph open topology}%
   {compact open topology}%
\Index{texDiffEq}%C%71
   {completely integrable system}%
   {completely integrable system}%
\Index{texBundleRelation}%C%188
   {composition of fibered correspondences}%
   {composition of fibered correspondences}%
\SubIndex{texVectorSpace}%74
   {of \drc linear equations}%
   {extended matrix, system of drc linear equations}%
\SubIndex{texVectorSpace}%75
   {of \rcd linear equations}%
   {extended matrix, system of rcd linear equations}%
\Index{texBundleRelation}%C%154
   {composition of reduced fibered correspondences}%
   {composition of reduced fibered correspondences}%
\Index{texBiring}%C%208
   {condition of reducibility of products}%
   {condition of reducibility of products}%
\Index{texBundleRelation}%C%339
   {continuous correspondence}%
   {continuous correspondence}%
\Index{texFiberedGroup}%C%317
   {contravariant \Ts representation of fibered group}%
   {Tstar contravariant representation of fibered group}%
\Index{texVectorSpace}%C%89
   {coordinate \drc isomorphism}%
   {coordinate drc isomorphism}%
\Index{texVectorField}%C%380
   {coordinate \Drc vector bundle}%
   {coordinate drc vector bundle}%
\Index{texVectorSpace}%C%87
   {coordinate \drc vector space}%
   {coordinate drc vector space}%
\Index{texBasis}%C%286
   {coordinate isomorphism}%
   {coordinate isomorphism}%
\Index{texVectorSpace}%C%21
   {coordinate matrix}%
   {}%
\SubIndex{texVectorSpace}%7
   {of set of vectors in \dcr rows vector space}%
   {coordinate matrix of set of vectors, dcr vector space}%
\SubIndex{texVectorSpace}%8
   {of set of vectors in \drc rows vector space}%
   {coordinate matrix of set of vectors, drc vector space}%
\SubIndex{texVectorSpace}%6
   {of vector in \drc basis}%
   {coordinate matrix of vector in drc basis}%
\Index{texVectorField}%C%377
   {coordinate matrix of vector field in \Drc basis}%
   {coordinate matrix of vector field in drc basis}%
\Index{texRefernceFrame}%C%86
   {coordinate reference frame}%
   {coordinate reference frame}%
\Index{texDrcBasis}%C%88
   {coordinate representation in \drc vector space}%
   {coordinate representation, vector space}%
\Index{texBasis}%C%287
   {coordinate representation of group in vector space}%
   {coordinate representation, vector space}%
\Index{texBasis}%C%285
   {coordinate vector space}%
   {coordinate vector space}%
\Index{texBasis}%C%291
   {coordinates of geometrical object}%
   {coordinates of geometrical object, vector space}%
\Index{texDrcBasis}%C%90
   {coordinates of geometrical object}%
   {}%
\SubIndex{texDrcBasis}%23
   {in coordinate vector space}%
   {coordinates of geometrical object, coordinate vector space}%
\SubIndex{texDrcBasis}%22
   {in vector space}%
   {coordinates of geometrical object, vector space}%
\Index{texBasis}%C%289
   {coordinates of geometrical object in coordinate representation}%
   {coordinates of geometrical object, coordinate vector space}%
\Index{texDrcBasis}%C%93
   {coordinates of representation}%
   {coordinates of representation}%
\Index{texBasis}%C%265
   {coordinates of representation}%
   {coordinates of representation}%
\Index{texVectorSpace}%C%91
   {coordinates of set of vectors in \dcr vector space}%
   {coordinates of set of vectors, dcr vector space}%
\Index{texVectorSpace}%C%92
   {coordinates of set of vectors in \drc vector space}%
   {coordinates of set of vectors, drc vector space}%
\Index{texVectorField}%C%378
   {coordinates of vector field in \Drc basis}%
   {coordinates of vector field in drc basis}%
\Index{texVectorSpace}%C%22
   {coordinates of vector in \drc basis}%
   {coordinates of vector in drc basis}%
\Index{texBundleRelation}%C%340
   {correspondence continuous on the set}%
   {correspondence continuous on the set}%
\Index{texBundleRelation}%C%297
   {correspondence of homomorphism}%
   {correspondence of homomorphism}%
\Index{texFiberedGroup}%C%318
   {covariant \Ts representation of fibered group}%
   {Tstar covariant representation of fibered group}%
\Index{texBiring}%C%6
   {\CR inverse element of biring}%
   {cr-inverse element}%
\Index{texVectorSpace}%C%5
   {\CR matrix group}%
   {cr-matrix group}%
\Index{texBiring}%C%8
   {\CR power}%
   {cr power}%
\Index{texBiring}%C%7
   {\CR product of matrices}%
   {cr-product of matrices}%
\Index{texVectorSpace}%C%9
   {\crd vector space}%
   {crd vector space}%
\SetIndexSpace%
\Index{texCalculus}%D%96
   {$D$\Hyph valued variable}%
   {D valued variable}%
\Index{texVectorSpace}%D%11
   {\dcr basis of \subs rows vector space}%
   {dcr basis, c rows vector space}%
\Index{texVectorSpace}%D%12
   {\dcr vector}%
   {dcr vector}%
\Index{texVectorSpace}%D%13
   {\dcr vector space}%
   {dcr vector space}%
\Index{texBiring}%D%134
   {determinant of matrix}%
   {determinant}%
\Index{texTidal}%D%137
   {deviation of trajectories}%
   {deviation of trajectories}%
\Index{texBundleRelation}%D%113
   {diagonal in bundle}%
   {diagonal in bundle}%
\Index{texBundleRelation}%D%227
   {diagram of correspondences}%
   {diagram of correspondences}%
\Index{texCalculus}%D%99
   {differentiable functions of \drc vector space to skew field $D$}%
   {differentiable functions, drc vector space to skew field}%
\Index{texCalculus}%D%105
   {differential of mapping of normed \drc vector space to valued skew field}%
   {differential, drc vector space to skew field}%
\Index{texVectorSpace}%D%183
   {dimension of \drc vector space}%
   {dimension of vector space}%
\Index{texFiberedGroup}%D%179
   {direct product of representations of fibered group}%
   {direct product of representations of fibered group}%
\Index{texRepresentation}%D%254
   {direct product of representations of group}%
   {direct product of representations of group}%
\Index{texTstarRepresentation}%D%178
   {direct product of \Ts representations of group}%
   {direct product of representations of group}%
\Index{texLieRepresentation}%D%177
   {direct sum of representations}%
   {direct sum of representations}%
\Index{texVectorField}%D%368
   {distributive law for $\mathcal D\star$\Hyph vector fields}%
   {distributive law, Dstar vector fields}%
\Index{texVectorSpace}%D%81
   {distributive law for $D\star$\Hyph vector space}%
   {distributive law, Dstar vector space}%
\Index{texVectorSpace}%D%55
   {\drc automorphism of vector space}%
   {automorphism of vector space}%
\Index{texVectorSpace}%D%375
   {\drc basis for \sups rows vector space}%
   {drc basis, r rows vector space}%
\Index{texVectorField}%D%376
   {\Drc basis for vector  bundle}%
   {drc basis, vector bundle}%
\Index{texVectorSpace}%D%17
   {\drc basis for vector space}%
   {drc basis, vector space}%
\Index{texVectorSpace}%D%83
   {\drc isomorphism of vector spaces}%
   {isomorphism of vector spaces}%
\Index{texVectorField}%D%383
   {\Drc linear map of vector bundles}%
   {drc linear map of vector bundles}%
\Index{texDrcMorphism}%D%25
   {\drc linear map of vector spaces}%
   {drc linear map of vector spaces}%
\Index{texVectorSpace}%D%23
   {\drc linear span in vector space}%
   {linear span, vector space}%
\Index{texVectorField}%D%373
   {\Drc linearly dependent vector fields}%
   {linearly dependent vector fields}%
\Index{texVectorSpace}%D%15
   {\drc linearly dependent vectors}%
   {linearly dependent, vector space}%
\Index{texVectorField}%D%374
   {\Drc linearly independent vector fields}%
   {linearly independent vector fields}%
\Index{texVectorSpace}%D%114
   {\drc linearly independent vectors}%
   {linearly independent, vector space}%
\Index{texVectorSpace}%D%193
   {\drc system of linear equations}%
   {system of linear equations}%
\Index{texVectorSpace}%D%18
   {\drc vector}%
   {drc vector}%
\Index{texCalculus}%D%46
   {\drc vector function}%
   {drc vector function}%
\Index{texVectorSpace}%D%19
   {\drc vector space}%
   {drc vector space}%
\Index{texVectorField}%D%366
   {$\mathcal D\star$\Hyph vector bundle}%
   {Dstar vector bundle}%
\Index{texVectorField}%D%365
   {$\mathcal D\star$\Hyph vector field}%
   {Dstar vector field}%
\Index{texVectorSpace}%D%27
   {$D\star$\hyph  vector space}%
   {Dstar vector space}%
\Index{texVectorField}%D%372
   {$\mathcal D\star$\hyph linear composition of vector fields}%
   {linear composition of vector fields}%
\Index{texVectorField}%D%370
   {$\mathcal D\star$\hyph product of vector field over scalar}%
   {Dstar product of vector field over scalar, vector space}%
\Index{texVectorSpace}%D%26
   {$D\star$\hyph product of vector over scalar}%
   {Dstar product of vector over scalar, vector space}%
\Index{texBiring}%D%169
   {duality principle for biring}%
   {duality principle for biring}%
\Index{texBiring}%D%170
   {duality principle for biring of matrices}%
   {duality principle for biring of matrices}%
\SetIndexSpace%
\Index{texTstarMorphism}%E%306
   {effective representation of $\mathfrak{F}$\Hyph algebra $A$}%
   {effective representation of algebra}%
\Index{texFiberedAlgebra}%E%326
   {effective representation of fibered $\mathfrak{F}$\Hyph algebra}%
   {effective representation of fibered F-algebra}%
\Index{texFiberedGroup}%E%315
   {effective \Ts representation of fibered group}%
   {effective representation of fibered group}%
\Index{texVectorField}%E%364
   {effective \Ts representation of fibered skew field}%
   {effective representation of fibered skew field}%
\Index{texRepresentation}%E%256
   {effective representation of group}%
   {effective representation of group}%
\Index{texVectorSpace}%E%210
   {effective representation of skew field}%
   {effective representation of skew field}%
\Index{texELie}%E%28
   {enhanced Lie group}%
   {enhanced Lie group}%
\Index{texDiffEq}%E%29
   {essential parameters}%
   {essential parameters}%
\Index{texBundleRelation}%E%235
   {extension of correspondence}%
   {extension of correspondence}%
\Index{texAffine}%E%209
   {extreme line}%
   {extreme line}%
\SetIndexSpace%
\Index{texVectorField}%F%379
   {fibered coordinate \Drc isomorphism}%
   {fibered coordinate drc isomorphism}%
\Index{texBundleRelation}%F%223
   {fibered correspondence from $\mathcal{A}$ to $\mathcal{B}$}%
   {fibered correspondence from A to B}%
\Index{texBundleRelation}%F%224
   {fibered correspondence in $\mathcal{A}$}%
   {fibered correspondence in A}%
\Index{texBundleRelation}%F%330
   {fibered correspondence of homomorphism}%
   {fibered correspondence of homomorphism}%
\Index{texBundleRelation}%F%238
   {fibered equivalence}%
   {fibered equivalence}%
\Index{texFiberedAlgebra}%F%184
   {fibered $\mathfrak{F}$\Hyph algebra}%
   {fibered F-algebra}%
\Index{texFiberedAlgebra}%F%186
   {fibered $\mathfrak{F}$\Hyph subalgebra}%
   {fibered F-subalgebra}%
\Index{texFiberedAlgebra}%F%185
   {fibered group}%
   {fibered group}%
\Index{texFiberedMorphism}%F%338
   {fibered identification morphism}%
   {fibered identification morphism}%
\Index{texFiberedMorphism}%F%346
   {fibered little group}%
   {fibered little group}%
\Index{texBundle}%F%353
   {fibered morphism from bundle $\mathcal{A}$ into $\mathcal{B}$}%
   {fibered morphism from A into B}%
\Index{texFiberedMorphism}%F%337
   {fibered natural morphism}%
   {fibered natural morphism}%
\Index{texBundleRelation}%F%237
   {fibered ordering}%
   {fibered ordering}%
\Index{texBundleRelation}%F%145
   {fibered preordering}%
   {fibered preordering}%
\Index{texFiberedAlgebra}%F%187
   {fibered ring}%
   {fibered ring}%
\Index{texFiberedMorphism}%F%347
   {fibered stability group}%
   {fibered stability group}%
\Index{texBundle}%F%219
   {fibered subset}%
   {fibered subset}%
\Index{texNewton}%F%202
   {field-strength tensor}%
   {field-strength tensor}%
\Index{texBundleRelation}%F%341
   {filter $\mathfrak{F}$ converges to $A$}%
   {filter converges}%
\Index{texNewton}%F%142
   {first Newton law}%
   {First Newton law}%
\Index{texFiberedMorphism}%F%348
   {free \Ts representation of fibered group}%
   {free representation of fibered group}%
\Index{texTstarRepresentation}%F%345
   {free \Ts representation of group}%
   {free representation of group}%
\Index{texAffine}%F%144
   {Frenet transport}%
   {Frenet transport}%
\Index{texCalculus}%F%68
   {function continuous with respect to set of arguments}%
   {function continuous with respect to set of arguments}%
\Index{texCalculus}%F%150
   {function of $\gi n$ $D$\Hyph valued variables}%
   {function of n D valued variables}%
\SetIndexSpace%
\Index{texRefernceFrame}%G%205
   {$G$\Hyph reference frame}%
   {G reference frame}%
\Index{texTypeBasis}%G%30
   {\Gbasis}%
   {G-basis}%
\Index{texBasis}%G%267
   {\Gbasis\ of vector space}%
   {G-basis}%
\Index{texBasis}%G%270
   {\Gcoords\ of basis}%
   {G-coordinates}%
\Index{texTypeBasis}%G%31
   {\Gcoords}%
   {G-coordinates}%
\Index{texTypeBasis}%G%32
   {\Gspace}%
   {GSpace}%
\Index{texDrcBasis}%G%73
   {geometrical object}%
   {}%
\SubIndex{texDrcBasis}%14
   {defined in vector space}%
   {geometrical object, vector space}%
\SubIndex{texDrcBasis}%13
   {in coordinate representation defined in vector space}%
   {geometrical object, coordinate vector space}%
\SubIndex{texDrcBasis}%15
   {of type $A$}%
   {geometrical object of type A, vector space}%
\Index{texBasis}%G%288
   {geometrical object in coordinate representation}%
   {geometrical object, coordinate vector space}%
\Index{texBasis}%G%290
   {geometrical object in vector space}%
   {geometrical object, vector space}%
\Index{texBasis}%G%293
   {geometrical object of type $A$ in vector space}%
   {geometrical object of type A, vector space}%
\Index{texGroupRing}%G%79
   {group algebra}%
   {group algebra}%
\Index{texBasis}%G%271
   {\Gspace}%
   {GSpace}%
\SetIndexSpace%
\Index{texBiring}%H%129
   {Hadamard inverse of matrix}%
   {Hadamard inverse of matrix}%
\Index{texRefernceFrame}%H%195
   {holonomic coordinates of connection}%
   {holonomic coordinates of connection}%
\Index{texRefernceFrame}%H%74
   {holonomic coordinates of vector}%
   {vector holonomic coordinates}%
\Index{texFiberedGroup}%H%132
   {homogeneous bundle of fibered group}%
   {homogeneous bundle of fibered group}%
\Index{texTstarRepresentation}%H%131
   {homogeneous space of group}%
   {homogeneous space of group}%
\Index{texRepresentation}%H%259
   {homogeneous space of group}%
   {homogeneous space of group}%
\Index{texFiberedAlgebra}%H%75
   {homomorphism of fibered $\mathfrak{F}$\Hyph algebras}%
   {homomorphism of fibered F-algebras}%
\Index{texFiberedGroup}%H%76
   {homomorphism of fibered groups}%
   {homomorphism of fibered groups}%
\SetIndexSpace%
\Index{texLieRepresentation}%I%64
   {infinitesimal generator}%
   {infinitesimal generator}%
\Index{texLinearLie}%I%85
   {infinitesimal generators of group Lie}%
   {infinitesimal generators of group Lie}%
\Index{texDrcBasis}%I%171
   {invariance principle}%
   {invariance principle}%
\Index{texBasis}%I%296
   {invariance principle in vector space}%
   {invariance principle, vector space}%
\Index{texBundleRelation}%I%115
   {inverse fibered correspondence}%
   {inverse fibered correspondence}%
\Index{texBundleRelation}%I%121
   {inverse reduced fibered correspondence}%
   {inverse reduced fibered correspondence}%
\Index{texFiberedAlgebra}%I%84
   {isomorphism of fibered $\mathfrak{F}$\Hyph algebras}%
   {isomorphism of fibered F-algebras}%
\SetIndexSpace%
\Index{texFiberedGroup}%K%212
   {kernel of inefficiency of representation of fibered group}%
   {kernel of inefficiency of representation of fibered group}%
\Index{texRepresentation}%K%255
   {kernel of inefficiency of representation of group}%
   {kernel of inefficiency of representation of group}%
\Index{texTstarRepresentation}%K%211
   {kernel of inefficiency of \Ts representation of group $G$}%
   {kernel of inefficiency of representation of group}%
\Index{texDiffProperty}%K%206
   {Killing equation}%
   {Killing equation}%
\Index{texDiffProperty}%K%207
   {Killing equation of second type}%
   {Killing equation second type}%
\Index{texDiffProperty}%K%69
   {Killing vector of second type}%
   {Killing vector second type}%
\Index{texBiring}%K%191
   {Kronecker symbol}%
   {Kronecker symbol}%
\SetIndexSpace%
\Index{texLie}%L%101
   {left invariant vector field}%
   {left invariant vector}%
\Index{texVectorSpace}%L%214
   {left module}%
   {left module}%
\Index{texTstarRepresentation}%L%108
   {left shift}%
   {left shift}%
\Index{texFiberedGroup}%L%109
   {left shift on fibered group}%
   {Tstar shift, fibered group}%
\Index{texRepresentation}%L%251
   {left shift on group}%
   {left shift, group}%
\Index{texLie}%L%103
   {left structural constant of Lie algebra}%
   {left structural constant of Lie algebra}%
\Index{texVectorSpace}%L%98
   {left vector space}%
   {left vector space}%
\Index{texRepresentation}%L%249
   {left-side contravariant representation of group}%
   {left-side contravariant representation}%
\Index{texRepresentation}%L%247
   {left-side covariant representation of group}%
   {left-side covariant representation}%
\Index{texTstarMorphism}%L%302
   {left-side representation of $\mathfrak{F}$\Hyph algebra $A$ in set $M$}%
   {left-side representation of algebra}%
\Index{texFiberedAlgebra}%L%322
   {left-side representation of fibered $\mathfrak{F}$\Hyph algebra}%
   {left-side representation of fibered F-algebra}%
\Index{texRepresentation}%L%244
   {left-side representation of group}%
   {left-side representation of group}%
\Index{texTstarMorphism}%L%164
   {left-side transformation}%
   {left-side transformation}%
\Index{texFiberedAlgebra}%L%328
   {left-side transformation on bundle}%
   {left-side transformation of bundle}%
\Index{texLie}%L%58
   {Lie algebra of Lie group}%
   {algebra Lie group Lie}%
\SubIndex{texLie}%9
   {left defined}%
   {left defined algebra Lie}%
\SubIndex{texLie}%10
   {right defined}%
   {right defined algebra Lie}%
\Index{texDiffProperty}%L%175
   {Lie derivative}%
   {Lie derivative}%
\SubIndex{texDiffProperty}%73
   {of connection}%
   {Lie derivative of connection}%
\SubIndex{texDiffProperty}%72
   {of metric}%
   {Lie derivative of metric}%
\Index{texLie}%L%63
   {Lie group basic operators}%
   {Lie group basic operators}%
\Index{texBundleRelation}%L%232
   {lift of correspondence}%
   {lift of correspondence}%
\Index{texBundle}%L%116
   {lift of map}%
   {lift of map}%
\Index{texBundleRelation}%L%342
   {limit of correspondence with respect to the filter}%
   {limit of correspondence with respect to the filter}%
\Index{texBundleRelation}%L%335
   {limit of filter}%
   {limit of filter}%
\Index{texBundleRelation}%L%334
   {limit set of filter}%
   {limit set of filter}%
\Index{texRepresentation}%L%240
   {linear representation of group}%
   {linear representation of group}%
\Index{texTstarRepresentation}%L%343
   {little group}%
   {little group}%
\Index{texRefernceFrame}%L%117
   {local reference frame}%
   {local reference frame}%
\Index{texBundle}%L%363
   {locally compact at point $p$ space}%
   {locally compact at point space}%
\Index{texBundle}%L%362
   {locally compact space}%
   {locally compact space}%
\Index{texRefernceFrame}%L%165
   {Lorentz transformation}%
   {Lorentz transformation}%
\SetIndexSpace%
\Index{texPolyvector}%M%361
   {$m$\Hyph dimensional parallelepiped}%
   {m dimensional parallelepiped}%
\Index{texPolyvector}%M%24
   {$m$\Hyph vector}%
   {m-vector}%
\Index{texDrcReferenceFrame}%M%138
   {map of type $G$ on manifold}%
   {map of type G on manifold}%
\Index{texCartesian}%M%333
   {mapping space}%
   {mapping space}%
\Index{texDrcMorphism}%M%118
   {matrix of \drc linear map}%
   {matrix of drc linear map}%
\Index{texVectorField}%M%384
   {matrix of fibered \Drc linear map}%
   {matrix of fibered drc linear map}%
\Index{texGeomObject}%M%119
   {metric-affine manifold}%
   {metric-affine manifold}%
\Index{texFiberedMorphism}%M%356
   {morphism of fibered \Ts representations from $\mathcal{F}$ into $\mathcal{G}$}%
   {morphism of fibered representations from f into g}%
\Index{texTstarMorphism}%M%307
   {morphism of \Ts representations from $f$ into $g$}%
   {morphism of representations from f into g}%
\Index{texTstarMorphism}%M%298
   {morphism of \Ts representations of $\mathfrak{F}$ algebra}%
   {morphism of representations of F algebra}%
\Index{texTstarMorphism}%M%300
   {morphism of \Ts representations of $\mathfrak{F}$\Hyph algebra in $\mathfrak{H}$\Hyph algebra}%
   {morphism of representations of F algebra in H algebra}%
\Index{texFiberedMorphism}%M%355
   {morphism of \Ts representations of fibered $\mathfrak{F}$\Hyph algebra}%
   {morphism of representations of fibered F algebra}%
\Index{}%M%283
   {movement on basis manifold}%
   {movement transformation}%
\SetIndexSpace%
\Index{texPolymodule}%N%231
   {$(n)$\hyph vector space}%
   {(n)-vector space}%
\Index{texBundleRelation}%N%218
   {$n$\Hyph ary fibered relation}%
   {fibered relation}%
\Index{texGeomObject}%N%128
   {nonmetricity}%
   {nonmetricity}%
\Index{texVectorSpace}%N%124
   {nonsingular \drc system of linear equations}%
   {nonsingular system of linear equations}%
\Index{texRepresentation}%N%241
   {nonsingular \Ts transformation}%
   {nonsingular transformation}%
\Index{texCalculus}%N%97
   {norm on \drc vector space}%
   {norm on drc vector space}%
\Index{texCalculus}%N%104
   {normed \drc vector space}%
   {normed drc vector space}%
\SetIndexSpace%
\Index{texFiberedAlgebra}%O%133
   {operation on bundle}%
   {operation on bundle}%
\Index{texBundleRelation}%O%236
   {opposite fibered preordering}%
   {opposite fibered preordering}%
\Index{texFiberedGroup}%O%136
   {orbit of representation of fibered group}%
   {orbit of representation of fibered group}%
\Index{texRepresentation}%O%253
   {orbit of representation of group}%
   {orbit of representation of group}%
\Index{texTstarRepresentation}%O%135
   {orbit of \Ts representation of group}%
   {orbit of representation of group}%
\Index{texBasis}%O%281
   {orthonornal basis}%
   {Orthonornal Basis}%
\SetIndexSpace%
\Index{texGeomObject}%P%139
   {parallelogram}%
   {parallelogram}%
\Index{texCalculus}%P%106
   {partial derivative of mapping $f$ with respect to variable $v^{\gi a}$}%
   {partial derivative of mapping with respect to variable, skew field}%
\Index{texCalculus}%P%107
   {partial derivative of mapping $\Vector f$ with respect to variable $v^{\gi a}$}%
   {partial derivative of mapping with respect to variable, drc vector space}%
\Index{texBasis}%P%273
   {passive representation}%
   {passive representation}%
\Index{texBasis}%P%272
   {passive transformation on basis manifold}%
   {passive transformation}%
\Index{texDrcBasis}%P%141
   {passive \Ts representation}%
   {passive representation}%
\Index{texRefernceFrame}%P%182
   {pfaffian derivative}%
   {pfaffian derivative}%
\Index{texPolyvector}%P%358
   {polyvector}%
   {polyvector}%
\Index{texNewton}%P%146
   {potential energy}%
   {potential energy}%
\Index{texDrcBasis}%P%173
   {product of geometrical object and constant}%
   {product of geometrical object and constant}%
\Index{texBasis}%P%295
   {product of geometrical object and constant in vector space}%
   {product of geometrical object and constant, vector space}%
\Index{texTstarMorphism}%P%299
   {product of morphisms of \Ts representations of $\mathfrak{F}$\Hyph algebra}%
   {product of morphisms of representations of F algebra}%
\Index{texVectorField}%P%382
   {product of morphisms of \Ts representations of fibered $\mathfrak{F}$\Hyph algebra}%
   {product of morphisms of representations of fibered F algebra}%
\Index{texBundle}%P%354
   {projection of bundle $\mathcal{E}$ along fiber $E$}%
   {projection of bundle along fiber}%
\SetIndexSpace%
\Index{texBasis}%Q%278
   {quasi affine transformation on basis manifold}%
   {quasi affine transformation}%
\Index{texBasis}%Q%284
   {quasi movement on basis manifold}%
   {quasi movement}%
\Index{texFiberedMorphism}%Q%336
   {quotient bundle}%
   {quotient bundle}%
\SetIndexSpace%
\Index{texBiring}%R%35
   {\sups row of matrix}%
   {r row}%
\Index{texVectorSpace}%R%20
   {\sups rows \drc vector space}%
   {sups rows drc vector space}%
\Index{texBiring}%R%47
   {$r$\hyph row of matrix}%
   {r-row}%
\Index{texBiring}%R%41
   {\RC inverse element of biring}%
   {rc-inverse element}%
\Index{texVectorSpace}%R%37
   {\RC major minor}%
   {RC-major minor}%
\Index{texVectorSpace}%R%39
   {\RC matrix group}%
   {rc-matrix group}%
\Index{texVectorSpace}%R%40
   {\RC nonsingular matrix}%
   {RC nonsingular matrix}%
\Index{texBiring}%R%44
   {\RC power}%
   {rc power}%
\Index{texBiring}%R%42
   {\RC product of matrices}%
   {rc-product of matrices}%
\Index{texBiring}%R%38
   {\RC quasideterminant}%
   {RC-quasideterminant}%
\Index{texVectorSpace}%R%43
   {\RC rank of matrix}%
   {rc-rank of matrix}%
\Index{texVectorSpace}%R%36
   {\RC singular matrix}%
   {RC singular matrix}%
\Index{texVectorSpace}%R%45
   {\rcd vector space}%
   {rcd vector space}%
\Index{texCartesian}%R%221
   {reduced Cartesian product of bundles}%
   {reduced Cartesian product of bundles}%
\Index{texCartesian}%R%222
   {reduced Cartesian product of total spaces}%
   {reduced Cartesian product of total spaces}%
\Index{texBundleRelation,texBundleRelation}%R%225
   {reduced fibered correspondence from $\mathcal{A}$ to $\mathcal{B}$}%
   {reduced fibered correspondence from A to B}%
\Index{texBundleRelation}%R%226
   {reduced fibered correspondence in $\mathcal{A}$}%
   {reduced fibered correspondence in A}%
\Index{texBiring}%R%168
   {reducible biring}%
   {reducible biring}%
\Index{texRefernceFrame}%R%194
   {reference frame in event space}%
   {reference frame in event space}%
\Index{texDrcReferenceFrame}%R%123
   {reference frame manifold}%
   {reference frame manifold}%
\Index{texBundleRelation}%R%157
   {reflexive $2$\Hyph ary fibered relation}%
   {reflexive 2 ary fibered relation}%
\Index{texTstarRepresentation}%R%161
   {representation of group}%
   {}%
\SubIndex{texTstarRepresentation}%43
   {contravariant \Ts}%
   {Tstar contravariant representation of group}%
\SubIndex{texTstarRepresentation}%42
   {covariant \Ts}%
   {Tstar covariant representation of group}%
\SubIndex{texDrcBasis}%38
   {\drc linear \sT}%
   {linear representation of group}%
\SubIndex{texTstarRepresentation}%44
   {effective}%
   {effective representation of group}%
\SubIndex{texDrcBasis}%39
   {\rcd}%
   {rcd linear representation of group}%
\SubIndex{texTstarRepresentation}%40
   {\sT}%
   {starT representation of group}%
\SubIndex{texTstarRepresentation}%41
   {\Ts}%
   {Tstar representation of group}%
\Index{texRepresentation}%R%246
   {representation of group}%
   {representation of group}%
\Index{texDrcBasis}%R%159
   {representative of geometrical object in vector space}%
   {representative of geometrical object, vector space}%
\Index{texBasis}%R%292
   {representative of geometrical object in vector space}%
   {representative of geometrical object, vector space}%
\Index{texBundleRelation}%R%234
   {restriction of correspondence $\Phi$ to set $C$}%
   {restriction of correspondence}%
\Index{texLie}%R%151
   {right invariant vector field}%
   {right invariant vector}%
\Index{texVectorSpace}%R%215
   {right module}%
   {right module}%
\Index{texTstarRepresentation}%R%158
   {right shift}%
   {right shift}%
\Index{texRepresentation}%R%252
   {right shift on group}%
   {right shift, group}%
\Index{texLie}%R%153
   {right structural constant of Lie algebra}%
   {right structural constant of Lie algebra}%
\Index{texVectorSpace}%R%149
   {right vector space}%
   {right vector space}%
\Index{texRepresentation}%R%250
   {right-side contravariant representation of group}%
   {right-side contravariant representation}%
\Index{texRepresentation}%R%248
   {right-side covariant representation of group}%
   {right-side covariant representation}%
\Index{texTstarMorphism}%R%304
   {right-side representation of $\mathfrak{F}$\Hyph algebra $A$ in set $M$}%
   {right-side representation of algebra}%
\Index{texFiberedAlgebra}%R%324
   {right-side representation of fibered $\mathfrak{F}$\Hyph algebra}%
   {right-side representation of fibered F-algebra}%
\Index{texRepresentation}%R%245
   {right-side representation of group}%
   {right-side representation of group}%
\Index{texTstarMorphism}%R%242
   {right-side transformation}%
   {right-side transformation}%
\Index{texRepresentation}%R%243
   {right-side transformation}%
   {right-side transformation}%
\Index{texRepresentation}%R%260
   {row vector}%
   {row vector}%
\Index{texVectorSpace}%R%213
   {$R\star$\Hyph module}%
   {Rstar-module}%
\SetIndexSpace%
\Index{texNewton}%S%196
   {scalar potential}%
   {scalar potential}%
\Index{texNewton}%S%72
   {second Newton law}%
   {Second Newton law}%
\Index{texPolyvector}%S%359
   {simple polyvector}%
   {simple polyvector}%
\Index{texTstarMorphism}%S%303
   {single transitive representation of $\mathfrak{F}$\Hyph algebra $A$}%
   {single transitive representation of algebra}%
\Index{texFiberedAlgebra}%S%323
   {single transitive representation of fibered $\mathfrak{F}$\Hyph algebra}%
   {single transitive representation of fibered F-algebra}%
\Index{texRepresentation}%S%258
   {single transitive representation of group}%
   {single transitive representation of group}%
\Index{texPolyvector}%S%360
   {skew product of vectors}%
   {skew product of vectors}%
\Index{texTstarRepresentation}%S%352
   {space of orbits of \Ts representation}%
   {space of orbits of Ts representation}%
\Index{texTidal}%S%197
   {speed of deviation}%
   {speed of deviation}%
\Index{texVectorField}%S%381
   {$(\mathcal S\RCstar,\mathcal T\RCstar)$\Hyph linear map of vector bundles}%
   {src trc linear map of vector bundles}%
\Index{texDrcMorphism}%S%314
   {$(S\RCstar,T\RCstar)$\Hyph linear map of vector spaces}%
   {src trc linear map of vector spaces}%
\Index{texTstarRepresentation}%S%344
   {stability group}%
   {stability group}%
\Index{texDrcBasis}%S%199
   {standard coordinates of basis}%
   {standard coordinates of basis}%
\Index{texBasis}%S%269
   {standard coordinates of basis}%
   {standard coordinates of basis}%
\Index{texBiring}%S%198
   {standard representation of matrix}%
   {Standard representation}%
\Index{texVectorSpace}%S%217
   {$\star D$\hyph vector space}%
   {starD-vector space}%
\Index{texLinearMap}%S%48
   {$\star D$\Hyph product of \drc linear map $A$ over scalar}%
   {starD product of drc linear map over scalar}%
\Index{texVectorSpace}%S%216
   {$\star R$\hyph module}%
   {starR-module}%
\Index{texTstarRepresentation}%S%49
   {\sT shift}%
   {starT shift}%
\Index{texFiberedGroup}%S%50
   {\sT shift on fibered group}%
   {starT shift, fibered group}%
\Index{texTstarMorphism}%S%160
   {\sT representation of $\mathfrak{F}$\Hyph algebra $A$ in set $M$}%
   {starT representation of algebra}%
\Index{texFiberedAlgebra}%S%162
   {\sT representation of fibered $\mathfrak{F}$\Hyph algebra}%
   {starT representation of fibered F-algebra}%
\Index{texFiberedGroup}%S%163
   {\sT representation of fibered group}%
   {starT representation of fibered group}%
\Index{texTstarMorphism}%S%308
   {\sT transformation}%
   {starT transformation}%
\Index{texFiberedAlgebra}%S%167
   {\sT transformation on bundle}%
   {starT transformation of bundle}%
\SubIndex{}%71
   {nonsingular}%
   {nonsingular transformation of bundle}%
\Index{texBundle}%S%220
   {subbundle}%
   {subbundle}%
\Index{texVectorField}%S%371
   {subbundle of $\mathcal D\star$\hyph vector space}%
   {subbundle of Dstar vector bundle}%
\Index{texLinearMap}%S%200
   {sum of \drc linear maps}%
   {sum of drc linear maps, drc vector spaces}%
\Index{texDrcBasis}%S%201
   {sum of geometrical objects}%
   {sum of geometrical objects}%
\Index{texBasis}%S%294
   {sum of geometrical objects in vector space}%
   {sum of geometrical objects, vector space}%
\Index{texBundleRelation}%S%148
   {symmetric $2$\Hyph ary fibered relation}%
   {symmetric 2 ary fibered relation}%
\Index{texBasis}%S%264
   {symmetry group}%
   {symmetry group}%
\Index{texDrcBasis}%S%78
   {symmetry group}%
   {SymmetryGroup}%
\Index{texGenRelativity}%S%192
   {synchronization of reference frame}%
   {synchronization of reference frame}%
\SetIndexSpace%
\Index{texLie}%T%203
   {tensor product of representations}%
   {tensor product of representations}%
\Index{texCalculus}%T%152
   {topological \drc vector space}%
   {topological drc vector space}%
\Index{texCalculus}%T%33
   {topological skew field}%
   {topological skew field}%
\Index{texGeomObject}%T%239
   {torsion form}%
   {torsion form}%
\Index{texGeomObject}%T%95
   {torsion tensor}%
   {torsion tensor}%
\Index{texFiberedMorphism}%T%351
   {tower of bundles}%
   {tower of bundles}%
\Index{texTstarMorphism}%T%313
   {transformation of set}%
   {transformation of set}%
\Index{texDrcBasis}%T%166
   {transformation on basis manifold}%
   {}%
\SubIndex{texDrcBasis}%62
   {active}%
   {active transformation, vector space}%
\SubIndex{texTypeBasis}%63
   {affine}%
   {affine transformation}%
\SubIndex{texTypeBasis}%64
   {movement}%
   {movement transformation}%
\SubIndex{texDrcBasis}%67
   {passive}%
   {passive transformation, vector space}%
\SubIndex{texTypeBasis}%65
   {quasi affine}%
   {quasi affine transformation}%
\SubIndex{texTypeBasis}%66
   {quasi movement}%
   {quasi movement}%
\Index{texFiberedAlgebra}%T%329
   {transformation on bundle}%
   {transformation of bundle}%
\Index{texBundleRelation}%T%147
   {transitive $2$\Hyph ary fibered relation}%
   {transitive 2 ary fibered relation}%
\Index{texTstarMorphism}%T%305
   {transitive representation of $\mathfrak{F}$\Hyph algebra $A$}%
   {transitive representation of algebra}%
\Index{texFiberedAlgebra}%T%325
   {transitive representation of fibered $\mathfrak{F}$\Hyph algebra}%
   {transitive representation of fibered F-algebra}%
\Index{texRepresentation}%T%257
   {transitive representation of group}%
   {transitive representation of group}%
\Index{texTstarMorphism}%T%386
   {\Ts transformation coordinated with equivalence}%
   {transformation coordinated with equivalence}%
\Index{texVectorSpace}%T%52
   {\Ts linear composition of  vectors}%
   {linear composition of  vectors}%
\Index{texVectorSpace}%T%51
   {\Ts matrices vector space}%
   {matrices vector space}%
\Index{texTstarRepresentation}%T%53
   {\Ts shift}%
   {Tstar shift}%
\Index{texTstarMorphism}%T%301
   {\Ts representation of $\mathfrak{F}$\Hyph algebra $A$ in set $M$}%
   {Tstar representation of algebra}%
\Index{texFiberedAlgebra}%T%321
   {\Ts representation of fibered $\mathfrak{F}$\Hyph algebra}%
   {Tstar representation of fibered F-algebra}%
\Index{texFiberedGroup}%T%316
   {\Ts representation of fibered group}%
   {Tstar representation of fibered group}%
\Index{texTstarMorphism}%T%309
   {\Ts transformation}%
   {Tstar transformation}%
\Index{texFiberedAlgebra}%T%327
   {\Ts transformation on bundle}%
   {Tstar transformation of bundle}%
\Index{texFiberedGroup}%T%319
   {twin representations of fibered group}%
   {twin representations of fibered group}%
\Index{texTstarRepresentation}%T%140
   {twin representations of group}%
   {twin representations of group}%
\Index{texLinearMap}%T%320
   {twin representations of skew field}%
   {twin representations of skew field}%
\SetIndexSpace%
\Index{texVectorField}%U%369
   {unitarity law for  $\mathcal D\star$\Hyph vector fields}%
   {unitarity law, Dstar vector fields}%
\Index{texVectorSpace}%U%82
   {unitarity law for $D\star$\Hyph vector space}%
   {unitarity law, Dstar vector space}%
\SetIndexSpace%
\Index{texPolymodule}%V%230
   {($D_1\RCstar$, ..., $D_n\RCstar$)\hyph vector space}%
   {(d1rc,dnrc)-vector space}%
\Index{texPolymodule}%V%229
   {($S\star$, $\star T$)\hyph vector space}%
   {(Sstar,starT)-vector space}%
\Index{texCalculus}%V%102
   {valued skew field}%
   {valued skew field}%
\Index{texFiberedAlgebra}%V%70
   {vector bundle}%
   {vector bundle}%
\Index{texNewton}%V%54
   {vector potential}%
   {vector potential}%
\Index{texVectorSpace}%V%204
   {vector space type}%
   {vector space type}%

\CloseIndex

%auto-ignore
\def\indexname{Special Symbols and Notations}
\OpenIndex

\SetIndexSpace%
\Symb{texBiring}%A/0
   {$(^{\gi a}_{\gi b})$\hyph\CR quasideterminant}%
   {a b CR quasideterminant definition}%
\Symb{texBiring}%A/0
   {$(^{\gi a}_{\gi b})$\hyph \RC quasideterminant}%
   {a b RC-quasideterminant definition}%
\Symb{texBiring}%A/0
   {minor}%
   {A from b a}%
\Symb{texBiring}%A/0
   {minor}%
   {A from columns T}%
\Symb{texBiring}%A/0
   {minor}%
   {A from rows S}%
\Symb{texBiring}%A/0
   {minor}%
   {A without column a}%
\Symb{texBiring}%A/0
   {minor}%
   {A without columns T}%
\Symb{texBiring}%A/0
   {minor}%
   {A without row b}%
\Symb{texBiring}%A/0
   {minor}%
   {A without rows S}%
\Symb{texPolymodule}%A/0
   {active transformation}%
   {active transformation}%
\Symb{texTypeBasis}%A/0
   {affine space}%
   {affine space}%
\Symb{texBasis}%A/0
   {affine space}%
   {An}%
\Symb{texBiring}%A/0
   {\subs row ($c$\hyph row) of matrix}%
   {c row}%
\Symb{texBiring}%A/0
   {\CR power of element $A$ of biring}%
   {cr power}%
\Symb{texBiring}%A/0
   {\CR inverse element of biring}%
   {cr-inverse element}%
\Symb{texBiring}%A/0
   {\CR product of matrices}%
   {cr-product of matrices}%
\Symb{texVectorSpace}%A/0
   {\dcr vector}%
   {dcr vector}%
\Symb{texLie}%A/0
   {derivative of left shift}%
   {derivative of left shift}%
\Symb{texLie}%A/0
   {derivative of left shift}%
   {derivative of left shift, 1-Parameter Group}%
\Symb{texLie}%A/0
   {derivative of right shift}%
   {derivative of right shift}%
\Symb{texLie}%A/0
   {}%
   {derivative of right shift}%
\Symb{texLie}%A/0
   {derivative of right shift}%
   {derivative of right shift, 1-Parameter Group}%
\Symb{texLie}%A/0
   {derivative of left shift}%
   {derivative of Tstar shift}%
\Symb{texVectorSpace}%A/0
   {\drc vector}%
   {drc vector}%
\Symb{texAffine}%A/0
   {derivative}%
   {overline nabla_l, definition 2}%
\Symb{texPolymodule}%A/0
   {passive transformation}%
   {passive transformation}%
\Symb{texBiring}%A/0
   {\sups row ($r$\hyph row) of matrix}%
   {r row}%
\Symb{texBiring}%A/0
   {\RC power of element $A$ of biring}%
   {rc power}%
\Symb{texBiring}%A/0
   {\RC inverse element of biring}%
   {rc-inverse element}%
\Symb{texBiring}%A/0
   {\RC product of matrices}%
   {rc-product of matrices}%
\Symb{texBiring}%A/0
   {\RC quasideterminant}%
   {RC-quasideterminant definition}%
\Symb{texTstarRepresentation}%A/0
   {right shift}%
   {right shift}%
\Symb{texPolyvector}%A/0
   {skew product of vectors $\Vector a_1$, ..., $\Vector a_m$}%
   {skew product of vectors}%
\Symb{texFiberedGroup}%A/0
   {\sT shift}%
   {starT shift, fibered group}%
\Symb{texTstarRepresentation}%A/0
   {left shift}%
   {Tstar shift}%
\Symb{texFiberedGroup}%A/0
   {\Ts shift}%
   {Tstar shift, fibered group}%
\Symb{texRefernceFrame}%A/0
   {anholonomic coordinates of vector}%
   {vector anholonomic coordinates}%
\Symb{texRefernceFrame}%A/0
   {holonomic coordinates of vector}%
   {vector holonomic coordinates}%

\SetIndexSpace%
\Symb{texBasis}%B/0
   {basis manifold of affine space}%
   {BAn}%
\Symb{texBasis}%B/0
   {basis manifold of vector space}%
   {basis manifold of vector space}%
\Symb{texBasis}%B/0
   {basis manifold of vector space $\mathcal{V}$}%
   {basis manifold of vector space}%
\Symb{texBasis}%B/0
   {basis manifold of central affine space}%
   {BCAn}%
\Symb{texBasis}%B/0
   {basis manifold of Euclid space}%
   {BEn}%
\Symb{texBundle}%B/0
   {Cartesian power $\mathcal{A}$ of bundle $\mathcal{B}$}%
   {Cartesian power of bundle}%
\Symb{texBundle}%B/0
   {Cartesian power $A$ of set $B$}%
   {Cartesian power of set}%
\Symb{texTypeBasis}%B/0
   {basis manifold of affine space}%
   {FAn}%
\Symb{texTypeBasis}%B/0
   {basis manifold of central affine space}%
   {FCAn}%
\Symb{texTypeBasis}%B/0
   {basis manifold of Euclid space}%
   {FEn}%

\SetIndexSpace%
\Symb{texBasis}%C/0
   {central affine space}%
   {CAn}%
\Symb{texTypeBasis}%C/0
   {central affine space}%
   {central affine space}%
\Symb{texLie}%C/0
   {left structural constant of Lie algebra}%
   {left structural constant of Lie algebra}%
\Symb{texLie}%C/0
   {right structural constant of Lie algebra}%
   {right structural constant of Lie algebra}%

\SetIndexSpace%
\Symb{texLieRepresentation}%D/0
   {basis vector of \sT representation}%
   {basis vector of starT representation}%
\Symb{texLieRepresentation}%D/0
   {basis vector of \sT representation}%
   {basis vector of starT representation, coordinates}%
\Symb{texLieRepresentation}%D/0
   {basis vector of \Ts representation}%
   {basis vector of Tstar representation}%
\Symb{texLieRepresentation}%D/0
   {basis vector of \Ts representation}%
   {basis vector of Tstar representation, coordinates}%
\Symb{texVectorSpace}%D/0
   {\subs rows \dcr vector space}%
   {c rows dcr vector space}%
\Symb{texVectorField}%D/0
   {coordinate \Drc vector bundle}%
   {coordinate drc vector bundle}%
\Symb{texVectorSpace}%D/0
   {coordinate \drc vector space}%
   {coordinate drc vector space}%
\Symb{texCalculus}%D/0
   {differential of function}%
   {differential, drc vector space to drc vector space}%
\Symb{texCalculus}%D/0
   {differential of function}%
   {differential, drc vector space to skew field}%
\Symb{texVectorSpace}%D/0
   {matrices vector space}%
   {matrices vector space}%
\Symb{texAffine}%D/0
   {Cartan derivative}%
   {overbrace D}%
\Symb{texAffine}%D/0
   {derivative}%
   {overline D}%
\Symb{texCalculus}%D/0
   {partial derivative of mapping $\Vector f$ with respect to variable $v^{\gi a}$}%
   {partial derivative of mapping, 1, drc vector space}%
\Symb{texCalculus}%D/0
   {partial derivative of mapping $f$ with respect to variable $v^{\gi a}$}%
   {partial derivative of mapping, 1, skew field}%
\Symb{texVectorSpace}%D/0
   {\sups rows \drc vector space}%
   {r rows drc vector space}%
\Symb{texTidal}%D/0
   {speed of deviation}%
   {speed of deviation}%
\Symb{texVectorSpace}%D/0
   {vector space type}%
   {vector space type}%

\SetIndexSpace%
\Symb{texTypeBasis}%E/0
   {affine basis}%
   {Affine Basis}%
\Symb{texBasis}%E/0
   {affine basis}%
   {Affine Basis}%
\Symb{texTypeBasis}%E/0
   {basis}%
   {basis}%
\Symb{texBasis}%E/0
   {basis of vector space}%
   {Basis e}%
\Symb{texBasis}%E/0
   {basis in vector space $\mathcal{V}$}%
   {basis in V}%
\Symb{texVectorSpace}%E/0
   {basis of vector space}%
   {basis, vector space}%
\Symb{texPolymodule}%E/0
   {basis of $(n)$\hyph vector space}%
   {basis,n vector space}%
\Symb{texCartesian}%E/0
   {Cartesian power of total spaces}%
   {Cartesian power of total spaces}%
\Symb{texCartesian}%E/0
   {Cartesian product of total spaces}%
   {Cartesian product of total spaces, definition 1}%
\Symb{texBasis}%E/0
   {central affine basis}%
   {Central Affine Basis}%
\Symb{texVectorField}%E/0
   {basis for \Drc vector bundle}%
   {drc basis, vector bundle}%
\Symb{texRefernceFrame}%E/0
   {form of reference frame}%
   {dual forms, reference frame}%
\Symb{texBasis}%E/0
   {Euclid space}%
   {En}%
\Symb{texTypeBasis}%E/0
   {Euclid space}%
   {En}%
\Symb{texTypeBasis}%E/0
   {pseudo Euclid space}%
   {Enm}%
\Symb{texBasis}%E/0
   {pseudo Euclid space}%
   {Enm}%
\Symb{texFiberedAlgebra}%E/0
   {identical transformation of bundle}%
   {identical transformation of bundle}%
\Symb{texBasis}%E/0
   {orthonornal basis}%
   {Orthonornal Basis}%
\Symb{texCartesian}%E/0
   {reduced Cartesian product of total spaces}%
   {reduced Cartesian product of total spaces, definition 1}%
\Symb{texFiberedAlgebra}%E/0
   {set of nonsingular \sT transformations of bundle $\mathcal{E}$}%
   {set of starT nonsingular transformations of bundle}%
\Symb{texFiberedAlgebra}%E/0
   {set of nonsingular \Ts transformations of bundle $\mathcal{E}$}%
   {set of Tstar nonsingular transformations of bundle}%
\Symb{texBasis}%E/0
   {standard coordinates of basis}%
   {standard coordinates of basis}%
\Symb{texRefernceFrame}%E/0
   {standard coordinates of reference frame}%
   {standard coordinates of reference frame}%
\Symb{texRefernceFrame}%E/0
   {vector field of reference frame}%
   {vector field of reference frame}%
\Symb{texBasis}%E/0
   {vector of basis}%
   {vector of basis}%

\SetIndexSpace%
\Symb{texVectorSpace}%F/0
   {coordinates of basis in \subs rows \dcr vector space}%
   {basis coordinates, c rows dcr vector space}%
\Symb{texVectorSpace}%F/0
   {coordinates of basis in \sups rows \drc vector space}%
   {basis coordinates, r rows drc vector space}%
\Symb{texVectorSpace}%F/0
   {basis for \subs rows \dcr vector space}%
   {basis, c rows dcr vector space}%
\Symb{texVectorSpace}%F/0
   {basis for \sups rows \drc vector space}%
   {basis, r rows drc vector space}%
\Symb{texDiffEq}%F/0
   {central affine basis}%
   {Central Affine Basis}%
\Symb{texBundle}%F/0
   {fibered morphism from bundle $\mathcal{A}$ into $\mathcal{B}$}%
   {fibered morphism from A into B}%
\Symb{texBundleRelation}%F/0
   {filter $\mathfrak{F}$ converges to set $A$}%
   {filter converges}%
\Symb{texFiberedAlgebra}%F/0
   {homomorphism of fibered $\mathfrak{F}$\Hyph algebras}%
   {homomorphism of fibered F-algebras}%
\Symb{texBundleRelation}%F/0
   {inverse fibered correspondence}%
   {inverse fibered correspondence, 1}%
\Symb{texBundleRelation}%F/0
   {inverse reduced fibered correspondence}%
   {inverse reduced fibered correspondence, 1}%
\Symb{texCartesian}%F/0
   {map to Cartesian product}%
   {map to Cartesian product}%
\Symb{texTstarRepresentation}%F/0
   {representation orbit of group $G$}%
   {orbit of representation of group}%
\Symb{texTypeBasis}%F/0
   {orthonornal basis}%
   {Orthonornal Basis}%
\Symb{texRefernceFrame}%F/0
   {reference frame}%
   {reference frame}%
\Symb{texRefernceFrame}%F/0
   {reference frame, extensive definition}%
   {reference frame, extensive definition}%
\Symb{texPolymodule}%F/0
   {standard coordinates of basis}%
   {standard coordinates of basis}%
\Symb{texPolymodule}%F/0
   {vector of basis}%
   {vector of basis}%

\SetIndexSpace%
\Symb{texVectorSpace}%G/0
   {\CR matrix group}%
   {cr-matrix group}%
\Symb{texFiberedMorphism}%G/0
   {fibered little group of section $h$}%
   {fibered little group}%
\Symb{texFiberedMorphism}%G/0
   {fibered stability group of section $h$}%
   {fibered stability group}%
\Symb{texLie}%G/0
   {algebra Lie of group Lie}%
   {g}%
\Symb{texLie}%G/0
   {left defined algebra Lie of group Lie}%
   {gl}%
\Symb{texTypeBasis}%G/0
   {affine transformation group}%
   {GLAn}%
\Symb{texBasis}%G/0
   {affine transformation group}%
   {GLAn}%
\Symb{texLie}%G/0
   {right defined algebra Lie of group Lie}%
   {gr}%
\Symb{texBasis}%G/0
   {group of homomorphisms of vector space $\mathcal{V}$}%
   {GV}%
\Symb{texTstarRepresentation}%G/0
   {little group of $x$}%
   {little group}%
\Symb{texFiberedGroup}%G/0
   {orbit of effective covariant \sT representation of fibered group}%
   {orbit of effective starT covariant representation of fibered group}%
\Symb{texTstarRepresentation}%G/0
   {orbit of effective covariant \sT representation of group}%
   {orbit of effective starT covariant representation of group}%
\Symb{texFiberedGroup}%G/0
   {orbit of effective covariant \Ts representation of fibered group}%
   {orbit of effective Tstar covariant representation of fibered group}%
\Symb{texTstarRepresentation}%G/0
   {orbit of effective covariant \Ts representation of group}%
   {orbit of effective Tstar covariant representation of group}%
\Symb{texVectorSpace}%G/0
   {\RC matrix group}%
   {rc-matrix group}%
\Symb{texTstarRepresentation}%G/0
   {stability group of $x$}%
   {stability group}%

\SetIndexSpace%
\Symb{texBiring}%H/0
   {Hadamard inverse of matrix}%
   {Hadamard inverse of matrix}%
\Symb{texLinearMap}%H/0
   {\rcd vector space of \drc linear maps}%
   {rcd vector space of drc linear maps}%

\SetIndexSpace%
\Symb{texLieRepresentation}%I/0
   {infinitesimal generator of representation}%
   {infinitesimal generator of representation}%
\Symb{texLinearLie}%I/0
   {Lie group infinitesimal generators}%
   {Lie group infinitesimal generators}%

\SetIndexSpace%
\Symb{texRepresentation}%L/0
   {left shift}%
   {left shift}%
\Symb{texDiffProperty}%L/0
   {Lie derivative of connection}%
   {Lie derivative of connection}%
\Symb{texDiffProperty}%L/0
   {Lie derivative of metric}%
   {Lie derivative of metric}%
\Symb{texBundleRelation}%L/0
   {limit of correspondence $\Phi$ with respect to the filter $\mathfrak{F}$}%
   {limit of correspondence with respect to the filter}%
\Symb{texBasis}%L/0
   {passive transformation}%
   {passive transformation}%
\Symb{texRepresentation}%L/0
   {set of left-side nonsingular transformations of set $M$}%
   {set of left-side nonsingular transformations}%

\SetIndexSpace%
\Symb{texTstarMorphism}%M/0
   {set of \sT transformations of set $M$}%
   {set of starT transformations}%
\Symb{texTstarMorphism}%M/0
   {set of \Ts transformations of set $M$}%
   {set of Tstar transformations}%
\Symb{texTstarRepresentation}%M/0
   {space of orbits of effective \sT covariant representation of the group}%
   {space of orbits of effective sT representation}%
\Symb{texTstarRepresentation}%M/0
   {space of orbits of effective \Ts covariant representation of the group}%
   {space of orbits of effective Ts representation}%
\Symb{texTstarRepresentation}%M/0
   {space of orbits of \Ts representation $f$ of group $G$ in set $M$}%
   {space of orbits of Ts representation}%

\SetIndexSpace%
\Symb{texBasis}%O/0
   {geometrical object in coordinate representation}%
   {geometrical object, coordinate vector space}%
\Symb{texBasis}%O/0
   {geometrical object}%
   {geometrical object, vector space}%
\Symb{texFiberedGroup}%O/0
   {orbit of representation of fibered group $\mathcal{G}$}%
   {orbit of representation of fibered group}%
\Symb{texRepresentation}%O/0
   {orbit of representation of the group $G$}%
   {orbit of representation of group}%

\SetIndexSpace%
\Symb{texBundle}%P/0
   {bundle}%
   {bundle}%
\Symb{texFiberedMorphism}%P/0
   {bundle of level $2$}%
   {bundle of level 2}%
\Symb{texFiberedMorphism}%P/0
   {bundle of level $n$}%
   {bundle of level n}%
\Symb{texCartesian}%P/0
   {Cartesian power of bundle}%
   {Cartesian power of bundle}%
\Symb{texCartesian}%P/0
   {Cartesian product of bundles}%
   {Cartesian product of bundles, definition 1}%
\Symb{texCartesian}%P/0
   {reduced Cartesian product of bundles}%
   {reduced Cartesian product of bundles, definition 1}%
\Symb{texFiberedAlgebra}%P/0
   {set of nonsingular \sT transformations of bundle $\bundle{}pE{}$}%
   {set of starT nonsingular transformations of bundle, projection}%
\Symb{texFiberedAlgebra}%P/0
   {set of nonsingular \Ts transformations of bundle $\bundle{}pE{}$}%
   {set of Tstar nonsingular transformations of bundle, projection}%

\SetIndexSpace%
\Symb{texBasis}%R/0
   {active transformation}%
   {active transformation}%
\Symb{texAffine}%R/0
   {Cartan curvature}%
   {Cartan curvature}%
\Symb{texVectorSpace}%R/0
   {\CR rank of matrix}%
   {cr-rank of matrix}%
\Symb{texBundleRelation}%R/0
   {diagonal in bundle  $\bundle{}pA{}$}%
   {diagonal in bundle, 2}%
\Symb{texBundleRelation}%R/0
   {diagonal in bundle $\mathcal{A}$}%
   {diagonal in reduced bundle, 2}%
\Symb{texAffine}%R/0
   {curvature}%
   {GLn curvature_overline}%
\Symb{texVectorSpace}%R/0
   {\RC rank of matrix}%
   {rc-rank of matrix}%
\Symb{texRepresentation}%R/0
   {right shift}%
   {right shift}%
\Symb{texRepresentation}%R/0
   {set of right-side nonsingular transformations of set $M$}%
   {set of right-side nonsingular transformations}%

\SetIndexSpace%
\Symb{texBundleRelation}%S/0
   {composition of fibered correspondences}%
   {composition of fibered correspondences}%
\Symb{texBundleRelation}%S/0
   {inverse fibered correspondence}%
   {inverse fibered correspondence, 2}%
\Symb{texBundleRelation}%S/0
   {inverse reduced fibered correspondence}%
   {inverse reduced fibered correspondence, 2}%
\Symb{texVectorSpace}%S/0
   {linear span in vector space}%
   {linear span, vector space}%

\SetIndexSpace%
\Symb{texLie}%T/0
   {tangent plane to group $G$}%
   {TaG}%

\SetIndexSpace%
\Symb{texBasis}%V/0
   {coordinate vector space}%
   {coordinate vector space}%
\Symb{texBasis}%V/0
   {coordinates in vector space}%
   {coordinates in vector space}%
\Symb{texVectorSpace}%V/0
   {\dcr vector space}%
   {left CR vector space}%
\Symb{texVectorSpace}%V/0
   {\drc vector space}%
   {left RC vector space}%
\Symb{texLinearMap}%V/0
   {($S$, $T$)\hyph bimodule}%
   {R S bimodule}%
\Symb{texVectorSpace}%V/0
   {\crd vector space}%
   {right CR vector space}%
\Symb{texVectorSpace}%V/0
   {\rcd vector space}%
   {right RC vector space}%
\Symb{texBasis}%V/0
   {vector space}%
   {V}%

\SetIndexSpace%
\Symb{texPolymodule}%W/0
   {geometrical object in coordinate representation		defined in vector space}%
   {geometrical object, coordinate vector space}%
\Symb{texPolymodule}%W/0
   {geometrical object in vector space}%
   {geometrical object, vector space}%

\SetIndexSpace%
\Symb{texRefernceFrame}%X/0
   {anholonomic coordinate}%
   {x(k)}%

\SetIndexSpace%
\Symb{texBundleRelation}%D/1
   {diagonal in bundle $\mathcal{A}$}%
   {diagonal in bundle, 1}%

\SetIndexSpace%
\Symb{texTidal}%Delta/1
   {deviation of trajectories}%
   {deviation of trajectories}%
\Symb{texRepresentation}%Delta/1
   {identical transformation}%
   {identical transformation}%
\Symb{texTstarMorphism}%Delta/1
   {identical transformation}%
   {identical transformation}%
\Symb{texBasis}%Delta/1
   {image of vector $\Vector e_k\in\Basis e$ under isomorphism to coordinate vector space}%
   {image of vector e_k, coordinate vector space}%
\Symb{texBiring}%Delta/1
   {Kronecker symbol}%
   {Kronecker symbol}%

\SetIndexSpace%
\Symb{texRefernceFrame}%Gamma/1
   {anholonomic coordinates of connection}%
   {anholonomic coordinates of connection}%
\Symb{texAffine}%Gamma/1
   {Cartan symbol}%
   {Cartan symbol}%
\Symb{texAffine}%Gamma/1
   {connection}%
   {conection overline}%
\Symb{texRefernceFrame}%Gamma/1
   {holonomic coordinates of connection}%
   {holonomic coordinates of connection}%
\Symb{texAffine}%Gamma/1
   {Cartan connection}%
   {overbrace Gamma i kl}%
\Symb{texBundle}%Gamma/1
   {set of sections of bundle}%
   {set of sections of bundle}%

\SetIndexSpace%
\Symb{texLie}%Lambda/1
   {inverse operator to operator $\psi_l$}%
   {inverse operator to operator psi l}%
\Symb{texLie}%Lambda/1
   {inverse operator to operator $\psi_r$}%
   {inverse operator to operator psi r}%

\SetIndexSpace%
\Symb{texRefernceFrame}%Omega/1
   {anholonomity object}%
   {anholonomity object}%

\SetIndexSpace%
\Symb{texLie}%Psi/1
   {basic operator of group Lie}%
   {Lie Basic Operator L}%
\Symb{texLie}%Psi/1
   {}%
   {Lie Basic Operator L}%
\Symb{texLie}%Psi/1
   {basic operator of group Lie}%
   {Lie Basic Operator L, 1-Parameter Group}%
\Symb{texLie}%Psi/1
   {basic operator of group Lie}%
   {Lie Basic Operator R}%
\Symb{texLie}%Psi/1
   {}%
   {Lie Basic Operator R}%
\Symb{texLie}%Psi/1
   {basic operator of group Lie}%
   {Lie Basic Operator R, 1-Parameter Group}%

\SetIndexSpace%
\Symb{texRefernceFrame}%D/2
   {coordinate reference frame}%
   {coordinate reference frame, extensive definition}%
\Symb{texCalculus}%D/2
   {partial derivative of mapping $\Vector f$ with respect to variable $v^{\gi a}$}%
   {partial derivative of mapping, 2, drc vector space}%
\Symb{texCalculus}%D/2
   {partial derivative of mapping $f$ with respect to variable $v^{\gi a}$}%
   {partial derivative of mapping, 2, skew field}%
\Symb{texRefernceFrame}%D/2
   {derivative $e_{(k)}$}%
   {partial(k)}%

\SetIndexSpace%
\Symb{texLie}%Fi/2
   {Lie group composition law}%
   {Lie group composition law}%

\SetIndexSpace%
\Symb{texAffine}%Nabla/2
   {Cartan derivative}%
   {overbrace nabla_l}%
\Symb{texAffine}%Nabla/2
   {derivative}%
   {overline nabla_l, definition 1}%

\SetIndexSpace%
\Symb{texBundleRelation}%Phi/2
   {restriction of correspondence $\Phi$ to set $C$}%
   {restriction of correspondence}%

\SetIndexSpace%
\Symb{texCartesian}%Pi/2
   {Cartesian product of bundles}%
   {Cartesian product of bundles, definition 2}%
\Symb{texCartesian}%Pi/2
   {Cartesian product of total spaces}%
   {Cartesian product of total spaces, definition 2}%
\Symb{texCartesian}%Pi/2
   {reduced Cartesian product of bundles}%
   {reduced Cartesian product of bundles, definition 2}%
\Symb{texCartesian}%Pi/2
   {reduced Cartesian product of total spaces}%
   {reduced Cartesian product of total spaces, definition 2}%

\SetIndexSpace%
\Symb{texBundle}%S/2
   {fibered subset}%
   {fibered subset}%
\Symb{texBundle}%S/2
   {subbundle}%
   {subbundle}%

\CloseIndex

\end{document}